\documentclass[11pt]{article}
\usepackage[top=2cm,bottom=2cm,left=2cm,right=2cm]{geometry}
\usepackage[utf8]{inputenc}
\usepackage[T1]{fontenc}
\usepackage{amssymb,amsfonts,amsmath,amsthm,mathtools}
\usepackage[mathscr]{euscript}
\usepackage{enumitem}
\usepackage[numbers,sort&compress]{natbib}
\usepackage{hyperref}
\usepackage[capitalise,noabbrev]{cleveref}
\usepackage{fancyhdr}
\fancypagestyle{plain}{%
\fancyhf{}
\fancyfoot[C]{\thepage}
}

\numberwithin{equation}{section} 
\setlist[enumerate,1]{label={\rm(\roman*)}} 
\newtheorem{Theorem}{Theorem}[section] 
\newtheorem{lemma}[Theorem]{Lemma} 
\newtheorem{corollary}[Theorem]{Corollary} 
\newtheorem{proposition}[Theorem]{Proposition} 
\theoremstyle{definition} 
\newtheorem{definition}[Theorem]{Definition} 
\newtheorem{example}[Theorem]{Example}

\crefname{Theorem}{Theorem}{Theorems} 
\Crefname{Theorem}{Theorem}{Theorems} 
\crefname{lemma}{Lemma}{Lemmas} 
\Crefname{lemma}{Lemma}{Lemmas} 
\crefname{corollary}{Corollary}{Corollaries} 
\Crefname{corollary}{Corollary}{Corollaries} 
\crefname{proposition}{Proposition}{Propositions} 
\Crefname{proposition}{Proposition}{Propositions} 
\crefname{definition}{Definition}{Definitions} 
\Crefname{definition}{Definition}{Definitions} 
\crefname{example}{Example}{Examples} 
\Crefname{example}{Example}{Examples} 
\crefname{section}{Section}{Sections} 
\Crefname{section}{Section}{Sections} 
\crefname{subsection}{Section}{Sections} 
\Crefname{subsection}{Section}{Sections} 
\newcommand{\R}{\mathbb{R}}
\newcommand{\N}{\mathbb{N}}
\newcommand{\PP}{\mathbb{P}}
\newcommand{\E}{\mathbb{E}}
\newcommand{\Z}{\mathcal{Z}}
\newcommand{\Zo}{\mathcal{Z}^{\circ}}
\DeclareMathOperator{\dist}{dist} 
\DeclareMathOperator{\Var}{Var} 
\DeclareMathOperator{\Cov}{Cov}
\newcommand{\pmin}{\varphi_{-}}
\newcommand{\pmax}{\varphi_{+}}

\begin{document}
\title{\Large Forward and pullback attraction under worst-case and ergodic exponents in a predator--prey model with bounded random refuge}
\author{\normalsize Ever F. Cruzado-Quispe \\  {\footnotesize  Departamento de Matem\'atica, Facultad de Ciencias} \\ {\footnotesize Universidad Nacional del Santa,} \\ {\footnotesize 021809 Nuevo Chimbote, \'Ancash, Per\'u}  \\ \\ \normalsize Javier L\'opez-de-la-Cruz  \\ {\footnotesize Departamento de Matem\'atica Aplicada a las Tecnolog\'ias de la Informaci\'on y las Comunicaciones,} \\ {\footnotesize Universidad Polit\'ecnica de Madrid,} \\ {\footnotesize Campus de Montegancedo, 28660 Boadilla del Monte, Madrid, Spain} \\ \\ \normalsize  Paulo N. Seminario-Huertas\thanks{\raggedright Corresponding author. \\ {\it E-mail addresses}:  2022825003@uns.edu.pe (E.F. Cruzado-Quispe, ORCID: 0000-0001-8045-6785), javier.lopez.delacruz@upm.es (J. L\'opez-de-la-Cruz, ORCID: 0000-0003-2081-5146), paulo.seminario.huertas@upm.es (P.N. Seminario-Huertas, ORCID: 0000-0002-9467-1299) } \\ {\footnotesize Departamento de Matem\'atica Aplicada a las Tecnolog\'ias de la Informaci\'on y las Comunicaciones,} \\ {\footnotesize Universidad Polit\'ecnica de Madrid,} \\ {\footnotesize Campus de Montegancedo, 28660 Boadilla del Monte, Madrid, Spain} }
\date{}
\maketitle
\begin{abstract}
\noindent We consider a predator--prey system with logistic prey, Holling type~II response and a prey refuge fluctuating randomly within a band, driven by a bounded transformation of an Ornstein--Uhlenbeck process. Bounded noise yields a deterministic compact set attracting forward uniformly in the realisation, alongside the random pullback attractor, and we compare the two. The smallest such deterministic set, the minimal uniform attracting set, is the kernel of the control system in which the accessible fraction is any measurable band-valued function, and it meets the coexistence region exactly when the upper worst-case invasion exponent of the predator is positive, whereas the pullback attractor is governed by the ergodic invasion exponent, which decides the predator's fate pathwise. When the predator becomes extinct, the closure of the union of the fibres of the pullback attractor equals the minimal uniform attracting set if and only if the upper worst-case exponent is nonpositive, and otherwise the inclusion is strict and the latter is invariant under no realisation. When the predator persists, the essential closure of the union of the fibres contains, for every frozen environment in which the predator invades, the equilibrium or a periodic orbit of the frozen system, and every state the control system can reach from them. The predator is uniformly persistent when the lower worst-case exponent is positive and almost surely not when it is negative. The ergodic exponent lies below its deterministic counterpart by an amount comparable to the curvature of the response times the variance of the accessible fraction.
\end{abstract}

\medskip
\begin{sloppypar}
\end{sloppypar}

\section{Introduction}

A prey refuge is any mechanism, such as a burrow, a thicket, a depth stratum or a behavioural pattern, that makes part of a prey population unavailable to its predator. Since Maynard Smith \cite{MaynardSmith} observed that, in the Lotka--Volterra model, a refuge protecting a fixed number of prey stabilises an otherwise oscillatory interaction while one protecting a fixed proportion does not, the stabilising role of refuges has been debated in the ecological literature \cite{McNair1986,Sih1987,Krivan1998}. Refuge terms have since become a standard ingredient of predator--prey modelling; see \cite{ChenChenChen2010,GonzalezOlivares2003,Huang2006,Kar2005,Ma2009} and the references therein. The most common formulation splits the prey population $x$ into a protected part $mx$ and an accessible part $(1-m)x$, with $m\in[0,1)$ constant, and lets the predator $y$ act on the accessible part only. Coupled with the logistic growth law of Verhulst \cite{Verhulst} for the prey and the Holling type~II functional response \cite{Holling1959}, this formulation yields
\begin{equation}\label{eq:det}
	\left\{
	\begin{aligned}
		\frac{dx}{dt}&=\alpha x\Bigl(1-\frac{x}{k}\Bigr)
		-\frac{\beta(1-m)\,xy}{1+a(1-m)\,x},\\[4pt]
		\frac{dy}{dt}&=-\gamma y+\frac{c\beta(1-m)\,xy}{1+a(1-m)\,x},
	\end{aligned}
	\right.
\end{equation}
which is the model analysed in \cite{ChenChenChen2010}. Here $\alpha$ is the intrinsic growth rate of the prey, $k$ its carrying capacity, $\beta$ the attack rate, $\beta/a$ the maximal ingestion rate per predator, $1/a$ the half-saturation density, $c$ the conversion efficiency and $\gamma$ the predator mortality; all parameters are positive. For constant $m$ the system \eqref{eq:det} is the classical logistic model with Holling type~II response and no refuge, with the attack rate and the saturation parameter replaced by $\beta(1-m)$ and $a(1-m)$, so that the proportional refuge is a reparametrisation of that model. What changes when $m$ fluctuates is that the effective attack rate and the effective saturation parameter then fluctuate together, through the same scalar, and it is this coupling that drives the quantitative effects of \cref{sec:cost}.

The constancy of $m$ is the weakest point of \eqref{eq:det}. A refuge is a property of the environment, and cover density, turbidity, vegetation, temperature and the behavioural state of the prey all vary in a way that is neither periodic nor predictable, so it is natural to let the accessible fraction $1-m$ be a stochastic process. In the stochastic refuge models known to us (see for instance \cite{QiMeng2021} and the references therein), the refuge level $m$ is kept constant and white noise acts on the growth rates. The difficulty in randomising the refuge itself is that the accessible fraction is a proportion, and any model of the fluctuation must keep it inside $(0,1)$ for every time and every realisation, which rules out additive It\^o perturbations of the coefficient $1-m$ and also multiplicative ones, since a lognormal factor is unbounded above. The device we use, an odd bounded transformation of the Ornstein--Uhlenbeck process, appears in Caraballo, Colucci and Han \cite{CaraballoColucciHan2016} for a predation model and was developed systematically in the work of Caraballo, Colucci, L\'opez-de-la-Cruz and Rapaport \cite{CCLR2019} and in the chemostat models of \cite{CaraballoLopez2020,CGLR2019,LopezGut2026}; see \cite{dOnofrio2013} for a general discussion of bounded noises. Given a mean refuge level $m\in(0,1)$ and an amplitude $\varepsilon\in(0,\min\{m,1-m\})$, the constant $1-m$ in \eqref{eq:det} is replaced by
\[
\psi\bigl(z^{*}(\theta_{t}\omega)\bigr),\qquad
\psi(z):=(1-m)+\frac{2\varepsilon}{\pi}\arctan z ,
\]
where $z^{*}$ is the stationary Ornstein--Uhlenbeck process. The range of $\psi$ is the interval $(1-m-\varepsilon,1-m+\varepsilon)$, hence the accessible fraction is bounded away from $0$ and $1$ by deterministic constants, which are optimal (\cref{prop:sharpbounds}), and the mean is preserved, since $z^{*}$ is centred and $\arctan$ is odd. As $\varepsilon\downarrow0$ every threshold obtained below converges to the threshold of \eqref{eq:det}.

It is useful to keep in mind the deterministic picture. With $p:=1-m$ and $\Phi(\varphi):=c\beta k\varphi/(1+ak\varphi)$, the system \eqref{eq:det} has the trivial equilibrium $(0,0)$, the prey-only equilibrium $(k,0)$ and, when $\Phi(p)>\gamma$, a unique interior equilibrium
\begin{equation}\label{eq:detequil}
	x^{*}=\frac{\gamma}{p\,(c\beta-a\gamma)},
	\qquad
	y^{*}=\frac{\alpha}{\beta p}\Bigl(1-\frac{x^{*}}{k}\Bigr)\bigl(1+ap\,x^{*}\bigr),
\end{equation}
the condition $\Phi(p)>\gamma$ being exactly $x^{*}<k$ and forcing $c\beta>a\gamma$. The number $\lambda_{\rm det}:=\Phi(p)-\gamma$ is the classical invasion criterion of \eqref{eq:det}, and the exponent $\lambda$ introduced below is its exact ergodic counterpart.

In the population models of \cite{CCLR2019,CaraballoLopez2020,CGLR2019,LopezGut2026}, the point of bounded noise is that it yields a deterministic compact set attracting every bounded set forward in time, uniformly in the realisation, which no unbounded perturbation provides. Nor does the random pullback attractor of the white-noise models provide it, since that attractor is reached only by starting further and further in the past. For the system with random refuge both descriptions are available at once, and the aim of this paper is to identify the smallest deterministic set attracting forward uniformly in the realisation, to compare it with the pullback attractor, and to decide which of the two describes the motion along a given realisation. The answer is the following. The smallest such set, the minimal uniform attracting set, which we denote by $\mathcal U$, is the kernel of the control system obtained by letting the accessible fraction be an arbitrary measurable function with values in the band, that is, the set of values at time zero of the bounded complete solutions of that control system (\cref{thm:control}). In light of this identification, $\mathcal U$ depends on the noise only through the band, it contains the global attractor of every autonomous system obtained by freezing the accessible fraction, and it meets the coexistence region exactly when the upper worst-case invasion exponent of the predator, the one obtained by freezing the accessible fraction at the top of the band, is positive (\cref{cor:interior}). The pullback attractor is governed instead by the ergodic invasion exponent, the average against the stationary law of the noise of the invasion exponents of the frozen systems, which decides the fate of the predator along every realisation. Forward attraction, uniform in the realisation, is thus a worst-case notion, because the set it produces retains every environmental history admissible for the band whatever its probability, whereas pullback attraction is an ergodic notion, because the set it produces is governed by the time averages of the noise. This is the sense of the title.

When the predator becomes extinct, the closure of the union of the fibres of the pullback attractor coincides with $\mathcal U$ if and only if the upper worst-case exponent is nonpositive, and otherwise $\mathcal U$ strictly contains it and is invariant under no individual realisation, although it is invariant under the family of all realisations (\cref{thm:charext,cor:invariance}). The essential closure of the union of the fibres, an intrinsic object that does not depend on the null set on which the attractor is modified, is positively invariant under every admissible control for all values of the parameters. When the predator persists it contains, for every frozen environment in which the predator invades, the equilibrium or a periodic orbit of the frozen system, hence every state reachable from them under the control system, so that the two descriptions can differ only at states that the control system cannot reach from these sets (\cref{thm:lowerbound}). The lower bound rests on the independence of the past and the future of the driving noise. To the best of our knowledge, neither this comparison nor the identification of the forward set as the kernel of a control system has been carried out for a population model driven by bounded real noise.

The identification is a representation of a uniform attractor, in the sense of Chepyzhov and Vishik \cite[Chapter IV]{ChepyzhovVishik2002}, by the kernel sections of the hull of the symbols \cite[Theorem IV.5.1]{ChepyzhovVishik2002}, and at the end of \cref{ssec:uniform} we say precisely which symbol space and which topology make it an instance of that theory. Three facts specific to the model make it work. First, the realisations shadow, on any compact time interval and with positive probability, every continuous function with values in the band (\cref{lem:longexc}), and continuous controls are dense among admissible controls in the weak-$*$ topology, so that $\mathcal U$ is as large as the band allows. Second, the solutions depend continuously on the control in that topology once the control is reparametrised affinely (\cref{lem:weakstar}), which is needed because the accessible fraction enters the equations nonlinearly. Third, the absorbing set is uniform over the controls (\cref{cor:absorbing}). The identification places $\mathcal U$ within the theory of control sets of systems driven by bounded noise \cite{Kliemann1987,ColoniusKliemann2000}, where the gap between the support of the stationary measure and the reachable sets is a familiar phenomenon of which \cref{thm:charext,thm:lowerbound} are the counterpart for the pullback attractor. It also places $\mathcal U$ within the theory of attractors of multivalued dynamical systems \cite{MelnikValero1998,CaraballoLangaMelnikValero2003,LambRasmussenRodrigues2015}, of which $\mathcal U$ is an instance (\cref{cor:multivalued}); how such sets change with the amplitude of the noise is the subject of the bifurcation theory of bounded-noise systems \cite{HomburgYoung2006,ZmarrouHomburg2007,BottsHomburgYoung2012}. The distinction between forward and pullback attraction is classical in nonautonomous dynamics; see \cite{CarvalhoLangaRobinson2013,KloedenYang2015,CuiKloedenYang2020} and, for random systems, \cite{CuiLanga2017} and \cite[Section 4.2]{CaraballoHan}. The dissipativity estimates of \cref{sec:noise,sec:wellposed,sec:dissipative}, which rest on the boundedness of the noise and on a deterministic compact absorbing set obtained from a single scalar functional, follow the pattern of \cite{CCLR2019,CaraballoLopez2020,CGLR2019,LopezGut2026}, and what is new begins with the minimal set.

The second part of the paper concerns the thresholds. The ergodic invasion exponent
\[
\lambda=-\gamma+\E\bigl[\Phi(\psi(z^{*}))\bigr]
\]
is the invasion rate of the predator at the prey-only state in the sense of the theory of stochastic persistence \cite{SchreiberBenaimAtchade2011,BenaimSchreiber2019,HeningNguyen2018,Benaim2018,HeningNguyenChesson2021}, a theory formulated for Markov processes with internal and external variables, the Ornstein--Uhlenbeck process being an external variable of the present model. In this setting the prey-only dynamics are unaffected by the noise, so the invasion rate reduces to a one-dimensional Gaussian integral and the dichotomy can be proved pathwise, for every realisation in a set of full measure. The number $\lambda$ is the Lyapunov exponent of the predator at the prey-only state; when $\lambda$ is negative it is the exact rate of extinction of the predator, and when $\lambda$ is positive it forces the pullback attractor to meet the coexistence region and gives persistence in the mean with an explicit constant. Since $\lambda$ depends on the driving process only through its stationary variance, the correlation time of the environment plays no role in the threshold once that variance is fixed. The two worst-case exponents, obtained by freezing the accessible fraction at the endpoints of the band, govern instead the uniform persistence of the predator, which holds with explicit constants when the lower one is positive and fails almost surely when it is negative, whereas the prey is uniformly persistent for all values of the parameters. The gap between the ergodic and the worst-case criteria is of first order in the amplitude of the fluctuation and is a genuine effect of the noise, since in the regime between them the predator persists in the mean while every trajectory comes arbitrarily close to extinction.

The comparison of $\lambda$ with the exponent $\lambda_{\rm det}$ of the deterministic model with the same mean refuge completes the picture. It is an instance of nonlinear averaging, that is, of Jensen's inequality for the concave functional response, a mechanism long recognised in theoretical ecology \cite{RuelAyres1999,MelbourneChesson2006,Denny2017}, and our contribution here is a rigorous two-sided version of it, in terms of the curvature of the response and the variance of the accessible fraction, together with its consequences for the dynamics. The deficit $\lambda_{\rm det}-\lambda$ is strictly increasing in the amplitude and in the stationary variance of the driving process, and its bounds are largest at an intermediate degree of saturation of the response. The set of predator mortalities at which a zero-mean fluctuation drives to extinction a predator that invades the deterministic model is exactly an interval of length equal to the deficit; in that regime no worst-case argument detects the phenomenon, which is unrelated to the switching-induced extinction of \cite{BenaimLobry2016}, as explained in \cref{sec:cost}.

The main features of the paper are the following.
\begin{enumerate}[label={\rm(\alph*)}]
\item The minimal uniform attracting set $\mathcal U$ and the pullback attractor $\mathcal A(\omega)$ are constructed in \cref{thm:uniformattract,thm:attractor}, ordered as in \eqref{eq:sandwich}, and $\mathcal U$ is identified with the kernel of the control system in \cref{thm:control}, with the consequences \cref{cor:invariance,cor:frozen,cor:multivalued,cor:interior}; the identification rests on the shadowing of continuous controls by the noise (\cref{lem:longexc}).
\item The fate of the predator along every realisation is decided by the sign of $\lambda$ in \cref{thm:extinction,thm:attractorext,thm:meanpersist,thm:attractorpersist}, and the two descriptions are compared in \cref{thm:charext} when the predator becomes extinct and in \cref{thm:lowerbound} when it persists.
\item The uniform persistence of the prey, that of the predator when $\lambda_{-}>0$ and its failure when $\lambda_{-}<0$ are proved in \cref{thm:preypersist,thm:uniformpersist,thm:nouniform}, with explicit constants and entrance times (\cref{lem:entrance,cor:coexattractor}).
\item The two-sided bound for the deficit $\lambda_{\rm det}-\lambda$, its monotonicity in the amplitude and in the stationary variance, and the interval of mortalities at which the fluctuation reverses the deterministic outcome are \cref{thm:cost,prop:amax,prop:monotone,prop:kapsig,cor:noiseinduced} and \cref{ssec:reversal}.
\end{enumerate}

The paper is organised as follows. \Cref{sec:noise} fixes the driving noise and proves the optimality of its bounds, the shadowing of continuous controls and the ergodic averages that we shall use. \Cref{sec:wellposed} proves global well-posedness and records the cocycle structure. \Cref{sec:dissipative} is devoted to dissipativity, to the pullback attractor, and to the minimal uniform attracting set and its identification with the kernel of the control system. \Cref{sec:extinction,sec:persistence} contain the extinction and persistence results, including the characterisation of the extinction regime. \Cref{sec:fibres} bounds the union of the fibres from below when the predator persists, \cref{sec:cost} quantifies the cost of the fluctuation, and \cref{sec:open} collects the questions left open.

\section{The random refuge and the invasion exponents}\label{sec:noise}

In this section we fix the driving noise and prove the three facts about it on which the rest of the paper rests, namely the optimality of its deterministic bounds, the positive probability with which it shadows any continuous function with values in the band, and the ergodic averages of nonlinear functions of it.

\subsection{The metric dynamical system and the Ornstein--Uhlenbeck process}

Let $\Omega:=\{\omega\in C(\R;\R)\,:\,\omega(0)=0\}$ be endowed with the compact-open topology, let $\mathscr F$ be its Borel $\sigma$-algebra and $\PP$ the two-sided Wiener measure on $(\Omega,\mathscr F)$, and define the Wiener shift
\begin{equation}\label{eq:shift}
	\theta_{t}\omega(\cdot):=\omega(\cdot+t)-\omega(t),\qquad t\in\R .
\end{equation}
Then $(\Omega,\mathscr F,\PP,\{\theta_{t}\}_{t\in\R})$ is an ergodic metric dynamical system; see \cite[Appendix A]{Arnold} or \cite[Chapter 4]{CaraballoHan}. Fix $\kappa>0$ and $\sigma>0$ and consider the Langevin equation
\begin{equation}\label{eq:langevin}
	dz+\kappa z\,dt=\sigma\,d\omega(t),
\end{equation}
whose unique stationary solution is the Ornstein--Uhlenbeck process \cite{UhlenbeckOrnstein}
\begin{equation}\label{eq:OU}
	z^{*}(\theta_{t}\omega)
	:=-\kappa\sigma\int_{-\infty}^{0}e^{\kappa s}\,(\theta_{t}\omega)(s)\,ds ,
	\qquad t\in\R .
\end{equation}

\begin{proposition}\label{prop:OU}
There exists a $\theta$-invariant set $\Omega_{1}\in\mathscr F$ with $\PP(\Omega_{1})=1$ such that, for every $\omega\in\Omega_{1}$, the integral in \eqref{eq:OU} converges and the map $t\mapsto z^{*}(\theta_{t}\omega)$ is continuous on $\R$ and solves \eqref{eq:langevin} pathwise. Moreover, the random variable $\omega\mapsto z^{*}(\omega)$ is Gaussian and centred, with
	\begin{equation}\label{eq:OUvar}
		\Var\bigl(z^{*}\bigr)=\sigma_{*}^{2}:=\frac{\sigma^{2}}{2\kappa},
	\end{equation}
and the process $t\mapsto z^{*}(\theta_{t}\omega)$ is stationary and ergodic.
\end{proposition}

\begin{proof}
Let $\Omega_{1}$ be the $\theta$-invariant set of full measure constructed in \cite[Proposition 4.1]{CaraballoHan}; see also \cite[Lemma 4.1]{CaraballoKloedenSchmalfuss}. On it $\omega(s)/s\to0$ as $|s|\to\infty$, the integral in \eqref{eq:OU} converges, and the pathwise assertions of the statement hold. Off $\Omega_{1}$ we set $z^{*}(\omega):=0$, which makes $z^{*}$ a random variable on $\Omega$ without affecting any statement made up to null sets. For each $n\in\N$, the integration by parts formula for the Wiener integral of the $C^{1}$ function $s\mapsto e^{\kappa s}$ on $[-n,0]$ holds almost surely and gives
	\[
	\sigma\int_{-n}^{0}e^{\kappa s}\,d\omega(s)=-\sigma e^{-\kappa n}\omega(-n) -\kappa\sigma\int_{-n}^{0}e^{\kappa s}\omega(s)\,ds .
	\]
As $n\to\infty$ the right-hand side converges to $z^{*}(\omega)$ on $\Omega_{1}$, because $\omega(s)/s\to0$, so that $z^{*}$ is measurable on $\Omega_{1}$ as a pointwise limit of continuous maps; the left-hand side converges in $L^{2}(\PP)$ to the Wiener integral $\sigma\int_{-\infty}^{0}e^{\kappa s}\,d\omega(s)$, so that the two limits coincide almost surely and $z^{*}$ is centred Gaussian with variance $\sigma^{2}\int_{-\infty}^{0}e^{2\kappa s}ds=\sigma^{2}/(2\kappa)$ by the It\^o isometry. Stationarity follows from the invariance of $\PP$ under $\theta$, and the process is ergodic because a factor of an ergodic flow is ergodic.
\end{proof}

\subsection{The bounded transformation and the standing hypotheses}

Throughout the paper we assume the following.

\medskip
\noindent\textbf{(H1)} $\alpha,k,\beta,a,c,\gamma>0$.

\noindent\textbf{(H2)} $m\in(0,1)$ and $\varepsilon\in\bigl(0,\min\{m,1-m\}\bigr)$.

\noindent\textbf{(H3)} $\kappa>0$ and $\sigma>0$.
\medskip

\noindent Write
\[
p:=1-m\in(0,1),\qquad \pmin:=p-\varepsilon,\qquad \pmax:=p+\varepsilon,
\]
so that \textbf{(H2)} is exactly $0<\pmin<p<\pmax<1$, and define $\psi:\R\to\R$ by
\begin{equation}\label{eq:psi}
	\psi(z):=p+\frac{2\varepsilon}{\pi}\arctan z .
\end{equation}
For $\omega\in\Omega_{1}$ the accessible fraction is the process
\begin{equation}\label{eq:phi}
	\varphi_{\omega}(t):=\psi\bigl(z^{*}(\theta_{t}\omega)\bigr),\qquad t\in\R .
\end{equation}
The map $\psi$ is a $C^{\infty}$, strictly increasing bijection from $\R$ onto the open interval $(\pmin,\pmax)$, and $\psi-p$ is odd. Consequently,
\begin{equation}\label{eq:phibounds}
	\pmin<\varphi_{\omega}(t)<\pmax\qquad\text{for every }t\in\R
	\text{ and every }\omega\in\Omega_{1}.
\end{equation}
Any $C^{1}$, strictly increasing, odd and bounded function could replace $\arctan$ in \eqref{eq:psi} after rescaling its range to $(\pmin,\pmax)$. The next proposition says that the constants in \eqref{eq:phibounds} cannot be improved; as we note after its proof, it also shows that a map $\psi$ whose range is all of $(0,1)$ admits no deterministic worst-case constants at all.

\begin{proposition}[Optimality of the bounds]\label{prop:sharpbounds}
There exists a $\theta$-invariant set $\Omega_{2}\subset\Omega_{1}$ with $\PP(\Omega_{2})=1$ such that, for every $\omega\in\Omega_{2}$,
	\begin{equation}\label{eq:zunbounded}
		\limsup_{t\to+\infty}z^{*}(\theta_{t}\omega)=+\infty,
		\qquad
		\liminf_{t\to+\infty}z^{*}(\theta_{t}\omega)=-\infty ,
	\end{equation}
and consequently
	\begin{equation}\label{eq:phisharp}
		\sup_{t\ge0}\varphi_{\omega}(t)=\pmax,
		\qquad
		\inf_{t\ge0}\varphi_{\omega}(t)=\pmin ,
	\end{equation}
neither value being attained.
\end{proposition}

\begin{proof}
Fix $M>0$ and set $B_{M}:=\{\omega\in\Omega_{1}\,:\,\limsup_{t\to+\infty}z^{*}(\theta_{t}\omega)\ge M\}$, a measurable set because $t\mapsto z^{*}(\theta_{t}\omega)$ is continuous and the $\limsup$ may be taken along rational times. For every $s\in\R$,
	\[
	\limsup_{t\to+\infty}z^{*}(\theta_{t}\theta_{s}\omega)
	=\limsup_{t\to+\infty}z^{*}(\theta_{t+s}\omega)
	=\limsup_{t\to+\infty}z^{*}(\theta_{t}\omega),
	\]
so that $\theta_{s}^{-1}B_{M}=B_{M}$ for every $s$. Hence $B_{M}$ is $\theta$-invariant and, by ergodicity, its probability is either $0$ or $1$. To see that it is positive, consider the events $E_{n}:=\{\omega\in\Omega_{1}\,:\,z^{*}(\theta_{n}\omega)>M\}$, $n\in\N$. By stationarity all of them have the same probability $q_{M}:=\PP(z^{*}>M)=\PP(\mathcal N(0,\sigma_{*}^{2})>M)$, which is positive because $\sigma_{*}^{2}>0$. The reverse Fatou lemma applied to the indicator functions gives
	\[
	\PP\Bigl(\limsup_{n\to\infty}E_{n}\Bigr)\ \ge\
	\limsup_{n\to\infty}\PP(E_{n})=q_{M}>0 ,
	\]
and since $\limsup_{n}E_{n}\subset B_{M}$ we obtain $\PP(B_{M})=1$. The intersection over $M\in\N$ gives the first identity in \eqref{eq:zunbounded} on a set of full measure, and the second follows in the same way from $\PP(z^{*}<-M)=q_{M}$. Let $\Omega_{2}$ be the intersection of these two sets with $\Omega_{1}$. It is $\theta$-invariant because each $B_{M}$ is. Finally, \eqref{eq:phisharp} follows from \eqref{eq:zunbounded} because $\psi$ is continuous and strictly increasing with $\psi(z)\to\pmax$ as $z\to+\infty$ and $\psi(z)\to\pmin$ as $z\to-\infty$, and the values are not attained because $\psi(\R)=(\pmin,\pmax)$.
\end{proof}

\Cref{prop:sharpbounds} applies verbatim to any strictly increasing $\psi$, in particular to $\psi_{0}(z)=\frac12+\frac1\pi\arctan z$, for which it gives $\inf_{t}\psi_{0}(z^{*}(\theta_{t}\omega))=0$ and $\sup_{t}\psi_{0}(z^{*}(\theta_{t}\omega))=1$ almost surely, and no worst-case constants exist for $\psi_{0}$. Boundedness of the range is not enough. What worst-case arguments require is that the range be bounded away from the endpoints, and this is exactly what the constraint on $\varepsilon$ in \textbf{(H2)} provides. The next lemma says considerably more than \cref{prop:sharpbounds}. For every continuous function $\varphi_{0}$ with values in the closed band, and for an arbitrarily long horizon, the accessible fraction remains uniformly close to $\varphi_{0}$ on an event of positive probability. It is a one-dimensional instance of the support theorem for diffusions, and only the law of $z^{*}(\omega)$, the support of the Wiener measure and the independence of its increments enter its proof. This is the only quantitative property of the noise used in \cref{ssec:uniform}.

\begin{lemma}[Shadowing of continuous controls]\label{lem:longexc}
Let $T>0$, let $\varphi_{0}:[0,T]\to[\pmin,\pmax]$ be continuous and let $\varrho>0$. Then
	\begin{equation}\label{eq:longexc}
		\PP\Bigl(\Bigl\{\omega\in\Omega_{1}\,:\ \sup_{t\in[0,T]}
		\bigl|\varphi_{\omega}(t)-\varphi_{0}(t)\bigr|<\varrho\Bigr\}\Bigr)\ >\ 0 .
	\end{equation}
In particular, every set of full measure contains some $\omega$ with $\sup_{t\in[0,T]}|\varphi_{\omega}(t)-\varphi_{0}(t)|<\varrho$.
\end{lemma}

\begin{proof}
The event in \eqref{eq:longexc} is measurable, since on $\Omega_{1}$ the map $t\mapsto\varphi_{\omega}(t)$ is continuous and the supremum may be taken over rational $t$. We first record a pathwise variation-of-constants formula. Using the identities $(\theta_{t}\omega)(s)=\omega(t+s)-\omega(t)$ and $\kappa\int_{-\infty}^{0}e^{\kappa s}ds=1$, together with the substitution $u=t+s$ in \eqref{eq:OU},
	\[
	z^{*}(\theta_{t}\omega)
	=\sigma\omega(t)-\kappa\sigma e^{-\kappa t}\!\int_{-\infty}^{t}\!\!e^{\kappa u}\omega(u)\,du .
	\]
Since $\omega(0)=0$, the same computation at $t=0$ gives $z^{*}(\omega)=-\kappa\sigma\int_{-\infty}^{0}e^{\kappa u}\omega(u)\,du$, whence, splitting the integral at $0$,
	\begin{equation}\label{eq:varconst}
		z^{*}(\theta_{t}\omega)=e^{-\kappa t}z^{*}(\omega)
		+\sigma\Bigl(\omega(t)-\kappa\!\int_{0}^{t}\!e^{-\kappa(t-u)}\omega(u)\,du\Bigr),
		\qquad t\ge0 .
	\end{equation}

Replacing $\varphi_{0}$ by $p+(1-\delta)(\varphi_{0}-p)$ with $\delta\in(0,1)$ and $\delta\varepsilon<\varrho/2$, and $\varrho$ by $\varrho/2$, we may assume that $\varphi_{0}$ takes values in the open interval $(\pmin,\pmax)$. Then $M_{0}:=\psi^{-1}\circ\varphi_{0}$ is continuous on $[0,T]$, and by uniform continuity of $\psi$ there is $\varrho'>0$ such that $|z-M_{0}(t)|<\varrho'$ implies $|\psi(z)-\varphi_{0}(t)|<\varrho$ for all $t\in[0,T]$. Put
	\[
	r(t):=M_{0}(t)-e^{-\kappa t}M_{0}(0),
	\qquad
	h(t):=\frac1\sigma\Bigl(r(t)+\kappa\!\int_{0}^{t}r(s)\,ds\Bigr),
	\qquad t\in[0,T],
	\]
so that $h$ is continuous, $h(0)=0$, and a direct differentiation of $e^{-\kappa t}\int_{0}^{t}e^{\kappa u}h(u)\,du$ shows that
	\begin{equation}\label{eq:volterra}
		\sigma\Bigl(h(t)-\kappa\!\int_{0}^{t}\!e^{-\kappa(t-u)}h(u)\,du\Bigr)=r(t),
		\qquad t\in[0,T].
	\end{equation}
Set $\delta_{0}:=\varrho'/(4\sigma)$ and consider the events
	\begin{align*}
	A&:=\Bigl\{\omega\in\Omega\,:\ \int_{-\infty}^{0}e^{\kappa u}\omega(u)\,du\ \text{converges and}\
	\Bigl|-\kappa\sigma\!\int_{-\infty}^{0}\!e^{\kappa u}\omega(u)\,du-M_{0}(0)\Bigr|<\tfrac{\varrho'}{2}\Bigr\},\\
	B&:=\Bigl\{\omega\in\Omega\,:\ \sup_{t\in[0,T]}\bigl|\omega(t)-h(t)\bigr|\le\delta_{0}\Bigr\}.
	\end{align*}
The event $A$ is measurable with respect to the $\sigma$-algebra generated by $\omega|_{(-\infty,0]}$, and on $\Omega_{1}$ it coincides with $\{|z^{*}(\omega)-M_{0}(0)|<\varrho'/2\}$. By \eqref{eq:OUvar} the random variable $z^{*}$ is Gaussian with positive variance, so its law charges every nonempty open interval and $\PP(A)>0$. The event $B$ is measurable with respect to the $\sigma$-algebra generated by the increments of $\omega$ after time $0$, and $\PP(B)>0$. Indeed, the topological support of the Wiener measure on $C([0,T])$ is the closure of its Cameron--Martin space \cite[Theorem 3.6.1]{Bogachev1998}, and the latter, which consists of the absolutely continuous functions vanishing at $0$ with square-integrable derivative \cite[Lemma 2.3.14]{Bogachev1998}, is dense in the space of continuous functions vanishing at $0$, to which $h$ belongs. Under the Wiener measure the two $\sigma$-algebras are independent, so $\PP(A\cap B)=\PP(A)\PP(B)>0$.

Let $\omega\in A\cap B\cap\Omega_{1}$, a set of positive probability since $\PP(\Omega_{1})=1$, and write $e(t):=\omega(t)-h(t)$, so that $\sup_{[0,T]}|e|\le\delta_{0}$. Substituting $\omega=h+e$ into \eqref{eq:varconst} and using \eqref{eq:volterra} we obtain, for $t\in[0,T]$,
	\[
	z^{*}(\theta_{t}\omega)-M_{0}(t)
	=e^{-\kappa t}\bigl(z^{*}(\omega)-M_{0}(0)\bigr)
	+\sigma\Bigl(e(t)-\kappa\!\int_{0}^{t}\!e^{-\kappa(t-u)}e(u)\,du\Bigr).
	\]
The first term is bounded in absolute value by $\varrho'/2$, and the second by $\sigma\delta_{0}(2-e^{-\kappa t})\le2\sigma\delta_{0}=\varrho'/2$. Hence $|z^{*}(\theta_{t}\omega)-M_{0}(t)|<\varrho'$ for every $t\in[0,T]$ and, by the choice of $\varrho'$, $|\varphi_{\omega}(t)-\varphi_{0}(t)|<\varrho$ there. This proves \eqref{eq:longexc}, and the last assertion holds because an event of positive probability meets every set of full measure.
\end{proof}

\subsection{Ergodic averages}

The following consequence of Birkhoff's theorem is the analytic engine of \cref{sec:extinction,sec:cost}. Observe that the linear statement $\frac1t\int_{0}^{t}z^{*}(\theta_{s}\omega)\,ds\to0$, which is part of \cite[Proposition 4.1]{CaraballoHan}, is not sufficient for our purposes, because the quantities to be averaged are nonlinear functions of the process.

\begin{lemma}[Ergodic averages]\label{lem:ergodic}
Let $g:[\pmin,\pmax]\to\R$ be continuous. Then there is a $\theta$-invariant set $\Omega_{g}\subset\Omega_{2}$ of full measure such that, for every $\omega\in\Omega_{g}$,
	\begin{equation}\label{eq:birkhoff}
		\lim_{T\to+\infty}\frac1T\int_{0}^{T}g\bigl(\varphi_{\omega}(t)\bigr)\,dt
		=\lim_{T\to+\infty}\frac1T\int_{-T}^{0}g\bigl(\varphi_{\omega}(t)\bigr)\,dt
		=\E\bigl[g(\psi(z^{*}))\bigr].
	\end{equation}
In particular, taking $g=\mathrm{id}$,
	\begin{equation}\label{eq:meanpreserved}
		\lim_{T\to+\infty}\frac1T\int_{0}^{T}\varphi_{\omega}(t)\,dt
		=\E\bigl[\psi(z^{*})\bigr]=p=1-m .
	\end{equation}
\end{lemma}

\begin{proof}
Set $G(\omega):=g(\psi(z^{*}(\omega)))$. Since $g$ is continuous on a compact interval and $\psi(\R)\subset[\pmin,\pmax]$, $G$ is measurable and bounded, hence integrable, and $G(\theta_{t}\omega)=g(\varphi_{\omega}(t))$ by \eqref{eq:phi}, and $(t,\omega)\mapsto G(\theta_{t}\omega)$ is jointly measurable because $\theta$ is a measurable flow. The Birkhoff--Khinchin ergodic theorem for the ergodic measure-preserving flow $\{\theta_{t}\}$, in the form \cite[Chapter 1, Section 2, Theorem 1]{CornfeldFominSinai1982}, gives both equalities in \eqref{eq:birkhoff} on a set of full measure, the averages over $[0,T]$ and over $[-T,0]$ having the same almost sure limit, which equals $\E[G]$ by ergodicity. The set on which \eqref{eq:birkhoff} holds is $\theta$-invariant, because replacing $\omega$ by $\theta_{s}\omega$ shifts the interval of integration by $s$ and changes the Ces\`aro means by at most $2|s|\sup|g|/T$. We intersect it with $\Omega_{2}$. For \eqref{eq:meanpreserved}, the law of $z^{*}$ is symmetric and $\arctan$ is odd and bounded, so $\E[\arctan z^{*}]=0$ and $\E[\psi(z^{*})]=p$.
\end{proof}

\subsection{The model}

Let $\omega\in\Omega_{2}$ be fixed. The system studied in this paper is
\begin{equation}\label{eq:sys}
	\left\{
	\begin{aligned}
		\frac{dx}{dt}&=\alpha x\Bigl(1-\frac{x}{k}\Bigr)
		-\frac{\beta\,\varphi_{\omega}(t)\,xy}{1+a\,\varphi_{\omega}(t)\,x},\\[4pt]
		\frac{dy}{dt}&=-\gamma y
		+\frac{c\beta\,\varphi_{\omega}(t)\,xy}{1+a\,\varphi_{\omega}(t)\,x},
	\end{aligned}
	\right.
\end{equation}
posed on the closed positive quadrant $\Z$, whose interior $\Zo$ is the coexistence region,
\[
\Z:=\bigl\{(x,y)\in\R^{2}\,:\,x\ge0,\ y\ge0\bigr\},
\qquad
\Zo:=\bigl\{(x,y)\in\R^{2}\,:\,x>0,\ y>0\bigr\}.
\]
For $x\ge0$ and $\varphi\in[\pmin,\pmax]$ we introduce the functional response and two normalisations of it,
\begin{equation}\label{eq:GF}
	G(x,\varphi):=\frac{\beta\varphi x}{1+a\varphi x},
	\qquad
	F(x,\varphi):=c\,G(x,\varphi),
	\qquad
	\Phi(\varphi):=F(k,\varphi)=\frac{c\beta k\varphi}{1+ak\varphi}.
\end{equation}

\begin{lemma}\label{lem:Fprops}
On $[0,\infty)\times[\pmin,\pmax]$ the function $F$ of \eqref{eq:GF} is of class $C^{\infty}$ and satisfies the following properties.
\begin{enumerate}
\item $\partial_{x}F=\dfrac{c\beta\varphi}{(1+a\varphi x)^{2}}>0$ and
$\partial_{\varphi}F=\dfrac{c\beta x}{(1+a\varphi x)^{2}}\ge0$, with equality only at $x=0$, so $F$ is increasing in each variable.
\item $0\le F(x,\varphi)<\dfrac{c\beta}{a}$ and
$F(x,\varphi)\le c\beta\varphi x\le c\beta\pmax x$.
\item $\Phi$ is strictly increasing and strictly concave on $[\pmin,\pmax]$, with
		\begin{equation}\label{eq:Phider}
			\Phi'(\varphi)=\frac{c\beta k}{(1+ak\varphi)^{2}}>0,
			\qquad
			\Phi''(\varphi)=-\frac{2ac\beta k^{2}}{(1+ak\varphi)^{3}}<0 .
		\end{equation}
\end{enumerate}
\end{lemma}

\begin{proof}
Everything follows by direct differentiation, since $1+a\varphi x\ge1$ on the stated domain. For (ii) we write $F(x,\varphi)=\frac{c\beta}{a}\cdot\frac{a\varphi x}{1+a\varphi x}$ and observe that the second factor lies in $[0,1)$.
\end{proof}

\begin{definition}[Invasion exponents]\label{def:lambda}
The \emph{ergodic invasion exponent of the predator at the prey-only state} is
	\begin{equation}\label{eq:lambda}
		\lambda:=-\gamma+\E\bigl[\Phi(\psi(z^{*}))\bigr],
	\end{equation}
and the \emph{lower} and \emph{upper worst-case exponents} are
	\begin{equation}\label{eq:lambdapm}
		\lambda_{-}:=-\gamma+\Phi(\pmin),
		\qquad
		\lambda_{+}:=-\gamma+\Phi(\pmax).
	\end{equation}
We also write $\lambda_{\rm det}:=-\gamma+\Phi(p)$ for the exponent of the deterministic model \eqref{eq:det}.
\end{definition}

Since $\psi(z^{*})$ is a nondegenerate random variable with values in $(\pmin,\pmax)$ and $\Phi$ is strictly increasing,
\begin{equation}\label{eq:lambdaorder}
	\lambda_{-}<\lambda<\lambda_{+}.
\end{equation}
The number $\lambda$ is a one-dimensional Gaussian integral,
\begin{equation}\label{eq:lambdaintegral}
	\lambda=-\gamma+\frac{1}{\sqrt{2\pi}\,\sigma_{*}}
	\int_{\R}\frac{c\beta k\,\psi(z)}{1+ak\,\psi(z)}\,
	e^{-z^{2}/(2\sigma_{*}^{2})}\,dz ,
	\qquad \sigma_{*}^{2}=\frac{\sigma^{2}}{2\kappa},
\end{equation}
and it depends on the Ornstein--Uhlenbeck process only through its stationary variance. In the parametrisation by $(\kappa,\sigma)$ this means that the correlation time $1/\kappa$ plays no role at fixed $\sigma_{*}$, whereas at fixed $\sigma$ a longer correlation time increases $\sigma_{*}$ and, by \cref{prop:kapsig} below, lowers $\lambda$. We shall also need the variance of the accessible fraction,
\begin{equation}\label{eq:Vdef}
	V:=\Var\bigl(\psi(z^{*})\bigr)=\E\bigl[(\psi(z^{*})-p)^{2}\bigr]
	=\frac{4\varepsilon^{2}}{\pi^{2}}\,\E\bigl[\arctan^{2}z^{*}\bigr]
	\ \in\ \bigl(0,\varepsilon^{2}\bigr),
\end{equation}
where the upper bound is strict because $\arctan^{2}<\pi^{2}/4$ pointwise, and the lower one because $z^{*}$ is nondegenerate.

\subsection{The four criteria}\label{ssec:criteria}

All the results of \cref{sec:extinction,sec:persistence} are governed by the sign of one of the exponents of \cref{def:lambda}. We display the four conditions once and for all,
\begin{align}
	&\lambda<0 , \tag{E}\label{eq:E}\\[2pt]
	&\lambda_{+}<0,\quad\text{that is}\quad
	\gamma>\Phi(\pmax)=\frac{c\beta k(1-m+\varepsilon)}{1+ak(1-m+\varepsilon)} ,
	\tag{E$'$}\label{eq:Eprime}\\[2pt]
	&\lambda>0 , \tag{P}\label{eq:P}\\[2pt]
	&\lambda_{-}>0,\quad\text{that is}\quad
	\gamma<\Phi(\pmin)=\frac{c\beta k(1-m-\varepsilon)}{1+ak(1-m-\varepsilon)} .
	\tag{P$'$}\label{eq:Pprime}
\end{align}
Conditions \eqref{eq:E} and \eqref{eq:P} are the \emph{ergodic} criteria. Together they cover all values of the parameters except the critical case $\lambda=0$. Conditions \eqref{eq:Eprime} and \eqref{eq:Pprime} are their \emph{worst-case} counterparts, strictly stronger but in closed form.

\begin{lemma}[Ordering of the criteria]\label{lem:ordering}
Assume \textbf{(H1)}--\textbf{(H3)}. Then
	\begin{equation}\label{eq:hierarchy}
		\gamma>\frac{c\beta}{a}\ \Longrightarrow\ \eqref{eq:Eprime}
		\ \Longrightarrow\ \eqref{eq:E},
		\qquad\text{and}\qquad
		\eqref{eq:Pprime}\ \Longrightarrow\ \eqref{eq:P},
	\end{equation}
and none of the three implications can be reversed.
\end{lemma}

\begin{proof}
By \cref{lem:Fprops}(ii), $0<\Phi(\pmax)<c\beta/a$, so $\gamma>c\beta/a$ gives $\gamma>\Phi(\pmax)$, and the implication cannot be reversed because $\Phi(\pmax)<c\beta/a$. The remaining two implications are \eqref{eq:lambdaorder}, which is strict.
\end{proof}

The following parameter set, at which \eqref{eq:E} holds and \eqref{eq:Eprime} fails, will serve as a running illustration.

\begin{example}\label{ex:nie}
Take
	\[
	\alpha=1,\quad k=40,\quad \beta=0.9,\quad a=0.9,\quad c=0.6,\quad
	m=0.5,\quad \varepsilon=0.45,\quad \kappa=1,\quad \sigma=2 ,
	\]
so that $p=0.5$, $[\pmin,\pmax]=[0.05,0.95]$ and $\sigma_{*}^{2}=2$. Then $\Phi(p)=0.56842$ and $\Phi(\pmax)=0.58295$ to five decimal places, and numerical quadrature of the smooth bounded integrands of \eqref{eq:lambdaintegral} and \eqref{eq:Vdef} gives, to the same precision,
	\[
	\E\bigl[\Phi(\psi(z^{*}))\bigr]=0.55963,\qquad V=0.05260 .
	\]
Choosing $\gamma=0.5679522$ therefore yields $\lambda_{\rm det}=+4.69\times10^{-4}$, $\lambda=-8.32\times10^{-3}$ and $\lambda_{+}=+1.50\times10^{-2}$. Hence \eqref{eq:E} holds while \eqref{eq:Eprime} fails, and moreover $\lambda_{\rm det}>0$, that is, the deterministic model \eqref{eq:det} with the same mean refuge possesses the interior equilibrium \eqref{eq:detequil}. These parameters fall within the regimes of \cref{thm:Ustrict,cor:noiseinduced}. The value of $\gamma$ is chosen so that $\lambda_{\rm det}$ lies inside the interval \eqref{eq:nie} required by \cref{cor:noiseinduced}, whose right endpoint is $9.38\times10^{-4}$ at these parameters; any $\gamma$ in that interval serves.
\end{example}

Finally, we fix the set of realisations on which every pathwise statement of this paper is made. For $\eta\ge0$ put
\begin{equation}\label{eq:Phieta}
	\Phi_{\eta}(\varphi):=F(k+\eta,\varphi)
	=\frac{c\beta(k+\eta)\varphi}{1+a(k+\eta)\varphi},
	\qquad
	\lambda_{\eta}:=-\gamma+\E\bigl[\Phi_{\eta}(\psi(z^{*}))\bigr],
\end{equation}
so that $\Phi_{0}=\Phi$ and $\lambda_{0}=\lambda$, and set
\begin{equation}\label{eq:Omega0}
	\Omega_{0}:=\Omega_{2}\cap\bigcap_{\eta\in\mathbb Q\cap[0,\infty)}
	\Omega_{\Phi_{\eta}},
\end{equation}
where $\Omega_{g}$ is the set provided by \cref{lem:ergodic}. Being a countable intersection of $\theta$-invariant sets of full measure, $\Omega_{0}$ is $\theta$-invariant and $\PP(\Omega_{0})=1$. From now on all pathwise statements are made for $\omega\in\Omega_{0}$, which inherits every property established for $\Omega_{1}$ and $\Omega_{2}$. We shall see after \cref{thm:control} that the deterministic objects constructed below do not depend on this choice.

\begin{lemma}\label{lem:lambdaeta}
The map $\eta\mapsto\lambda_{\eta}$ of \eqref{eq:Phieta} is continuous and strictly increasing on $[0,\infty)$, and $\lambda_{\eta}\downarrow\lambda$ as $\eta\downarrow0$.
\end{lemma}

\begin{proof}
For fixed $\varphi\in[\pmin,\pmax]$ the map $\eta\mapsto F(k+\eta,\varphi)$ is continuous and strictly increasing by \cref{lem:Fprops}(i), and it is bounded by $c\beta/a$ uniformly in $\eta$ and $\varphi$. Dominated convergence gives the continuity, and the strict monotonicity passes to the expectation because the difference of the integrands is positive.
\end{proof}

We close this section by recording the sense in which $\lambda$ is a Lyapunov exponent. The point $(k,0)$ is a common zero of the vector fields $\mathcal F(\cdot,\varphi)$ of \eqref{eq:field} below, $\varphi\in[\pmin,\pmax]$, and the linearisation of the predator equation along the prey-only solution $x\equiv k$ is the scalar equation
\begin{equation}\label{eq:linearised}
	\dot\zeta=\bigl[-\gamma+\Phi\bigl(\varphi_{\omega}(t)\bigr)\bigr]\zeta ,
	\qquad\text{so that}\qquad
	\lim_{t\to+\infty}\frac1t\log|\zeta(t)|=\lambda
	\quad\text{for every }\omega\in\Omega_{0}\text{ and }\zeta(0)\ne0 ,
\end{equation}
by \cref{lem:ergodic} with $g=\Phi$, which applies because $\Omega_{0}\subset\Omega_{\Phi_{0}}=\Omega_{\Phi}$ since $0\in\mathbb Q$. Thus $\lambda$ is, for every realisation in $\Omega_{0}$ and whatever its sign, the Lyapunov exponent of the predator at the prey-only state. \Cref{thm:extinction} below shows that when $\lambda<0$ this linear rate is also the exact rate of decay of the nonlinear predator, and \cref{thm:meanpersist,thm:attractorpersist} show what survives of the linear instability when $\lambda>0$. The role of the Lyapunov exponent of the linearisation at a common zero in deciding the behaviour of a randomly switched system is studied in \cite{BenaimStrickler2019}.

\section{Well-posedness and the random dynamical system}\label{sec:wellposed}

This section is standard. Its purpose is to fix notation and to record the representation \eqref{eq:representation} and the a priori bound \eqref{eq:apriori}, which are used repeatedly later on. We write \eqref{eq:sys} as $\dot v=\mathcal F(v,\varphi_{\omega}(t))$, where $v=(x,y)$ and, for $\varphi\in[\pmin,\pmax]$,
\begin{equation}\label{eq:field}
	\mathcal F(v,\varphi):=\Bigl(\alpha x\Bigl(1-\frac{x}{k}\Bigr)-G(x,\varphi)\,y,\
	-\gamma y+F(x,\varphi)\,y\Bigr).
\end{equation}
The function $\mathcal F$ is of class $C^{\infty}$ on $\mathcal V\times[\pmin,\pmax]$, where $\mathcal V:=\{(x,y)\in\R^{2}\,:\,x>-1/(2a\pmax)\}\supset\Z$ is open, because $1+a\varphi x>1-\varphi/(2\pmax)\ge\frac12$ there.

\begin{Theorem}[Global well-posedness]\label{thm:wellposed}
Assume \textbf{(H1)}--\textbf{(H3)} and let $\omega\in\Omega_{0}$. For every $v_{0}=(x_{0},y_{0})\in\Z$ the system \eqref{eq:sys} with initial condition $v(0)=v_{0}$ has a unique solution
	\[
	v(\cdot\,;0,\omega,v_{0})=\bigl(x(\cdot\,;0,\omega,v_{0}),
	y(\cdot\,;0,\omega,v_{0})\bigr)\in C^{1}\bigl([0,+\infty);\R^{2}\bigr),
	\]
which satisfies $v(t;0,\omega,v_{0})\in\Z$ for every $t\ge0$ together with
	\begin{equation}\label{eq:representation}
		\begin{aligned}
			x(t)&=x_{0}\exp\!\left(\int_{0}^{t}\Bigl[\alpha\Bigl(1-\frac{x(s)}{k}\Bigr)
			-\frac{\beta\varphi_{\omega}(s)\,y(s)}{1+a\varphi_{\omega}(s)\,x(s)}
			\Bigr]ds\right),\\[4pt]
			y(t)&=y_{0}\exp\!\left(\int_{0}^{t}\bigl[-\gamma
			+F\bigl(x(s),\varphi_{\omega}(s)\bigr)\bigr]ds\right).
		\end{aligned}
	\end{equation}
Consequently, $\Z$, its interior $\Zo$ and each of its two faces $\{x=0\}$ and $\{y=0\}$ are positively invariant, and
	\begin{equation}\label{eq:apriori}
		x(t)\le\max\{x_{0},k\}\quad\text{and}\quad
		y(t)\le y_{0}\,e^{(c\beta/a)t}\qquad\text{for all }t\ge0 .
	\end{equation}
The same statements hold, with the same proof, if $\varphi_{\omega}$ is replaced by any measurable function $\varphi:[0,\infty)\to[\pmin,\pmax]$, the solution being then absolutely continuous and the equation holding almost everywhere.
\end{Theorem}

\begin{proof}
Since $t\mapsto\varphi_{\omega}(t)$ is continuous and $\mathcal F$ is $C^{1}$ in $v$ on the open set $\mathcal V$, the Picard--Lindel\"of theorem yields a unique maximal solution $v\in C^{1}([0,T_{\max});\mathcal V)$. On $[0,T_{\max})$ both components of the field factor as $x\,g_{1}(t)$ and $y\,g_{2}(t)$ with $g_{1},g_{2}$ continuous, so $x$ and $y$ solve linear scalar equations, which gives \eqref{eq:representation}, and since the exponentials are positive, $x(t)$ and $y(t)$ have the signs of $x_{0}$ and $y_{0}$. This proves the invariance of $\Z$, $\Zo$ and the two faces. For $v_{0}\in\Z$ the predation term in the first equation is nonnegative, hence
	\begin{equation}\label{eq:logisticcomp}
		\frac{dx}{dt}\ \le\ \alpha x\Bigl(1-\frac{x}{k}\Bigr),\qquad t\in[0,T_{\max}),
	\end{equation}
and the comparison theorem for scalar equations gives $x(t)\le u(t)$, where $u$ is the logistic solution starting from $x_{0}$, which is nondecreasing and bounded by $k$ if $x_{0}\le k$ and nonincreasing if $x_{0}>k$. This is the first bound in \eqref{eq:apriori}. The second follows from $\dot y\le(c\beta/a)\,y$, which is \cref{lem:Fprops}(ii), and Gronwall's inequality. Both bounds are finite on bounded intervals, so $T_{\max}=+\infty$. For a measurable $\varphi$ the field is a Carath\'eodory function, locally Lipschitz in $v$ uniformly in $t$, and the same argument applies to the Carath\'eodory solution.
\end{proof}

The bound $x(t)\le\max\{x_{0},k\}$ in \eqref{eq:apriori} is considerably stronger than the exponential bound obtained by discarding the quadratic term, and it is what makes the sharp asymptotic estimate of \cref{thm:dissipative}(i) available. Throughout, $|\cdot|$ denotes the Euclidean norm of $\R^{2}$, and for a point $v$ and nonempty sets $C,D\subset\R^{2}$ we write
\[
\dist(v,C):=\inf_{c\in C}|v-c|,
\qquad
\dist(D,C):=\sup_{d\in D}\dist(d,C)
\]
for the distance from a point to a set and for the Hausdorff semidistance between sets. We now record the cocycle structure. Define
\[
\Psi(t,\omega)v_{0}:=v(t;0,\omega,v_{0}),\qquad
t\ge0,\ \omega\in\Omega_{0},\ v_{0}\in\Z ,
\]
and $\Psi(t,\omega):=\mathrm{id}_{\Z}$ for $\omega\in\Omega\setminus\Omega_{0}$, which is a $\theta$-invariant null set.

\begin{proposition}\label{prop:rds}
The pair $(\{\theta_{t}\}_{t\in\R},\Psi)$ is a continuous random dynamical system on $\Z$. That is, $\Psi(0,\omega)=\mathrm{id}$, the map $(t,v_{0})\mapsto\Psi(t,\omega)v_{0}$ is continuous for every $\omega$, the map $(t,\omega,v_{0})\mapsto\Psi(t,\omega)v_{0}$ is measurable, and
	\begin{equation}\label{eq:cocycleprop}
		\Psi(t+s,\omega)=\Psi(t,\theta_{s}\omega)\circ\Psi(s,\omega),
		\qquad t,s\ge0,\ \omega\in\Omega .
	\end{equation}
Moreover, for every $t>0$ and every compact $K\subset\Z$ there is $L_{t,K}>0$ with the following property. If $\varphi,\varphi':[0,t]\to[\pmin,\pmax]$ are measurable and $v,v'$ are the solutions of $\dot v=\mathcal F(v,\varphi(s))$ and $\dot v=\mathcal F(v,\varphi'(s))$ starting from $v_{0},v_{0}'\in K$ at time $0$, and both remain in $K$ on $[0,t]$, then
	\begin{equation}\label{eq:contdep}
		\sup_{s\in[0,t]}\bigl|v(s)-v'(s)\bigr|
		\ \le\ L_{t,K}\Bigl(|v_{0}-v_{0}'|+\operatorname*{ess\,sup}_{s\in[0,t]}
		\bigl|\varphi(s)-\varphi'(s)\bigr|\Bigr).
	\end{equation}
In particular, this applies to $v=\Psi(\cdot,\omega)v_{0}$ and $v'=\Psi(\cdot,\omega')v_{0}'$ with $\omega,\omega'\in\Omega_{0}$.
\end{proposition}

\begin{proof}
Continuity in $(t,v_{0})$ is the continuous dependence on initial data for locally Lipschitz vector fields. By \eqref{eq:shift} and \eqref{eq:OU} we have $z^{*}(\theta_{t}(\theta_{s}\omega))=z^{*}(\theta_{t+s}\omega)$, hence $\varphi_{\theta_{s}\omega}(t)=\varphi_{\omega}(t+s)$ for all $t,s\in\R$ by \eqref{eq:phi}. Therefore $t\mapsto v(t+s;0,\omega,v_{0})$ solves \eqref{eq:sys} with $\omega$ replaced by $\theta_{s}\omega$ and initial datum $v(s;0,\omega,v_{0})$, and uniqueness gives \eqref{eq:cocycleprop} on $\Omega_{0}$, which is $\theta$-invariant. Off $\Omega_{0}$ both sides of \eqref{eq:cocycleprop} are the identity; see \cite[Remark 2.2.3]{Arnold} for this standard extension. For \eqref{eq:contdep}, the partial derivatives of $\mathcal F$ in $v$ and in $\varphi$ are continuous on $\operatorname{conv}(K)\times[\pmin,\pmax]$, so there are $L,L'$ with $|\mathcal F(v,\varphi)-\mathcal F(v',\varphi')|\le L|v-v'|+L'|\varphi-\varphi'|$ on the convex hull of $K$, which is contained in $\Z$, and Gronwall's inequality gives \eqref{eq:contdep} with $L_{t,K}:=\max\{1,L't\}\,e^{Lt}$. Finally, fix $t\ge0$ and $v_{0}\in\Z$, and set $\varphi_{\omega}\equiv p$ for $\omega\notin\Omega_{0}$. The map $\omega\mapsto\varphi_{\omega}|_{[0,t]}\in C([0,t])$ is measurable, because the Borel $\sigma$-algebra of $C([0,t])$ is generated by the evaluations, $\Omega_{0}$ is measurable and each $\omega\mapsto\varphi_{\omega}(s)=\psi(z^{*}(\theta_{s}\omega))$ is measurable on $\Omega_{0}$. By \eqref{eq:contdep}, applied in the compact set $K=[0,\max\{x_{0},k\}]\times[0,y_{0}e^{(c\beta/a)t}]$, which contains all the trajectories starting from $v_{0}$ on $[0,t]$ by \eqref{eq:apriori}, the solution at time $t$ is a Lipschitz function of $\varphi_{\omega}|_{[0,t]}$ on the closed subset of $C([0,t])$ of functions with values in $[\pmin,\pmax]$. Hence $\omega\mapsto\Psi(t,\omega)v_{0}$, which is this composition on $\Omega_{0}$ and the constant $v_{0}$ on the measurable set $\Omega\setminus\Omega_{0}$, is measurable for each $(t,v_{0})$. Joint measurability follows by writing $\Psi(t,\omega)v_{0}$ as the pointwise limit, by continuity in $(t,v_{0})$, of the measurable maps $(t,\omega,v_{0})\mapsto\Psi(t_{n},\omega)v_{0,n}$, where $t_{n}$ and $v_{0,n}$ denote the dyadic approximations of $t$ and $v_{0}$ of order $n$, projected onto $\Z$.
\end{proof}

\section{Dissipativity, attracting sets and the pullback attractor}
\label{sec:dissipative}

Throughout this section $\omega\in\Omega_{0}$ is fixed and we abbreviate $v(t)=(x(t),y(t))=\Psi(t,\omega)v_{0}$. Every estimate of \cref{ssec:bounds,ssec:absorbing} uses only $\pmin\le\varphi_{\omega}\le\pmax$, so each of them holds verbatim, with the same constants, when $\varphi_{\omega}$ is replaced by an arbitrary measurable function $\varphi:[0,\infty)\to[\pmin,\pmax]$ in \eqref{eq:sys}. This observation is used in \cref{ssec:uniform}. We shall need the following elementary lemma.

\begin{lemma}\label{lem:gronwall}
Let $u:[0,\infty)\to\R$ be absolutely continuous and suppose that $\dot u(t)\le A-\delta u(t)$ for a.e. $t\ge0$, where $A\in\R$ and $\delta>0$. Then
	\begin{equation}\label{eq:gronwall}
		u(t)\le u(0)e^{-\delta t}+\frac{A}{\delta}\bigl(1-e^{-\delta t}\bigr)
		\qquad\text{for all }t\ge0 .
	\end{equation}
\end{lemma}

\begin{proof}
The function $t\mapsto e^{\delta t}u(t)$ is absolutely continuous with derivative $e^{\delta t}(\dot u+\delta u)\le Ae^{\delta t}$ a.e., and an integration on $[0,t]$ gives \eqref{eq:gronwall}.
\end{proof}

\subsection{Asymptotic bounds}\label{ssec:bounds}

\begin{Theorem}[Dissipativity]\label{thm:dissipative}
Assume \textbf{(H1)}--\textbf{(H3)}. For every $v_{0}\in\Z$ the solution of \eqref{eq:sys} satisfies the following properties.
\begin{enumerate}
\item $\displaystyle\limsup_{t\to+\infty}x(t)\le k$.
\item For every $\nu\in(0,\gamma]$, with $W:=x+\frac1c y$,
		\begin{equation}\label{eq:Wbound}
			W(t)\ \le\ W(0)\,e^{-\nu t}
			+\frac{k(\alpha+\nu)^{2}}{4\alpha\nu}\bigl(1-e^{-\nu t}\bigr)
			\qquad\text{for all }t\ge0 .
		\end{equation}
\item Consequently, setting
		\begin{equation}\label{eq:Lambda}
			\nu^{*}:=\min\{\alpha,\gamma\},
			\qquad
			\Lambda:=\frac{k(\alpha+\nu^{*})^{2}}{4\alpha\nu^{*}}
			=\begin{cases} k, & \gamma\ge\alpha,\\[4pt]
				\dfrac{k(\alpha+\gamma)^{2}}{4\alpha\gamma}, & \gamma<\alpha,
			\end{cases}
		\end{equation}
we have $\displaystyle\limsup_{t\to+\infty}\Bigl(x(t)+\frac1c y(t)\Bigr)\le\Lambda$, and in particular $\displaystyle\limsup_{t\to+\infty}y(t)\le c\Lambda$.
\end{enumerate}
Moreover, $\Lambda\ge k$, with equality if and only if $\gamma\ge\alpha$, and none of the bounds depends on $\omega$.
\end{Theorem}

\begin{proof}
(i) By \eqref{eq:logisticcomp} and comparison, $x(t)\le u(t)$ with $u$ the logistic solution starting from $x_{0}$, and $u(t)\to k$ if $x_{0}>0$, while $u\equiv0$ if $x_{0}=0$.

(ii) If the second equation of \eqref{eq:sys} is multiplied by $1/c$ and added to the first, the predation terms cancel exactly and
	\begin{equation}\label{eq:Wdot}
		\frac{dW}{dt}=\alpha x\Bigl(1-\frac{x}{k}\Bigr)-\frac{\gamma}{c}\,y .
	\end{equation}
Let $\nu\in(0,\gamma]$. Since $y\ge0$,
	\begin{equation}\label{eq:WdotnuW}
		\frac{dW}{dt}+\nu W
		=g_{\nu}(x)+\Bigl(\frac{\nu-\gamma}{c}\Bigr)y\ \le\ g_{\nu}(x)
		\ \le\ \frac{k(\alpha+\nu)^{2}}{4\alpha},
		\qquad g_{\nu}(x):=(\alpha+\nu)x-\frac{\alpha}{k}x^{2},
	\end{equation}
the last inequality being the maximum of the parabola $g_{\nu}$, attained at $x_{\nu}=(\alpha+\nu)k/(2\alpha)$. \Cref{lem:gronwall} gives \eqref{eq:Wbound}.

(iii) Letting $t\to+\infty$ in \eqref{eq:Wbound} gives $\limsup_{t}W(t)\le h(\nu)$ for every $\nu\in(0,\gamma]$, where $h(\nu):=k(\alpha+\nu)^{2}/(4\alpha\nu)$. Since $h'(\nu)=k(\alpha+\nu)(\nu-\alpha)/(4\alpha\nu^{2})$, the function $h$ is strictly decreasing on $(0,\alpha]$ and strictly increasing on $[\alpha,\infty)$, so its minimum over $(0,\gamma]$ is attained at $\nu^{*}$, and $h(\nu^{*})=\Lambda$ is \eqref{eq:Lambda}. Since $y=c(W-x)\le cW$, the bound on $y$ follows. Finally, $h(\nu)/k-1=(\alpha-\nu)^{2}/(4\alpha\nu)\ge0$, with equality if and only if $\nu=\alpha$, which is admissible exactly when $\alpha\le\gamma$.
\end{proof}

For $\gamma\ge\alpha$, \cref{thm:dissipative} says that the triangle $\{x+y/c\le k\}$ attracts, that is, the total biomass measured in prey units is asymptotically bounded by the carrying capacity of the prey alone. One cannot simply take $\nu=\alpha$ in \eqref{eq:WdotnuW}, because the constraint $\nu\le\gamma$ is what makes $(\nu-\gamma)y/c$ nonpositive.

\subsection{The deterministic attracting set}\label{ssec:absorbing}

The next statement is the analogue, for the present model, of the deterministic attracting set that is the standard asymptotic object for population models driven by bounded real noise; see \cite{CCLR2019,CaraballoLopez2020,CGLR2019,LopezGut2026}.

\begin{corollary}[Deterministic compact absorbing and attracting sets]
\label{cor:absorbing} For every $\eta>0$ the set
	\begin{equation}\label{eq:Keta}
		K_{\eta}:=\Bigl\{(x,y)\in\Z\,:\ x\le k+\eta\ \ \text{and}\ \
		x+\frac1c y\le\Lambda+\eta\Bigr\}
	\end{equation}
is compact and positively invariant under $\Psi(t,\omega)$ for every $\omega\in\Omega_{0}$, and it absorbs bounded sets uniformly in $\omega$, that is, for every bounded $B\subset\Z$ there is $T_{B}(\eta)>0$, independent of $\omega$, such that
	\begin{equation}\label{eq:absorb}
		\Psi(t,\omega)B\subset K_{\eta}\qquad\text{for all }t\ge T_{B}(\eta),\
		\omega\in\Omega_{0}.
	\end{equation}
Moreover, $K_{0}=\bigcap_{\eta>0}K_{\eta}$ is compact, positively invariant and forward attracting, uniformly in $\omega$, in the sense that, for every bounded $B\subset\Z$,
	\begin{equation}\label{eq:K0attract}
		\lim_{t\to+\infty}\ \sup_{\omega\in\Omega_{0}}\ \sup_{v_{0}\in B}\
		\dist\bigl(\Psi(t,\omega)v_{0},\,K_{0}\bigr)=0 .
	\end{equation}
The same statements hold for the solutions of \eqref{eq:sys} with $\varphi_{\omega}$ replaced by any measurable $\varphi:[0,\infty)\to[\pmin,\pmax]$, with the same $T_{B}(\eta)$.
\end{corollary}

\begin{proof}
The set $K_{\eta}$ is closed and contained in $[0,k+\eta]\times[0,c(\Lambda+\eta)]$, hence compact. For the positive invariance let $\eta\ge0$ and $v_{0}\in K_{\eta}$. Then $x(t)\le\max\{x_{0},k\}\le k+\eta$ for all $t\ge0$ by \eqref{eq:apriori}, and \eqref{eq:WdotnuW} with $\nu=\nu^{*}$, \eqref{eq:Lambda} and \cref{lem:gronwall} give $W(t)\le W(0)e^{-\nu^{*}t}+\Lambda(1-e^{-\nu^{*}t})\le\Lambda+\eta$, so $v(t)\in K_{\eta}$ for all $t\ge0$. For the absorption let $B\subset\{x\le R_{1},\,y\le R_{2}\}$. By \eqref{eq:logisticcomp} and the monotonicity of the logistic flow with respect to the initial datum, $x(t)\le u_{R_{1}}(t)$, where $u_{R_{1}}$ is the logistic solution starting from $R_{1}$, which converges to $k$, so there is $T^{(1)}=T^{(1)}(R_{1},\eta)$ with $x(t)\le k+\eta$ for $t\ge T^{(1)}$. By \eqref{eq:Wbound} with $\nu=\nu^{*}$ and $W(0)\le R_{1}+R_{2}/c$, there is $T^{(2)}=T^{(2)}(R_{1},R_{2},\eta)$ with $W(t)\le\Lambda+\eta$ for $t\ge T^{(2)}$, and $T_{B}(\eta):=\max\{T^{(1)},T^{(2)}\}$ is the required time. The set $K_{0}$ is the intersection of the decreasing family of compact sets $K_{\eta}$, hence compact, it is the set \eqref{eq:Keta} with $\eta=0$, and it is positively invariant by the case $\eta=0$ of the argument above. Next, $\dist(K_{\eta},K_{0})\to0$ as $\eta\downarrow0$. Otherwise, since $K_{\eta}$ increases with $\eta$, there would be $\varrho>0$ and points $u_{n}\in K_{1/n}$ with $\dist(u_{n},K_{0})\ge\varrho$, and a limit point $u$ of the sequence $(u_{n})$, which exists because $K_{1}$ is compact, would belong to every $K_{\eta}$ and hence to $K_{0}$, while $\dist(u,K_{0})\ge\varrho$, a contradiction. Given $\varrho>0$, we choose $\eta>0$ with $\dist(K_{\eta},K_{0})\le\varrho$, and \eqref{eq:absorb} gives $\dist(\Psi(t,\omega)v_{0},K_{0})\le\varrho$ for all $t\ge T_{B}(\eta)$, $v_{0}\in B$ and $\omega\in\Omega_{0}$, which is \eqref{eq:K0attract}. The last assertion holds because only $\pmin\le\varphi\le\pmax$ was used.
\end{proof}

\subsection{The pullback attractor}\label{ssec:pullback}

Let $\mathcal D$ denote the family of all nonempty bounded subsets of $\Z$. We use the standard notions of the theory of random dynamical systems, for which we refer to \cite{Arnold,CaraballoHan,CrauelFlandoli}. A random compact set $\{\mathcal A(\omega)\}_{\omega\in\Omega}$, that is, a family of nonempty compact subsets of $\Z$ such that $\omega\mapsto\dist(v,\mathcal A(\omega))$ is measurable for every $v\in\Z$, is a $\mathcal D$-pullback attractor for $\Psi$ if it is invariant, $\Psi(t,\omega)\mathcal A(\omega)=\mathcal A(\theta_{t}\omega)$ for all $t\ge0$ and $\omega\in\Omega_{0}$, and it pullback attracts every $B\in\mathcal D$, that is,
\[
\lim_{t\to\infty}\dist\bigl(\Psi(t,\theta_{-t}\omega)B,\,
\mathcal A(\omega)\bigr)=0\qquad\text{for every }\omega\in\Omega_{0}.
\]
Since $\Omega_{0}$ is $\theta$-invariant and of full measure, the restriction of $\Psi$ to $\Omega_{0}$ is a random dynamical system over the metric dynamical system $(\Omega_{0},\mathscr F\cap\Omega_{0},\PP,\theta)$, and it is to this restriction that the results quoted below are applied; off $\Omega_{0}$ the maps $\Psi(t,\omega)$ are the identity and absorb nothing, which is immaterial for statements made up to null sets. Random sets constructed on $\Omega_{0}$ are extended to $\Omega$ by a fixed deterministic compact set; for the pullback attractor we put $\mathcal A(\omega):=K_{1}$ for $\omega\notin\Omega_{0}$, which keeps it measurable and, since $\Psi(t,\omega)$ is the identity and $\theta_{t}\omega\notin\Omega_{0}$ for such $\omega$, invariant on all of $\Omega$. For a bounded set $B\subset\Z$ and $\omega\in\Omega_{0}$ the pullback $\omega$-limit set of $B$ is
\begin{equation}\label{eq:pullbackomega}
	\mathcal O_{B}(\omega):=\bigcap_{t\ge0}\overline{\bigcup_{s\ge t}\Psi(s,\theta_{-s}\omega)B} .
\end{equation}
We shall use the uniqueness theorem of Crauel \cite[Corollaries 5.7 and 5.8]{Crauel1999}, according to which, on a Polish phase space, an invariant random compact set that pullback attracts every deterministic compact set is almost surely contained in every other invariant random compact set with the same attraction property, so that two such sets coincide up to a null set.

\begin{Theorem}[Existence of the pullback attractor]\label{thm:attractor}
Assume \textbf{(H1)}--\textbf{(H3)}. Then $\Psi$ possesses a $\mathcal D$-pullback attractor $\{\mathcal A(\omega)\}$, which is unique, up to a $\PP$-null set, among the invariant random compact sets that pullback attract every compact deterministic set. For every $\omega\in\Omega_{0}$,
	\begin{equation}\label{eq:attractorformula}
		\mathcal A(\omega)=\bigcap_{t\ge0}\Psi\bigl(t,\theta_{-t}\omega\bigr)K_{\eta}
		\qquad(\text{any fixed }\eta>0),
	\end{equation}
and $\mathcal A(\omega)\subset K_{0}$. Moreover, $v_{0}\in\mathcal A(\omega)$ if and only if there exists a complete trajectory $\xi:\R\to K_{0}$, that is, a map with $\Psi(t-s,\theta_{s}\omega)\xi(s)=\xi(t)$ for all $s\le t$, such that $\xi(0)=v_{0}$.
\end{Theorem}

\begin{proof}
By \cref{prop:rds} the pair $(\theta,\Psi)$ is a random dynamical system, and by \cref{cor:absorbing} the compact deterministic set $K_{\eta}$ absorbs every bounded deterministic set, with a deterministic absorption time. The existence of the pullback attractor is \cite[Theorems 4.2 and 4.3]{CaraballoHan}. For the uniqueness, let $\{\mathcal A'(\omega)\}$ be an invariant random compact set that pullback attracts every compact deterministic set. The inclusion $\mathcal A(\omega)\subset\mathcal A'(\omega)$ holds for every $\omega\in\Omega_{0}$ without any measure theory, since $\mathcal A(\omega)=\Psi(t,\theta_{-t}\omega)\mathcal A(\theta_{-t}\omega)\subset\Psi(t,\theta_{-t}\omega)K_{0}$ for every $t\ge0$ and the deterministic compact set $K_{0}$ is pullback attracted by $\mathcal A'$; the reverse inclusion, up to a null set, is Crauel's theorem recalled above, since $\mathcal A$ attracts every compact deterministic set.

The attractor is the closure of the union over $B\in\mathcal D$ of the pullback $\omega$-limit sets\footnote{The classical term, unrelated to the realisation $\omega$.} \eqref{eq:pullbackomega}; see \cite[Theorem 3.11]{CrauelFlandoli}. For $\omega\in\Omega_{0}$ and $s\ge T:=T_{B}(\eta)$ the cocycle property and \eqref{eq:absorb} give $\Psi(s,\theta_{-s}\omega)B=\Psi(s-T,\theta_{-(s-T)}\omega)\Psi(T,\theta_{-s}\omega)B \subset\Psi(s-T,\theta_{-(s-T)}\omega)K_{\eta}$, so that $\mathcal O_{B}(\omega)\subset\mathcal O_{K_{\eta}}(\omega)\subset\mathcal A(\omega)$ and $\mathcal A(\omega)=\mathcal O_{K_{\eta}}(\omega)$. Since $K_{\eta}$ is positively invariant, the cocycle property gives $\Psi(s,\theta_{-s}\omega)K_{\eta}=\Psi(t,\theta_{-t}\omega)\Psi(s-t,\theta_{-s}\omega)K_{\eta} \subset\Psi(t,\theta_{-t}\omega)K_{\eta}$ for $s\ge t$, so the family on the right of \eqref{eq:attractorformula} is decreasing, consists of compact sets, and its intersection is the pullback $\omega$-limit set, which proves \eqref{eq:attractorformula}.

Formula \eqref{eq:attractorformula} also yields directly, for every $\omega\in\Omega_{0}$, the three properties of the attractor. The invariance follows from the cocycle property, since $\Psi(t,\omega)$ is continuous and commutes with decreasing intersections of compact sets. The pullback attraction of $K_{\eta}$, hence of every $B\in\mathcal D$ by the inclusion displayed above, follows because the Hausdorff semidistance from a decreasing family of compact sets to its intersection tends to zero. The measurability of $\mathcal A$ follows because $\dist(v,\mathcal A(\omega))=\lim_{n}\dist(v,\Psi(n,\theta_{-n}\omega)K_{\eta})$, where each term is the infimum over a countable dense subset of $K_{\eta}$ of the measurable functions $\omega\mapsto|v-\Psi(n,\theta_{-n}\omega)v'|$ of \cref{prop:rds}. The inclusion $\mathcal A(\omega)\subset K_{0}$ holds because \eqref{eq:attractorformula} is valid for every $\eta>0$. For the last assertion, if $v_{0}\in\mathcal A(\omega)$, then by invariance we may pick, for each $t\ge0$, $v_{-t}\in\mathcal A(\theta_{-t}\omega)\subset K_{0}$ with $\Psi(t,\theta_{-t}\omega)v_{-t}=v_{0}$. The maps $\Psi(s,\theta_{-s}\omega)$ are injective, because the field is $C^{1}$ on the open set $\mathcal V$ and solutions are unique backwards as well, so the family $\{v_{-t}\}$ is consistent and defines the required $\xi$. Conversely, a complete trajectory contained in $K_{0}$ has, by \eqref{eq:attractorformula}, its value at $0$ in $\mathcal A(\omega)$.
\end{proof}

\subsection{The minimal uniform attracting set}\label{ssec:uniform}

\begin{definition}\label{def:uniform}
A closed set $C\subset\Z$ is called a \emph{uniform attracting set} for $\Psi$ if
	\begin{equation}\label{eq:unifattract}
		\lim_{t\to+\infty}\ \sup_{\omega\in\Omega_{0}}\ \sup_{v_{0}\in B}\
		\dist\bigl(\Psi(t,\omega)v_{0},\,C\bigr)=0
		\qquad\text{for every bounded }B\subset\Z .
	\end{equation}
\end{definition}

The adjective \emph{uniform} refers to uniformity in the realisation $\omega$; a uniform attracting set is thus a deterministic set attracting forward in time uniformly in the realisation, and the minimal one is constructed below. In the language of \cite[Chapter IV]{ChepyzhovVishik2002}, the maps $\Pi_{\omega}(t,\tau):=\Psi(t-\tau,\theta_{\tau}\omega)$, $t\ge\tau$, form a family of processes indexed by the symbol $\omega\in\Omega_{0}$, the symbol space carrying no topology, and \eqref{eq:unifattract} is the uniform attraction of \cite[Definition IV.3.4]{ChepyzhovVishik2002}, because $\{\theta_{\tau}\omega\,:\,\omega\in\Omega_{0}\}=\Omega_{0}$ for every $\tau$. The same words are used in \cite{CuiLanga2017} for a random object that is uniform with respect to a nonautonomous symbol; here the symbol is the realisation itself and the attractor is deterministic. By \cref{cor:absorbing}, $K_{0}$ is a uniform attracting set. Fix the compact positively invariant set $K_{1}$ of \eqref{eq:Keta} and set
\begin{equation}\label{eq:STdef}
	\mathcal S_{T}:=\overline{\bigcup_{t\ge T}\ \bigcup_{\omega\in\Omega_{0}}
		\Psi(t,\omega)K_{1}}\ \subset K_{1},
	\qquad
	\mathcal U:=\bigcap_{T\ge0}\mathcal S_{T}.
\end{equation}

\begin{Theorem}[The minimal uniform attracting set]\label{thm:uniformattract}
Assume \textbf{(H1)}--\textbf{(H3)}. Then $\mathcal U$ is a nonempty compact subset of $K_{0}$, it is a uniform attracting set, and it is contained in every uniform attracting set. That is, $\mathcal U$ is the uniform attractor of the family of processes $\{\Pi_{\omega}\}_{\omega\in\Omega_{0}}$ in the sense of \cite[Definition IV.3.5]{ChepyzhovVishik2002}. Moreover,
	\begin{equation}\label{eq:sandwich}
		\bigcup_{\omega\in\Omega_{0}}\mathcal A(\omega)\ \subset\ \mathcal U\
		\subset\ K_{0}.
	\end{equation}
\end{Theorem}

\begin{proof}
Since $K_{1}$ is positively invariant, each $\mathcal S_{T}$ is a nonempty closed subset of the compact set $K_{1}$, and $T\mapsto\mathcal S_{T}$ is nonincreasing, so $\mathcal U$ is nonempty and compact. For every $\eta>0$, \eqref{eq:absorb} with $B=K_{1}$ gives $\mathcal S_{T}\subset K_{\eta}$ for $T\ge T_{K_{1}}(\eta)$, whence $\mathcal U\subset\bigcap_{\eta>0}K_{\eta}=K_{0}$.

We claim that
	\begin{equation}\label{eq:attractsK1}
		\sup_{\omega\in\Omega_{0}}\ \sup_{v_{0}\in K_{1}}\
		\dist\bigl(\Psi(t,\omega)v_{0},\mathcal U\bigr)\ \xrightarrow[t\to+\infty]{}\ 0 .
	\end{equation}
Otherwise there are $\varrho>0$, $t_{n}\to+\infty$, $\omega_{n}\in\Omega_{0}$ and $v_{n}\in K_{1}$ with $\dist(u_{n},\mathcal U)\ge\varrho$, where $u_{n}:=\Psi(t_{n},\omega_{n})v_{n}\in K_{1}$, and we may assume $u_{n}\to u$. For every $T\ge0$ we have $u_{n}\in\mathcal S_{T}$ for $n$ large, hence $u\in\mathcal S_{T}$, so $u\in\mathcal U$, a contradiction. Now let $B\subset\Z$ be bounded and $\tau_{B}:=T_{B}(1)$. For $t\ge\tau_{B}$ and $v_{0}\in B$ the cocycle property gives $\Psi(t,\omega)v_{0}\in\Psi(t-\tau_{B},\theta_{\tau_{B}}\omega)K_{1}$ with $\theta_{\tau_{B}}\omega\in\Omega_{0}$, and \eqref{eq:attractsK1} at time $t-\tau_{B}$ yields \eqref{eq:unifattract} with $C=\mathcal U$.

Let $C$ be a uniform attracting set and $u\in\mathcal U$. By \eqref{eq:STdef} there are $t_{n}\to+\infty$, $\omega_{n}\in\Omega_{0}$ and $v_{n}\in K_{1}$ with $\Psi(t_{n},\omega_{n})v_{n}\to u$, and \eqref{eq:unifattract} with $B=K_{1}$ gives $\dist(\Psi(t_{n},\omega_{n})v_{n},C)\to0$. Hence $u\in C$, as $C$ is closed. This is the minimality property of \cite[Definition IV.3.5]{ChepyzhovVishik2002}. Finally, let $\omega\in\Omega_{0}$ and $t\ge T$. By the invariance of the pullback attractor and $\mathcal A(\theta_{-t}\omega)\subset K_{0}\subset K_{1}$,
	\[
	\mathcal A(\omega)=\Psi\bigl(t,\theta_{-t}\omega\bigr)
	\mathcal A(\theta_{-t}\omega)\ \subset\
	\Psi\bigl(t,\theta_{-t}\omega\bigr)K_{1}\ \subset\ \mathcal S_{T} ,
	\]
and intersecting over $T\ge0$ gives $\mathcal A(\omega)\subset\mathcal U$.
\end{proof}

We now identify $\mathcal U$. Let $\mathcal C$ denote the set of Lebesgue measurable functions $\varphi:\R\to[\pmin,\pmax]$, which we call \emph{admissible controls}. For $\varphi\in\mathcal C$ we consider the control system
\begin{equation}\label{eq:control}
	\dot v=\mathcal F\bigl(v,\varphi(t)\bigr),
\end{equation}
with $\mathcal F$ as in \eqref{eq:field}. By \cref{thm:wellposed} its Carath\'eodory solutions starting from points of $\Z$ at any initial time are unique, global and remain in $\Z$, and the representation \eqref{eq:representation} shows that each coordinate of a solution vanishes either identically or nowhere. A \emph{bounded complete solution} of \eqref{eq:control} is a bounded absolutely continuous function $\xi:\R\to\Z$ that satisfies \eqref{eq:control} almost everywhere for some $\varphi\in\mathcal C$. The \emph{kernel} of the control system is
\begin{equation}\label{eq:kernel}
	\mathcal K:=\bigl\{\xi(0)\,:\ \xi\text{ is a bounded complete solution of \eqref{eq:control}
	for some }\varphi\in\mathcal C\bigr\}.
\end{equation}
By \cref{cor:absorbing} applied to the control system, a bounded complete solution takes values in $K_{0}$, since $\xi(s)$ is the value at time $t$ of the solution starting from $\xi(s-t)$, which lies in $K_{\eta}$ for $t\ge T_{B}(\eta)$ with $B$ the range of $\xi$. The proof of the next theorem rests on an affine reparametrisation of the control. For $x\ge0$ and $w\in[0,1]$ put
\begin{equation}\label{eq:gw}
	\widehat G(x,w):=G(x,\pmin)+w\bigl[G(x,\pmax)-G(x,\pmin)\bigr],
	\qquad
	\widehat{\mathcal F}(v,w):=\Bigl(\alpha x\Bigl(1-\frac{x}{k}\Bigr)-\widehat G(x,w)\,y,\
	-\gamma y+c\,\widehat G(x,w)\,y\Bigr).
\end{equation}
Since $\varphi\mapsto G(x,\varphi)$ is continuous and nondecreasing from $[\pmin,\pmax]$ onto $[G(x,\pmin),G(x,\pmax)]$, strictly increasing if $x>0$, the velocity sets coincide,
\begin{equation}\label{eq:velocity}
	\bigl\{\mathcal F(v,\varphi)\,:\,\varphi\in[\pmin,\pmax]\bigr\}
	=\bigl\{\widehat{\mathcal F}(v,w)\,:\,w\in[0,1]\bigr\}
	\qquad\text{for every }v\in\Z,
\end{equation}
and $\widehat{\mathcal F}(v,w)=\mathcal F_{0}(v)+w\,\mathcal F_{1}(v)$ is affine in $w$, with $\mathcal F_{0}$ and $\mathcal F_{1}$ of class $C^{\infty}$ on $\mathcal V$. For $x>0$ the inverse of $\varphi\mapsto G(x,\varphi)$ is explicit,
\begin{equation}\label{eq:Gamma}
	\Gamma(x,\varsigma):=\frac{\varsigma}{x(\beta-a\varsigma)},\qquad 0\le\varsigma<\frac{\beta}{a},
\end{equation}
and it is continuous on $\{x>0\}\times[0,\beta/a)$.

\begin{lemma}\label{lem:reparam}
Let $I\subset\R$ be an interval.
\begin{enumerate}
\item If $\varphi\in\mathcal C$ and $v=(x,y)$ solves \eqref{eq:control} on $I$, then
$v$ solves $\dot v=\widehat{\mathcal F}(v,w(t))$ a.e. on $I$ for the measurable function $w:I\to[0,1]$ given by
\[
w(t):=\frac{G(x(t),\varphi(t))-G(x(t),\pmin)}{G(x(t),\pmax)-G(x(t),\pmin)}\ \text{ if }x(t)>0,\qquad w(t):=0\ \text{ otherwise.}
\]
\item If $w:I\to[0,1]$ is measurable and $v=(x,y)$ is an absolutely continuous solution of
$\dot v=\widehat{\mathcal F}(v,w(t))$ a.e. on $I$ with values in $\Z$, then $v$ solves \eqref{eq:control} on $I$ for the measurable function $\varphi:I\to[\pmin,\pmax]$ given by $\varphi(t):=\Gamma\bigl(x(t),\widehat G(x(t),w(t))\bigr)$ if $x\not\equiv0$ on $I$ and $\varphi\equiv p$ otherwise. If $w$ and $x$ are continuous, so is $\varphi$.
\end{enumerate}
\end{lemma}

\begin{proof}
In both cases $x$ vanishes either identically or nowhere on $I$, by the representation \eqref{eq:representation} of the first coordinate, which is valid for absolutely continuous solutions. If $x\equiv0$ the two systems reduce to $\dot y=-\gamma y$ and there is nothing to prove. If $x>0$ on $I$, then in (i) $w$ is measurable as a composition of continuous maps with the measurable map $t\mapsto(x(t),\varphi(t))$, it takes values in $[0,1]$ by monotonicity of $G$ in $\varphi$, and $\widehat G(x(t),w(t))=G(x(t),\varphi(t))$ by construction, so the two right-hand sides coincide. In (ii), $\widehat G(x(t),w(t))\in[G(x(t),\pmin),G(x(t),\pmax)]\subset[0,\beta/a)$, so $\varphi(t)$ is well defined, lies in $[\pmin,\pmax]$, satisfies $G(x(t),\varphi(t))=\widehat G(x(t),w(t))$, and is measurable, or continuous when $w$ and $x$ are, because $\Gamma$ and $\widehat G$ are continuous.
\end{proof}

\begin{lemma}\label{lem:weakstar}
Let $T>0$, let $w_{n},w:[0,T]\to[0,1]$ be measurable with $w_{n}\to w$ weak-$*$ in $L^{\infty}(0,T)$, and let $v_{0,n}\to v_{0}$ in $K_{1}$. Let $v_{n}$ and $v$ be the solutions of $\dot v=\widehat{\mathcal F}(v,w_{n}(t))$ and $\dot v=\widehat{\mathcal F}(v,w(t))$ on $[0,T]$ starting from $v_{0,n}$ and $v_{0}$ at time $0$. Then $\sup_{[0,T]}|v_{n}-v|\to0$.
\end{lemma}

\begin{proof}
The components of $\widehat{\mathcal F}(v,w)$ factor as $x\,g_{1}$ and $y\,g_{2}$ with $g_{1},g_{2}$ smooth on $\mathcal V\times[0,1]$, exactly as those of $\mathcal F$, so \cref{thm:wellposed} applies verbatim to $\dot v=\widehat{\mathcal F}(v,w(t))$ with measurable $w$, and its solutions starting from $\Z$ are unique, global and remain in $\Z$. All the solutions in the statement remain in $K_{1}$, because by \cref{lem:reparam}(ii) they solve \eqref{eq:control} for suitable admissible controls and $K_{1}$ is positively invariant for the control system by \cref{cor:absorbing}. Let $L$ be a Lipschitz constant of $\mathcal F_{0}$ and $\mathcal F_{1}$ on $K_{1}$ and $M:=\max_{K_{1}}|\mathcal F_{1}|$. For $t\in[0,T]$,
	\[
	v_{n}(t)-v(t)=v_{0,n}-v_{0}+\int_{0}^{t}\bigl[\mathcal F_{0}(v_{n})-\mathcal F_{0}(v)\bigr]ds
	+\int_{0}^{t}w_{n}\bigl[\mathcal F_{1}(v_{n})-\mathcal F_{1}(v)\bigr]ds+e_{n}(t),
	\]
where $e_{n}(t):=\int_{0}^{t}(w_{n}-w)\mathcal F_{1}(v)\,ds$. For each fixed $t$ we have $e_{n}(t)\to0$, because $\mathbf 1_{[0,t]}\mathcal F_{1}(v)\in L^{1}(0,T)$, and $|e_{n}(t)-e_{n}(t')|\le M|t-t'|$, so the sequence $(e_{n})$ is equicontinuous and converges to $0$ uniformly on $[0,T]$. Since $0\le w_{n}\le1$, Gronwall's inequality gives $\sup_{[0,T]}|v_{n}-v|\le(|v_{0,n}-v_{0}|+\sup_{[0,T]}|e_{n}|)e^{2LT}\to0$.
\end{proof}

\begin{Theorem}[$\mathcal U$ is the kernel of the control system]\label{thm:control}
Assume \textbf{(H1)}--\textbf{(H3)}. Then $\mathcal U=\mathcal K$.
\end{Theorem}

\begin{proof}
\emph{The inclusion $\mathcal K\subset\mathcal U$.} Let $\xi$ be a bounded complete solution of \eqref{eq:control} with control $\varphi\in\mathcal C$, and let $w$ be the function given by \cref{lem:reparam}(i) on $I=\R$. Fix $T>0$ and $t\ge T$, and put $\tilde w(s):=w(s-t)$ for $s\in[0,t]$, so that $\xi(0)$ is the value at time $t$ of the solution of $\dot v=\widehat{\mathcal F}(v,\tilde w(s))$ starting from $\xi(-t)\in K_{0}\subset K_{1}$. Since continuous functions are dense in $L^{1}(0,t)$ and the projection of $\R$ onto $[0,1]$ is $1$-Lipschitz, there are continuous $w_{n}:[0,t]\to[0,1]$ with $w_{n}\to\tilde w$ in $L^{1}(0,t)$. Since $\|w_{n}\|_{\infty}\le1$, the convergence also holds weak-$*$ in $L^{\infty}(0,t)$, because for $f\in L^{1}(0,t)$ and a bounded $f_{M}$ with $\|f-f_{M}\|_{L^{1}(0,t)}<\delta$ one has $|\int_{0}^{t}(w_{n}-\tilde w)f|\le\delta+\|f_{M}\|_{\infty}\|w_{n}-\tilde w\|_{L^{1}(0,t)}$. Let $v_{n}$ be the solution of $\dot v=\widehat{\mathcal F}(v,w_{n}(s))$ starting from $\xi(-t)$. By \cref{lem:weakstar}, $v_{n}(t)\to\xi(0)$. By \cref{lem:reparam}(ii), $v_{n}$ solves \eqref{eq:control} on $[0,t]$ for a continuous $\varphi_{n}:[0,t]\to[\pmin,\pmax]$, and by \cref{lem:longexc} there is $\omega_{n}\in\Omega_{0}$ with $\sup_{[0,t]}|\varphi_{\omega_{n}}-\varphi_{n}|<1/n$. Both $v_{n}$ and $s\mapsto\Psi(s,\omega_{n})\xi(-t)$ remain in $K_{1}$, so \eqref{eq:contdep} gives $|\Psi(t,\omega_{n})\xi(-t)-v_{n}(t)|\le L_{t,K_{1}}/n$. Hence $\Psi(t,\omega_{n})\xi(-t)\to\xi(0)$ with $\xi(-t)\in K_{1}$, so that $\xi(0)\in\overline{\bigcup_{\omega\in\Omega_{0}}\Psi(t,\omega)K_{1}}\subset\mathcal S_{T}$. As $T$ was arbitrary, $\xi(0)\in\mathcal U$.

\emph{The inclusion $\mathcal U\subset\mathcal K$.} Let $u\in\mathcal U$. Since $u\in\mathcal S_{n}$ for every $n\in\N$, there are $t_{n}\ge n$, $\omega_{n}\in\Omega_{0}$ and $v_{n}\in K_{1}$ with $|\Psi(t_{n},\omega_{n})v_{n}-u|<1/n$. Define $\xi_{n}(s):=\Psi(s+t_{n},\omega_{n})v_{n}$ for $s\ge-t_{n}$ and $\xi_{n}(s):=v_{n}$ for $s<-t_{n}$, and let $w_{n}:\R\to[0,1]$ be the function given by \cref{lem:reparam}(i) for the control $s\mapsto\varphi_{\omega_{n}}(s+t_{n})$ on $[-t_{n},\infty)$, extended by $0$. Then $\xi_{n}$ takes values in $K_{1}$, it solves $\dot v=\widehat{\mathcal F}(v,w_{n}(s))$ on $[-t_{n},\infty)$, and it is Lipschitz there with constant $M':=\max_{K_{1}\times[0,1]}|\widehat{\mathcal F}|$. Fix $j\in\N$. For $n\ge j$ the functions $\xi_{n}|_{[-j,j]}$ are uniformly bounded and equi-Lipschitz, and the functions $w_{n}|_{[-j,j]}$ are bounded in $L^{\infty}(-j,j)=L^{1}(-j,j)^{*}$. By the Arzel\`a--Ascoli theorem, the Banach--Alaoglu theorem and the separability of $L^{1}$, a diagonal extraction provides a subsequence, not relabelled, a continuous $\xi:\R\to K_{1}$ and a measurable $w:\R\to\R$ such that $\xi_{n}\to\xi$ uniformly on $[-j,j]$ and $w_{n}\to w$ weak-$*$ in $L^{\infty}(-j,j)$ for every $j$. For every nonnegative $f\in L^{1}(-j,j)$ we have $\int wf=\lim\int w_{n}f\in[0,\int f]$, so $0\le w\le1$ a.e., and we modify $w$ on a null set so that it takes values in $[0,1]$. For $s\in[-j,j]$ and $n\ge j$,
	\[
	\xi_{n}(s)=\xi_{n}(-j)+\int_{-j}^{s}\bigl[\mathcal F_{0}(\xi_{n})+w_{n}\mathcal F_{1}(\xi_{n})\bigr]dr .
	\]
The terms $\xi_{n}(s)$, $\xi_{n}(-j)$ and $\int_{-j}^{s}\mathcal F_{0}(\xi_{n})\,dr$ converge to the corresponding terms with $\xi$, by uniform convergence, and
	\[
	\int_{-j}^{s}w_{n}\mathcal F_{1}(\xi_{n})\,dr-\int_{-j}^{s}w\,\mathcal F_{1}(\xi)\,dr
	=\int_{-j}^{s}w_{n}\bigl[\mathcal F_{1}(\xi_{n})-\mathcal F_{1}(\xi)\bigr]dr+\int_{-j}^{s}(w_{n}-w)\mathcal F_{1}(\xi)\,dr
	\ \longrightarrow\ 0,
	\]
the first integral by uniform convergence and $0\le w_{n}\le1$, the second by weak-$*$ convergence. Hence $\xi$ is a bounded absolutely continuous solution of $\dot v=\widehat{\mathcal F}(v,w(s))$ on $\R$, and $\xi(0)=\lim\xi_{n}(0)=u$. By \cref{lem:reparam}(ii), $\xi$ is a bounded complete solution of \eqref{eq:control} for some $\varphi\in\mathcal C$, so $u\in\mathcal K$.
\end{proof}

For $\varphi\in\mathcal C$ and $s\le t$ let $S_{\varphi}(t,s)$ denote the solution operator of \eqref{eq:control} from time $s$ to time $t$, so that $\Psi(t,\omega)=S_{\varphi_{\omega}}(t,0)$ for $\omega\in\Omega_{0}$.

\begin{corollary}[Invariance properties of $\mathcal U$]\label{cor:invariance}
Assume \textbf{(H1)}--\textbf{(H3)}. Then $S_{\varphi}(t,0)\mathcal U\subset\mathcal U$ for every $\varphi\in\mathcal C$ and $t\ge0$, in particular $\Psi(t,\omega)\mathcal U\subset\mathcal U$ for every $\omega\in\Omega_{0}$, and
	\begin{equation}\label{eq:familyinv}
		\mathcal U=\bigcup_{\varphi\in\mathcal C}S_{\varphi}(t,0)\,\mathcal U
		=\overline{\bigcup_{\omega\in\Omega_{0}}\Psi(t,\omega)\,\mathcal U}
		\qquad\text{for every }t\ge0 .
	\end{equation}
\end{corollary}

\begin{proof}
Let $u=\xi(0)$ with $\xi$ a bounded complete solution for the control $\varphi'\in\mathcal C$, and let $\varphi\in\mathcal C$. The function equal to $\xi$ on $(-\infty,0]$ and to $S_{\varphi}(s,0)u$ on $[0,\infty)$ is absolutely continuous and solves \eqref{eq:control} a.e. for the admissible control equal to $\varphi'$ on $(-\infty,0)$ and to $\varphi$ on $[0,\infty)$, and it is bounded because $u\in K_{0}$ and $K_{0}$ is positively invariant for the control system by \cref{cor:absorbing}. Its translate by $t$ is again a bounded complete solution, so $S_{\varphi}(t,0)u\in\mathcal K=\mathcal U$. Conversely, the translate of $\xi$ by $-t$ shows that $\xi(-t)\in\mathcal K$, and $u=S_{\varphi'(\cdot-t)}(t,0)\xi(-t)$, which gives the first equality in \eqref{eq:familyinv}. For the second, the inclusion $\supset$ follows from the positive invariance and the closedness of $\mathcal U$, and the first half of the proof of \cref{thm:control} shows that $u=\lim\Psi(t,\omega_{n})\xi(-t)$ for suitable $\omega_{n}\in\Omega_{0}$, with $\xi(-t)\in\mathcal U$.
\end{proof}

For $\bar\varphi\in[\pmin,\pmax]$ let $S_{\bar\varphi}(t)$, $t\ge0$, denote the semiflow on $\Z$ generated by the autonomous system obtained from \eqref{eq:sys} by freezing $\varphi_{\omega}\equiv\bar\varphi$, so that $S_{\bar\varphi}(t)=S_{\varphi}(t,0)$ for the constant control $\varphi\equiv\bar\varphi$. Since $K_{1}$ is compact, positively invariant and absorbing for $S_{\bar\varphi}$ by \cref{cor:absorbing}, the semiflow $S_{\bar\varphi}$ possesses a global attractor $\mathcal A_{\bar\varphi}$ \cite[Corollary 2.13]{CarvalhoLangaRobinson2013}, which is the set of values at time $0$ of the bounded complete solutions of the frozen system \cite[Theorem 1.7]{CarvalhoLangaRobinson2013}.

\begin{corollary}\label{cor:frozen}
Assume \textbf{(H1)}--\textbf{(H3)}. Then
	\begin{equation}\label{eq:frozenincluded}
		[0,k]\times\{0\}\ \subset\
		\bigcup_{\bar\varphi\in[\pmin,\pmax]}\mathcal A_{\bar\varphi}\ \subset\ \mathcal U .
	\end{equation}
If moreover $\lambda_{+}>0$, which forces $c\beta>a\gamma$, then the point
	\begin{equation}\label{eq:vplus}
		v^{*}_{\pmax}:=\bigl(x^{*}_{\pmax},y^{*}_{\pmax}\bigr),
		\qquad
		x^{*}_{\pmax}=\frac{\gamma}{\pmax(c\beta-a\gamma)},
		\qquad
		y^{*}_{\pmax}=\frac{\alpha}{\beta\pmax}\Bigl(1-\frac{x^{*}_{\pmax}}{k}\Bigr)
		\bigl(1+a\pmax x^{*}_{\pmax}\bigr),
	\end{equation}
lies in $\Zo\cap\mathcal U$, and consequently
	\begin{equation}\label{eq:Ustrictly}
		[0,k]\times\{0\}\ \subsetneq\ \mathcal U .
	\end{equation}
\end{corollary}

\begin{proof}
A bounded complete solution of the frozen system is a bounded complete solution of \eqref{eq:control} with the constant control $\varphi\equiv\bar\varphi$, so $\mathcal A_{\bar\varphi}\subset\mathcal K=\mathcal U$ by \cref{thm:control}. On the invariant face $\{y=0\}$ every frozen system reduces to the logistic equation $\dot x=\alpha x(1-x/k)$, whose bounded complete orbits fill $[0,k]$, so $[0,k]\times\{0\}\subset\mathcal A_{\bar\varphi}$. If $\lambda_{+}>0$, then $\gamma<\Phi(\pmax)<c\beta/a$ by \cref{lem:Fprops}(ii), so $x^{*}_{\pmax}>0$ is well defined, and the condition $\Phi(\pmax)>\gamma$ is equivalent to $x^{*}_{\pmax}<k$, as recalled after \eqref{eq:detequil}, whence $y^{*}_{\pmax}>0$. By construction $v^{*}_{\pmax}$ is an equilibrium of the frozen system with $\bar\varphi=\pmax$, hence a point of $\mathcal A_{\pmax}$, and the inclusion \eqref{eq:Ustrictly} is strict because $y^{*}_{\pmax}>0$.
\end{proof}

In general a uniform attractor can be strictly larger than the union of the attractors of the individual processes, as Haraux's example shows \cite[Chapter IV, Section 4]{ChepyzhovVishik2002}, because a suitably timed switching between symbols, that is, between controls in the present setting, can sustain states that no fixed symbol sustains. The inclusions \eqref{eq:frozenincluded} are what we use in the sequel.

The set $\mathcal K$ involves neither the probability space nor the parameters $\kappa$ and $\sigma$ of the Ornstein--Uhlenbeck process, but only the band $[\pmin,\pmax]$, and in light of \cref{thm:control} the same is true of $\mathcal U$. In particular, the definition \eqref{eq:STdef} does not depend on the choice of $\Omega_{0}$, since the proofs of \cref{thm:uniformattract,thm:control} use only that $\Omega_{0}$ is a $\theta$-invariant subset of $\Omega_{1}$ of full measure, and $\mathcal U$ is the same for all values of $\kappa$ and $\sigma$, whereas $\lambda$ and the pullback attractor do depend on $\sigma_{*}$. This is the precise sense in which the forward description is a worst-case one, and \cref{lem:longexc} is where the noise enters the proof of \cref{thm:control}. 

In the terminology of \cite{ChepyzhovVishik2002}, the bounded complete solutions of \eqref{eq:control} with a given control form the kernel of the process generated by that control, and $\mathcal K$ is the union over the admissible controls of the sections of these kernels at time $0$, so that \cref{thm:control} is a representation of the uniform attractor by kernel sections of the type of \cite[Theorem IV.5.1]{ChepyzhovVishik2002}, with the symbol space enlarged from the realisations to all admissible controls. 

The precise relation is the following. After the affine reparametrisation \eqref{eq:gw}, the natural symbol space is the set $\mathfrak S$ of measurable functions $w:\R\to[0,1]$ with the weak-$*$ topology of $L^{\infty}_{\rm loc}(\R)$. This space is compact and metrisable, and the translations act continuously on it. The family of processes generated by \eqref{eq:gw} is $(K_{1}\times\mathfrak S,K_{1})$-continuous in the sense of \cite{ChepyzhovVishik2002} by \cref{lem:weakstar} and has the uniform absorbing set of \cref{cor:absorbing}, so that the representation theorem of \cite[Chapter IV]{ChepyzhovVishik2002} gives the uniform attractor of the family over $\mathfrak S$ as the union of the kernel sections, which is $\mathcal K$ by \cref{lem:reparam}. The family of processes generated by the realisations is not a subfamily of this one, because the reparametrisation of \cref{lem:reparam} depends on the trajectory. \Cref{thm:control} proves independently that its uniform attractor $\mathcal U$ is also $\mathcal K$. The link between the two families is made along trajectories, by combining \cref{lem:reparam}(ii) with the shadowing of \cref{lem:longexc}, which shows that on each compact interval the uniform closure of the realisations is the set of continuous band-valued functions, and with the weak-$*$ density of continuous controls. What is specific to the model is therefore \cref{lem:longexc} and the affine structure after \eqref{eq:gw}; the self-contained proof above avoids the general theory at no extra cost.

The identification of \cref{thm:control} also places $\mathcal U$ within the theory of multivalued dynamical systems. The control system \eqref{eq:control} generates on $\Z$ the multivalued semiflow
\begin{equation}\label{eq:multisemiflow}
	\mathcal G(t)v:=\bigl\{S_{\varphi}(t,0)v\,:\,\varphi\in\mathcal C\bigr\},\qquad t\ge0,\ v\in\Z,
\end{equation}
which is a strict multivalued semiflow, $\mathcal G(t+s)=\mathcal G(t)\mathcal G(s)$, because controls can be concatenated and translated, and whose values are the reachable sets of the differential inclusion $\dot v\in\{\mathcal F(v,\varphi)\,:\,\varphi\in[\pmin,\pmax]\}$, the velocity sets of which are segments by \eqref{eq:velocity}.

\begin{corollary}[$\mathcal U$ as a multivalued global attractor]\label{cor:multivalued}
Assume \textbf{(H1)}--\textbf{(H3)}. Then $\mathcal U$ is the global attractor of the multivalued semiflow \eqref{eq:multisemiflow} in the sense of \cite{MelnikValero1998,CaraballoLangaMelnikValero2003}. That is, $\mathcal U$ is compact, $\mathcal G(t)\,\mathcal U=\mathcal U$ for every $t\ge0$, it satisfies
	\begin{equation}\label{eq:multiattract}
		\lim_{t\to+\infty}\ \sup_{\varphi\in\mathcal C}\ \sup_{v_{0}\in B}\
		\dist\bigl(S_{\varphi}(t,0)v_{0},\mathcal U\bigr)=0\qquad\text{for every bounded }B\subset\Z ,
	\end{equation}
and it is contained in every closed set $C\subset\Z$ with the attraction property \eqref{eq:multiattract}.
\end{corollary}

\begin{proof}
The invariance $\mathcal G(t)\mathcal U=\mathcal U$ is the first identity in \eqref{eq:familyinv}. For \eqref{eq:multiattract} it suffices, by the absorption of \cref{cor:absorbing} applied to the control system and the semiflow property, to treat $B=K_{1}$. If \eqref{eq:multiattract} failed for $B=K_{1}$, there would be $\varrho>0$, $t_{n}\to+\infty$, $\varphi_{n}\in\mathcal C$ and $v_{n}\in K_{1}$ with $\dist(S_{\varphi_{n}}(t_{n},0)v_{n},\mathcal U)\ge\varrho$. The second half of the proof of \cref{thm:control}, applied verbatim with $\xi_{n}(s):=S_{\varphi_{n}}(s+t_{n},0)v_{n}$ for $s\ge-t_{n}$ and with the function $w_{n}$ given by \cref{lem:reparam}(i) for the control $\varphi_{n}(\cdot+t_{n})$, in which nothing depends on the controls being realisations of the noise, provides a subsequence along which $\xi_{n}(0)$ converges to the value at time $0$ of a bounded complete solution of \eqref{eq:control}, that is, to a point of $\mathcal K=\mathcal U$, a contradiction. Finally, if $C$ is closed and has the attraction property \eqref{eq:multiattract}, then every $u\in\mathcal U$, being of the form $u=S_{\varphi(\cdot-t)}(t,0)\xi(-t)$ with $\xi(-t)\in K_{0}$ for every $t\ge0$, satisfies $\dist(u,C)\le\sup_{\varphi}\sup_{v_{0}\in K_{0}}\dist(S_{\varphi}(t,0)v_{0},C)\to0$, so $\mathcal U\subset C$ and $\mathcal U$ is minimal.
\end{proof}

\subsection{Forward versus pullback}\label{ssec:forwardvspullback}

We now have three nested objects, ordered as in \eqref{eq:sandwich}. The set $K_{0}$ is explicit but not minimal. The set $\mathcal U$ is deterministic and minimal, it is invariant under the family of all realisations by \eqref{eq:familyinv}, and it is governed by the upper worst-case exponent $\lambda_{+}$ (\cref{cor:interior}). The pullback attractor $\mathcal A(\omega)$ has the complementary profile. It is invariant in the cocycle sense and made of complete trajectories, but it is $\omega$-dependent, it is reached in the pullback sense, forward in time it attracts only in probability \cite[Section 4.2]{CaraballoHan}, and it is governed by ergodic averages of the noise, hence by $\lambda$ (\cref{thm:attractorext,thm:attractorpersist}).

\Cref{thm:charext,thm:lowerbound} make the comparison precise in the two regimes. Neither object determines the other. The set $\mathcal U$ is not determined by the fibres $\mathcal A(\omega)$, since by \cref{thm:Uext,thm:Ustrict} the two coincide when $\lambda_{+}<0$ and differ whenever $\lambda<0<\lambda_{+}$, the fibres being the segment $[0,k]\times\{0\}$ in both cases. Nor are the fibres determined by $\mathcal U$, since $\mathcal U$ depends on the noise only through the band, whereas on a band on which $\lambda$ changes sign as $\sigma_{*}$ varies, and such bands exist by \cref{prop:kapsig} whenever $-\gamma+\frac12[\Phi(\pmin)+\Phi(\pmax)]<0<\lambda_{\rm det}$, the fibres pass from the segment to sets meeting the coexistence region (\cref{thm:attractorext,thm:attractorpersist}).

\section{Extinction of the predator}\label{sec:extinction}

In this section we prove that a negative value of the single number $\lambda$ of \cref{def:lambda} forces the extinction of the predator at the exact rate $\lambda$, and we compare the forward and the pullback descriptions in the extinction regime.

\begin{Theorem}[Extinction of the predator]\label{thm:extinction}
Assume \textbf{(H1)}--\textbf{(H3)}, let $\omega\in\Omega_{0}$ and suppose that \eqref{eq:E} holds. Then, for every $v_{0}=(x_{0},y_{0})\in\Z$ with $y_{0}>0$,
	\begin{equation}\label{eq:yexp}
		\limsup_{t\to+\infty}\frac1t\log y(t)\ \le\ \lambda\ <\ 0 ,
	\end{equation}
so that $y(t)\to0$ exponentially fast as $t\to+\infty$. If moreover $x_{0}>0$, then
	\begin{equation}\label{eq:xtok}
		\lim_{t\to+\infty}x(t)=k
		\qquad\text{and}\qquad
		\lim_{t\to+\infty}\frac1t\log y(t)=\lambda .
	\end{equation}
\end{Theorem}

\begin{proof}
By \cref{lem:lambdaeta} we may pick $\eta\in\mathbb Q\cap(0,\infty)$ with $\lambda_{\eta}<0$, and by \cref{thm:dissipative}(i) there is $T\ge0$ with $x(t)\le k+\eta$ for all $t\ge T$. Since $F$ is increasing in its first variable, $F(x(t),\varphi_{\omega}(t))\le\Phi_{\eta}(\varphi_{\omega}(t))$ for $t\ge T$. Since $y_{0}>0$ we have $y(t)>0$ for all $t\ge0$, and the second identity in \eqref{eq:representation} gives
	\[
	\log y(t)-\log y(T)=\int_{T}^{t}
	\Bigl[-\gamma+F\bigl(x(s),\varphi_{\omega}(s)\bigr)\Bigr]ds
	\ \le\ \int_{T}^{t}\Bigl[-\gamma+\Phi_{\eta}
	\bigl(\varphi_{\omega}(s)\bigr)\Bigr]ds ,\qquad t\ge T .
	\]
After division by $t$, \cref{lem:ergodic} applied to $g=\Phi_{\eta}$, which is legitimate because $\omega\in\Omega_{0}\subset\Omega_{\Phi_{\eta}}$ and the finite interval $[0,T]$ does not affect the Ces\`aro limit, gives as $t\to+\infty$
	\[
	\limsup_{t\to+\infty}\frac1t\log y(t)\ \le\
	-\gamma+\E\bigl[\Phi_{\eta}(\psi(z^{*}))\bigr]=\lambda_{\eta}.
	\]
This holds for every rational $\eta>0$ with $\lambda_{\eta}<0$, and letting $\eta\downarrow0$ along the rationals we obtain \eqref{eq:yexp} by \cref{lem:lambdaeta}. In particular, there are $\lambda'\in(\lambda,0)$ and $T'\ge0$ with $y(t)\le y(T')e^{\lambda'(t-T')}$ for $t\ge T'$.

Assume now $x_{0}>0$, whence $x(t)>0$ for all $t\ge0$, and put
	\[
	q(t):=\frac{\beta\varphi_{\omega}(t)y(t)}{1+a\varphi_{\omega}(t)x(t)}
	\in\bigl[0,\ \beta\pmax\,y(t)\bigr],
	\]
so that $q(t)\to0$. Let $\varrho\in(0,\alpha)$ be arbitrary and choose $T_{\varrho}$ with $q(t)\le\varrho$ for $t\ge T_{\varrho}$. The first equation of \eqref{eq:sys} gives, for $t\ge T_{\varrho}$,
	\[
	\frac{dx}{dt}=x\Bigl[\alpha\Bigl(1-\frac{x}{k}\Bigr)-q(t)\Bigr]
	\ \ge\ (\alpha-\varrho)\,x\Bigl(1-\frac{x}{k_{\varrho}}\Bigr),
	\qquad k_{\varrho}:=k\Bigl(1-\frac{\varrho}{\alpha}\Bigr)>0,
	\]
and comparison with the logistic equation starting from $x(T_{\varrho})>0$ yields $\liminf_{t\to\infty}x(t)\ge k_{\varrho}$. Together with \cref{thm:dissipative}(i) this gives $k(1-\varrho/\alpha)\le\liminf_{t}x(t)\le\limsup_{t}x(t)\le k$, and letting $\varrho\downarrow0$ proves the first limit in \eqref{eq:xtok}. For the second, we use that $|\partial_{x}F|\le c\beta\pmax$ on $[0,\infty)\times[\pmin,\pmax]$ by \cref{lem:Fprops}(i), so that
	\begin{equation}\label{eq:unifclose}
		\Bigl|F\bigl(x(t),\varphi_{\omega}(t)\bigr)-\Phi\bigl(\varphi_{\omega}(t)\bigr)
		\Bigr|\ \le\ c\beta\pmax\,\bigl|x(t)-k\bigr|\ \xrightarrow[t\to+\infty]{}\ 0 .
	\end{equation}
By \eqref{eq:representation},
	\[
	\frac1t\log\frac{y(t)}{y_{0}}
	=\frac1t\int_{0}^{t}\Bigl[-\gamma+\Phi\bigl(\varphi_{\omega}(s)\bigr)\Bigr]ds
	+\frac1t\int_{0}^{t}\Bigl[F\bigl(x(s),\varphi_{\omega}(s)\bigr)
	-\Phi\bigl(\varphi_{\omega}(s)\bigr)\Bigr]ds ,
	\]
where the first term converges to $\lambda$ by \cref{lem:ergodic} and the second converges to $0$ because its integrand does, by \eqref{eq:unifclose}.
\end{proof}

\begin{Theorem}[The attractor in the extinction regime]\label{thm:attractorext}
Assume \textbf{(H1)}--\textbf{(H3)} and \eqref{eq:E}. Then, for every $\omega\in\Omega_{0}$,
	\[
	\mathcal A(\omega)=[0,k]\times\{0\} .
	\]
In particular, the pullback attractor is deterministic, and it consists of the two equilibria $(0,0)$ and $(k,0)$ together with the heteroclinic orbit joining them inside $\{y=0\}$.
\end{Theorem}

\begin{proof}
The face $\{y=0\}$ is invariant in both time directions, and on it the first equation of \eqref{eq:sys} becomes the autonomous logistic equation $\dot x=\alpha x(1-x/k)$. For each $x_{\sharp}\in[0,k]$ the logistic equation has a unique complete solution through $x_{\sharp}$, with values in $[0,k]$, and since $\Lambda\ge k$ we have $[0,k]\times\{0\}\subset K_{0}$. Hence \cref{thm:attractor} gives $[0,k]\times\{0\}\subset\mathcal A(\omega)$.

Conversely, let $v_{0}\in\mathcal A(\omega)$ and let $\xi=(x,y):\R\to K_{0}$ be a complete trajectory with $\xi(0)=v_{0}$. Suppose, for contradiction, that $y(0)>0$. Then $y(t)>0$ for every $t\in\R$, forward by \cref{thm:wellposed}, and backward because if $y(s)=0$ for some $s<0$ then uniqueness of solutions forward from $s$ would force $y\equiv0$ afterwards. Since $\xi(t)\in K_{0}$ we have $x(t)\le k$ for all $t\in\R$, so $F(x(t),\varphi_{\omega}(t))\le\Phi(\varphi_{\omega}(t))$, and for every $T>0$
	\[
	\log y(0)-\log y(-T)
	=\int_{-T}^{0}\Bigl[-\gamma+F\bigl(x(t),\varphi_{\omega}(t)\bigr)\Bigr]dt
	\ \le\ \int_{-T}^{0}\Bigl[-\gamma+\Phi\bigl(\varphi_{\omega}(t)\bigr)\Bigr]dt .
	\]
After division by $T$, the backward statement of \cref{lem:ergodic} applied to $g=\Phi$ gives as $T\to+\infty$
	\[
	\limsup_{T\to+\infty}\frac1T\bigl[\log y(0)-\log y(-T)\bigr]\ \le\ \lambda<0 .
	\]
Hence there is $T_{0}>0$ such that $y(-T)\ge y(0)e^{|\lambda|T/2}$ for $T\ge T_{0}$, which contradicts $y(-T)\le c\Lambda$, valid because $\xi(-T)\in K_{0}$. Therefore $y(0)=0$, so $y\equiv0$ along $\xi$, and $\xi$ is a bounded complete solution of the logistic equation, whence $v_{0}\in[0,k]\times\{0\}$.
\end{proof}

We now turn to the forward picture. Under the closed-form criterion \eqref{eq:Eprime} the sets $\mathcal U$ and $\mathcal A(\omega)$ coincide and the convergence is forward, uniform in the realisation and exponential.

\begin{Theorem}[The uniform attracting set under the worst-case criterion]
\label{thm:Uext} Assume \textbf{(H1)}--\textbf{(H3)} and \eqref{eq:Eprime}. Then there exist $\delta>0$ and, for each bounded $B\subset\Z$, a time $T'_{B}<\infty$ and a constant $C_{B}<\infty$, all three deterministic, such that
	\begin{equation}\label{eq:unifexp}
		y(t;0,\omega,v_{0})\ \le\ C_{B}\,e^{-\delta(t-T'_{B})}
		\qquad\text{for all }t\ge T'_{B},\ v_{0}\in B,\ \omega\in\Omega_{0}.
	\end{equation}
Consequently,
	\begin{equation}\label{eq:Uext}
		\mathcal U=\mathcal A(\omega)=[0,k]\times\{0\}\qquad\text{for every }
		\omega\in\Omega_{0},
	\end{equation}
and the segment $[0,k]\times\{0\}$ attracts every bounded subset of $\Z$ forward in time, uniformly in $\omega$.
\end{Theorem}

\begin{proof}
The map $\eta\mapsto\Phi_{\eta}(\pmax)$ is continuous and increasing with value $\Phi(\pmax)<\gamma$ at $\eta=0$, so we may fix $\eta>0$ with $\delta:=\gamma-\Phi_{\eta}(\pmax)>0$. Let $B\subset\Z$ be bounded and let $T'_{B}:=T_{B}(\eta)$ be the deterministic absorption time of \eqref{eq:absorb}, so that $x(t)\le k+\eta$ for $t\ge T'_{B}$, every $v_{0}\in B$ and every $\omega\in\Omega_{0}$. By \cref{lem:Fprops}(i) and \eqref{eq:phibounds},
	\[
	F\bigl(x(t),\varphi_{\omega}(t)\bigr)\le\Phi_{\eta}\bigl(\varphi_{\omega}(t)\bigr)
	\le\Phi_{\eta}(\pmax)=\gamma-\delta,\qquad t\ge T'_{B},
	\]
so that $\dot y\le-\delta y$ on $[T'_{B},+\infty)$, and \eqref{eq:unifexp} follows with $C_{B}:=(\sup_{B}y_{0})e^{(c\beta/a)T'_{B}}$, which bounds $y(T'_{B})$ by \eqref{eq:apriori}.

For \eqref{eq:Uext}, \cref{lem:ordering} and \cref{thm:attractorext} give $\mathcal A(\omega)=[0,k]\times\{0\}$, and \eqref{eq:sandwich} gives $[0,k]\times\{0\}\subset\mathcal U$. Conversely, let $u=(x_{\sharp},y_{\sharp})\in\mathcal U$ and let $t_{n}\to+\infty$, $\omega_{n}\in\Omega_{0}$ and $v_{n}\in K_{1}$ be such that $\Psi(t_{n},\omega_{n})v_{n}\to u$. By \eqref{eq:unifexp} with $B=K_{1}$, the second components satisfy $y(t_{n};0,\omega_{n},v_{n})\le C_{K_{1}}e^{-\delta(t_{n}-T'_{K_{1}})}\to0$, whence $y_{\sharp}=0$, and $x_{\sharp}\le k$ because $\mathcal U\subset K_{0}$. The final assertion is \eqref{eq:unifattract} for $C=\mathcal U$.
\end{proof}

The next theorem shows that the coincidence \eqref{eq:Uext} is tied to the worst-case criterion, and that in the regime $\lambda<0<\lambda_{+}$ the forward and the pullback descriptions differ.

\begin{Theorem}[The two notions differ]\label{thm:Ustrict}
Assume \textbf{(H1)}--\textbf{(H3)} and
	\begin{equation}\label{eq:strictregime}
		\lambda<0<\lambda_{+}.
	\end{equation}
Then, for every $\omega\in\Omega_{0}$,
	\begin{equation}\label{eq:strict}
		\mathcal A(\omega)=[0,k]\times\{0\}\ \subsetneq\ \mathcal U ,
	\end{equation}
and $\mathcal U$ contains the interior point $v^{*}_{\pmax}$ of \eqref{eq:vplus}. In particular, the first inclusion in \eqref{eq:sandwich} is strict and $\mathcal U\ne\overline{\bigcup_{\omega\in\Omega_{0}}\mathcal A(\omega)}$. Moreover, for every $\omega\in\Omega_{0}$,
	\begin{equation}\label{eq:Uforward}
		\lim_{t\to+\infty}\dist\bigl(\Psi(t,\omega)\,\mathcal U,\ [0,k]\times\{0\}\bigr)=0 ,
	\end{equation}
and consequently there is $t_{\omega}\ge0$ such that
	\begin{equation}\label{eq:notinvariant}
		\Psi(t,\omega)\,\mathcal U\ \subsetneq\ \mathcal U\qquad\text{for all }t\ge t_{\omega},
	\end{equation}
so that $\mathcal U$ is invariant under no individual realisation, although it satisfies \eqref{eq:familyinv}.
\end{Theorem}

\begin{proof}
Since $\lambda<0$, \cref{thm:attractorext} gives $\mathcal A(\omega)=[0,k]\times\{0\}$ for every $\omega\in\Omega_{0}$, whence also $\overline{\bigcup_{\omega\in\Omega_{0}}\mathcal A(\omega)}=[0,k]\times\{0\}$. Since $\lambda_{+}>0$, \cref{cor:frozen} gives $v^{*}_{\pmax}\in\mathcal U$ with $y^{*}_{\pmax}>0$, which is \eqref{eq:strict}. For \eqref{eq:Uforward}, let $\omega\in\Omega_{0}$ and $v_{0}=(x_{0},y_{0})\in\mathcal U$. Since $\mathcal U\subset K_{0}$ and $K_{0}$ is positively invariant, $x(t)\le k$ and hence $F(x(t),\varphi_{\omega}(t))\le\Phi(\varphi_{\omega}(t))$ for all $t\ge0$, while $y_{0}\le c\Lambda$, so \eqref{eq:representation} gives
	\[
	\sup_{v_{0}\in\mathcal U}\ y(t;0,\omega,v_{0})\ \le\ c\Lambda\,
	\exp\Bigl(\int_{0}^{t}\bigl[-\gamma+\Phi\bigl(\varphi_{\omega}(s)\bigr)\bigr]ds\Bigr)
	\ \xrightarrow[t\to+\infty]{}\ 0 ,
	\]
because the exponent divided by $t$ converges to $\lambda<0$ by \cref{lem:ergodic}. Since $0\le x(t)\le k$ along these solutions, the distance from $\Psi(t,\omega)v_{0}$ to the segment $[0,k]\times\{0\}$ is $y(t;0,\omega,v_{0})$, and \eqref{eq:Uforward} follows. Hence there is $t_{\omega}$ such that $\Psi(t,\omega)\mathcal U\subset\{y<y^{*}_{\pmax}\}$ for $t\ge t_{\omega}$, so $v^{*}_{\pmax}\notin\Psi(t,\omega)\mathcal U$, and the inclusion $\Psi(t,\omega)\mathcal U\subset\mathcal U$ of \cref{cor:invariance} is strict. \end{proof}

In the regime \eqref{eq:strictregime} the forward motion on $\mathcal U$ is thus attracted, realisation by realisation, by the pullback attractor. The convergence in \eqref{eq:Uforward} is not uniform in $\omega$, and in particular the time after which $\Psi(t,\omega)\mathcal U\subset\{y<y^{*}_{\pmax}/2\}$ cannot be chosen independently of $\omega$. Indeed, given $T>0$, the last assertion of \cref{lem:longexc} provides $\omega\in\Omega_{0}$ with $\sup_{[0,T]}|\varphi_{\omega}-\pmax|<\varrho$ for any prescribed $\varrho>0$, and \eqref{eq:contdep}, applied in $K_{1}$ to $\Psi(\cdot,\omega)v^{*}_{\pmax}$ and to the constant solution $v^{*}_{\pmax}$ of the frozen system with accessible fraction $\pmax$, gives $\sup_{[0,T]}|\Psi(s,\omega)v^{*}_{\pmax}-v^{*}_{\pmax}|\le L_{T,K_{1}}\varrho$, which is smaller than $y^{*}_{\pmax}/2$ for $\varrho$ small enough depending on $T$. For fixed values of the remaining parameters, the regime \eqref{eq:strictregime} consists exactly of the predator mortalities $\gamma\in(\E[\Phi(\psi(z^{*}))],\Phi(\pmax))$, an interval of positive length by \eqref{eq:lambdaorder}, since by \cref{def:lambda} the two conditions in \eqref{eq:strictregime} read $\gamma>\E[\Phi(\psi(z^{*}))]$ and $\gamma<\Phi(\pmax)$. In that regime the forward $\omega$-limit set of $\mathcal U$ under each realisation is the segment, by \eqref{eq:Uforward} and the invariance of the segment, whereas $\mathcal U$ itself is the forward $\omega$-limit set of the system uniformly over the realisations, which is the situation described in \cite{KloedenYang2015,CuiKloedenYang2020} for forward $\omega$-limit sets of nonautonomous systems.

\Cref{thm:Uext,thm:Ustrict} leave out the critical case $\lambda_{+}=0$. The next theorem fills it and shows that, in the extinction regime, the condition $\lambda_{+}\le0$, which is the worst-case criterion \eqref{eq:Eprime} together with its critical case, is exactly the condition under which the forward description is faithful.

\begin{Theorem}[The two descriptions in the extinction regime]\label{thm:charext}
Assume \textbf{(H1)}--\textbf{(H3)} and \eqref{eq:E}. Then
	\begin{equation}\label{eq:charext}
		\overline{\bigcup_{\omega\in\Omega_{0}}\mathcal A(\omega)}=\mathcal U
		\qquad\text{if and only if}\qquad \lambda_{+}\le0 ,
	\end{equation}
and in that case both sets equal $[0,k]\times\{0\}$. If $\lambda_{+}>0$ the inclusion $\overline{\bigcup_{\omega\in\Omega_{0}}\mathcal A(\omega)}\subset\mathcal U$ is strict, $\mathcal U$ contains the interior point $v^{*}_{\pmax}$, and $\mathcal U$ is invariant under no individual realisation.
\end{Theorem}

\begin{proof}
If $\lambda_{+}>0$ the assertions are those of \cref{thm:Ustrict}. Assume $\lambda_{+}\le0$. By \cref{thm:attractorext}, $\mathcal A(\omega)=[0,k]\times\{0\}$ for every $\omega\in\Omega_{0}$, and $[0,k]\times\{0\}\subset\mathcal U$ by \cref{cor:frozen}. By \cref{thm:control}, it remains to show that no bounded complete solution $\xi=(x,y)$ of \eqref{eq:control} satisfies $y(0)>0$. Suppose that one does, for some control $\varphi\in\mathcal C$. Then $y(s)>0$ for every $s\in\R$, by \eqref{eq:representation} and uniqueness, as in the proof of \cref{thm:attractorext}, and $\xi(s)\in K_{0}$ for every $s$, as observed after \eqref{eq:kernel}, so that $0\le x\le k$ and $0<y\le c\Lambda$ on $\R$. Since $F$ is nondecreasing in both variables,
	\[
	F\bigl(x(s),\varphi(s)\bigr)\le F(k,\pmax)=\Phi(\pmax)=\gamma+\lambda_{+}\le\gamma
	\qquad\text{for a.e. }s\in\R ,
	\]
hence $\log y$ is nonincreasing and, for every $T>0$,
	\begin{equation}\label{eq:charext1}
		\int_{-T}^{0}\Bigl[\gamma-F\bigl(x(s),\varphi(s)\bigr)\Bigr]ds
		=\log y(-T)-\log y(0)\ \le\ \log\frac{c\Lambda}{y(0)} .
	\end{equation}
If $\lambda_{+}<0$ the integrand is at least $-\lambda_{+}>0$ and the left-hand side of \eqref{eq:charext1} tends to $+\infty$ as $T\to\infty$, a contradiction. Let $\lambda_{+}=0$, so that $F(k,\pmax)=\gamma$. By \cref{lem:Fprops}(i), for $x\in[0,k]$,
	\[
	\gamma-F(x,\varphi)\ \ge\ F(k,\pmax)-F(x,\pmax)
	=\int_{x}^{k}\partial_{x}F(u,\pmax)\,du\ \ge\ m_{+}(k-x),
	\qquad m_{+}:=\frac{c\beta\pmax}{(1+ak\pmax)^{2}}>0 ,
	\]
and \eqref{eq:charext1} gives
	\[
	I:=\int_{-\infty}^{0}\bigl(k-x(s)\bigr)\,ds\ \le\ \frac{1}{m_{+}}\log\frac{c\Lambda}{y(0)}<\infty .
	\]
Let $\mathcal T:=\{s\le0\,:\,x(s)<k/2\}$, a measurable set of Lebesgue measure $|\mathcal T|\le2I/k$. On $(-\infty,0]\setminus\mathcal T$ we have $x\ge k/2$, hence $G(x,\varphi)\ge G(k/2,\pmin)=:g_{0}>0$ by the monotonicity of $G$, and, since $\alpha x(1-x/k)\le\alpha(k-x)$ on $[0,k]$ and $y\ge y(0)$ on $(-\infty,0]$, the first equation of \eqref{eq:control} yields
	\[
	\dot x\ \le\ \alpha(k-x)-g_{0}\,y(0)\quad\text{a.e. on }(-\infty,0]\setminus\mathcal T,
	\qquad
	\dot x\ \le\ \frac{\alpha k}{4}\quad\text{a.e. on }\mathcal T .
	\]
Integrating on $[s,0]$ for $s<0$ we obtain
	\[
	x(0)-x(s)\ \le\ \alpha I+\frac{\alpha k}{4}\,|\mathcal T|-g_{0}\,y(0)\bigl(|s|-|\mathcal T|\bigr),
	\]
so that $x(s)\ge g_{0}y(0)|s|-\alpha I-(\alpha k/4+g_{0}y(0))|\mathcal T|$, which exceeds $k$ for $|s|$ large, contradicting $x\le k$. Hence $y\equiv0$ along every bounded complete solution of \eqref{eq:control}, $\mathcal U\subset[0,k]\times\{0\}$, and \eqref{eq:charext} follows.
\end{proof}

Since $\lambda<\lambda_{+}$ by \eqref{eq:lambdaorder}, the condition $\lambda_{+}\le0$ implies \eqref{eq:E}, and \cref{thm:charext} together with \cref{cor:frozen} gives a precise form to the statement that $\mathcal U$ is governed by the upper worst-case exponent.

\begin{corollary}\label{cor:interior}
Assume \textbf{(H1)}--\textbf{(H3)}. Then $\mathcal U\cap\Zo\neq\emptyset$ if and only if $\lambda_{+}>0$, and $\mathcal U=[0,k]\times\{0\}$ if and only if $\lambda_{+}\le0$.
\end{corollary}

\begin{proof}
If $\lambda_{+}\le0$, then $\lambda<0$ and \cref{thm:charext} gives $\mathcal U=[0,k]\times\{0\}$, which does not meet $\Zo$. If $\lambda_{+}>0$, \cref{cor:frozen} gives $v^{*}_{\pmax}\in\mathcal U\cap\Zo$, so $\mathcal U\neq[0,k]\times\{0\}$.
\end{proof}

In light of \cref{thm:charext}, when the predator becomes extinct the two descriptions differ exactly when $\lambda_{+}>0$, and then the difference is the whole interior part of $\mathcal U$. The parameter set of \cref{ex:nie} lies in the regime \eqref{eq:strictregime}, in which the predator becomes extinct for every $\omega\in\Omega_{0}$ by \cref{thm:extinction} while $\mathcal U$ contains a compact set bounded away from the face $\{y=0\}$. That set is sustained by the histories along which the accessible fraction remains near $\pmax$ for an arbitrarily long time, which have positive probability on every finite horizon, by \cref{lem:longexc}, and probability zero on an infinite one. The gap between the two descriptions is the gap between $\lambda_{+}$ and $\lambda$, and this gap and the deficit $\lambda_{\rm det}-\lambda$ are of different orders in the amplitude. Since $\Phi$ is concave, its graph lies below the tangent line at $\pmax$, so that $\Phi(\pmax)-\Phi(u)\ge\Phi'(\pmax)(\pmax-u)$ for $u\in[\pmin,\pmax]$, and taking expectations at $u=\psi(z^{*})$, whose mean is $p$, gives the lower bound in
\[
\Phi'(\pmax)\,\varepsilon\ \le\ \lambda_{+}-\lambda\ \le\ \Phi(\pmax)-\Phi(\pmin)
\ \le\ 2\Phi'(\pmin)\,\varepsilon ,
\]
the upper bounds coming from \eqref{eq:lambdaorder} and the mean value theorem, whereas by \cref{thm:cost} below and $V<\varepsilon^{2}$ the cost of the fluctuation itself satisfies
\[
\lambda_{\rm det}-\lambda\ \le\ \frac{ac\beta k^{2}}{(1+ak\pmin)^{3}}\,\varepsilon^{2}.
\]
That is, the cost of the fluctuation is of second order in the amplitude, whereas the price of a guarantee valid for every admissible history rather than for almost every one is of first order.

\section{Persistence}\label{sec:persistence}

\Cref{sec:extinction} concerns the failure of the predator to invade. We now turn to the opposite conclusion, and we treat the two species separately because they behave differently. The prey persists with no hypothesis beyond the positivity of the parameters, whereas the predator persists only above the threshold, and then in two senses of unequal strength.

\subsection{The prey persists}

We use the notation
\[
x_{\flat}:=\frac{\gamma}{2c\beta\pmax},\qquad
x_{\natural}:=\min\Bigl\{x_{\flat},\ \frac{k}{4}\Bigr\},\qquad
\hat y:=\frac{\alpha}{2\beta\pmax},
\]
with $\Lambda$ as in \eqref{eq:Lambda}, and we set
\begin{equation}\label{eq:xlower}
	\underline x_{\eta}\ :=\ x_{\natural}\,
	\exp\Bigl(-\frac{2\beta\pmax\,c\,(\Lambda+\eta)}{\gamma}\Bigr),\qquad \underline x:=\underline x_{0}\ >\ 0 ,
\end{equation}
for $\eta\ge0$.

\begin{Theorem}[Uniform persistence of the prey]\label{thm:preypersist}
Assume \textbf{(H1)}--\textbf{(H3)} and let $\omega\in\Omega_{0}$. Then, for every $v_{0}=(x_{0},y_{0})\in\Z$ with $x_{0}>0$,
	\begin{equation}\label{eq:preypersist}
		\liminf_{t\to+\infty}x(t)\ \ge\ \underline x ,
	\end{equation}
with $\underline x$ given by \eqref{eq:xlower}. The constant $\underline x$ depends only on $\alpha,k,\beta,c,\gamma,m,\varepsilon$, not on the saturation parameter $a$, and it is independent of $v_{0}$ and of $\omega$. The same conclusion holds, with the same constant, for the solutions of \eqref{eq:control} with any measurable control $\varphi:[0,\infty)\to[\pmin,\pmax]$.
\end{Theorem}

\begin{proof}
Fix $\eta>0$ and write $y^{+}_{\eta}:=c(\Lambda+\eta)$. By \cref{cor:absorbing} applied to $B=\{v_{0}\}$, for $T_{0}:=T_{\{v_{0}\}}(\eta)$ we have $x(t)\le k+\eta$ and $y(t)\le y^{+}_{\eta}$ for all $t\ge T_{0}$, and $x(t)>0$ for all $t\ge0$ by \cref{thm:wellposed}. We record two elementary implications, valid for every $t\ge T_{0}$.

\smallskip
\noindent\emph{(a) If $x(t)\le x_{\flat}$ then $\dot y(t)\le-\frac{\gamma}{2}y(t)$.} Indeed, by \cref{lem:Fprops}(ii) and \eqref{eq:phibounds}, $F(x,\varphi_{\omega})\le c\beta\pmax x\le c\beta\pmax x_{\flat}=\frac{\gamma}{2}$.

\smallskip
\noindent\emph{(b) If $y(t)\le\hat y$ and $x(t)\le\frac{k}{4}$ then $\frac{\dot x}{x}(t)\ge\frac{\alpha}{4}$.} Indeed, using $1+a\varphi_{\omega}x\ge1$, \eqref{eq:phibounds}, $x\le k/4$ and the definition of $\hat y$,
	\[
	\frac{\dot x}{x}=\alpha\Bigl(1-\frac{x}{k}\Bigr)
	-\frac{\beta\varphi_{\omega}y}{1+a\varphi_{\omega}x}
	\ \ge\ \alpha\Bigl(1-\frac{x}{k}\Bigr)-\beta\pmax\hat y
	\ \ge\ \frac{3\alpha}{4}-\frac{\alpha}{2}=\frac{\alpha}{4}.
	\]

\smallskip
Consider the closed set $J:=\{t\ge T_{0}\,:\,x(t)\le x_{\natural}\}$ and let $[t_{1},t_{2}]$ be any connected component of $J$, with $t_{2}\le+\infty$ a priori. On $[t_{1},t_{2}]$ we have $x\le x_{\natural}\le x_{\flat}$, so (a) gives
	\begin{equation}\label{eq:ydecay}
		y(t)\ \le\ y(t_{1})\,e^{-\gamma(t-t_{1})/2}\ \le\ y^{+}_{\eta}\,
		e^{-\gamma(t-t_{1})/2},\qquad t\in[t_{1},t_{2}].
	\end{equation}
Since also $x\le x_{\natural}<k$ on that interval, $\alpha(1-x/k)\ge0$, and therefore $\dot x/x\ge-\beta\pmax\,y\ge-\beta\pmax\,y^{+}_{\eta}\,e^{-\gamma(t-t_{1})/2}$ there. Integrating from $t_{1}$ to $t$ and bounding the integral by its value over $[t_{1},+\infty)$,
	\begin{equation}\label{eq:noworse2}
		x(t)\ \ge\ x(t_{1})\,
		\exp\Bigl(-\frac{2\beta\pmax y^{+}_{\eta}}{\gamma}\Bigr),
		\qquad t\in[t_{1},t_{2}].
	\end{equation}
Next we show that $t_{2}<+\infty$. Put $\tau_{\eta}:=\frac{2}{\gamma}\bigl(\log\frac{y^{+}_{\eta}}{\hat y}\bigr)^{+}$. By \eqref{eq:ydecay}, $y(t)\le\hat y$ for $t\in[t_{1}+\tau_{\eta},t_{2}]$, so (b) applies and gives $x(t)\ge x(t_{1}+\tau_{\eta})\,e^{\alpha(t-t_{1}-\tau_{\eta})/4}$ there, which exceeds $x_{\natural}$ in finite time since $x(t_{1}+\tau_{\eta})>0$. Therefore $t_{2}<+\infty$, and $x(t_{2})=x_{\natural}$ by maximality.

If $T_{0}\notin J$, then every component $[t_{1},t_{2}]$ of $J$ has $t_{1}>T_{0}$, hence $x(t_{1})=x_{\natural}$, and \eqref{eq:noworse2} gives $x(t)\ge x_{\natural}\exp(-2\beta\pmax y^{+}_{\eta}/\gamma)$ on that component, while $x>x_{\natural}$ outside $J$. If $T_{0}\in J$, the same conclusion holds for all $t\ge t_{2}^{0}$, where $[T_{0},t_{2}^{0}]$ is the finite component of $J$ containing $T_{0}$. In either case
	\[
	\liminf_{t\to+\infty}x(t)\ \ge\ x_{\natural}\,
	\exp\Bigl(-\frac{2\beta\pmax\,c(\Lambda+\eta)}{\gamma}\Bigr)=\underline x_{\eta},
	\]
and letting $\eta\downarrow0$ yields \eqref{eq:preypersist}. The constants $\Lambda$, $x_{\natural}$ and $\hat y$ depend only on the parameters, so $\underline x$ is independent of $v_{0}$ and of $\omega$.
\end{proof}

The mechanism of the proof also shows that prey extinction is incompatible with a positive lower bound $y(t)\ge y_{*}>0$ for the predator at all large times. Indeed, if $x(t)\to0$ then \cref{lem:Fprops}(ii) gives $F(x(t),\varphi_{\omega}(t))\le c\beta\pmax\,x(t)\to0$, hence $\dot y\le-\frac{\gamma}{2}y$ for all large $t$, whence $y(t)\to0$ exponentially.

\subsection{The predator persists in the mean}

Under the ergodic criterion \eqref{eq:P} the predator does not die out, and its time average is bounded below by an explicit positive constant. Throughout we write
\begin{equation}\label{eq:means}
	\bar x_{t}:=\frac1t\int_{0}^{t}x(s)\,ds,
	\qquad
	\bar y_{t}:=\frac1t\int_{0}^{t}y(s)\,ds .
\end{equation}

\begin{Theorem}[Persistence in the mean]\label{thm:meanpersist}
Assume \textbf{(H1)}--\textbf{(H3)}, let $\omega\in\Omega_{0}$ and suppose \eqref{eq:P}. Then, for every $v_{0}\in\Zo$,
	\begin{equation}\label{eq:meanpersist}
		\liminf_{t\to+\infty}\ \frac1t\int_{0}^{t}y(s)\,ds
		\ \ \ge\ \ \frac{\alpha\,\lambda}{c\,k\,\beta^{2}\pmax^{2}}\ >\ 0 .
	\end{equation}
\end{Theorem}

\begin{proof}
Write $C_{0}:=ck\beta^{2}\pmax^{2}/\alpha$, so that the right-hand side of \eqref{eq:meanpersist} is $\lambda/C_{0}$.

\emph{Step 1. Boundedness of $\log x$.} Since $x_{0}>0$ we have $x(t)>0$ for every $t\ge0$ and $x(t)\le\max\{x_{0},k\}$ by \eqref{eq:apriori}. By \cref{thm:preypersist} there is $T_{1}$ with $x(t)\ge\underline x/2$ for $t\ge T_{1}$, while $\min_{[0,T_{1}]}x>0$ by continuity, so $\inf_{t\ge0}x(t)>0$. Therefore $\log x$ is bounded on $[0,\infty)$ and
	\[
	r_{t}:=\frac1t\log\frac{x(t)}{x_{0}}\ \xrightarrow[t\to+\infty]{}\ 0 .
	\]

\emph{Step 2. Control of the mean of $|x-k|$ by $\bar y_{t}$.} Dividing the first equation of \eqref{eq:sys} by $x$ and integrating on $[0,t]$,
	\[
	r_{t}=\alpha-\frac{\alpha}{k}\,\bar x_{t}
	-\frac1t\int_{0}^{t}\frac{\beta\varphi_{\omega}(s)\,y(s)}
	{1+a\varphi_{\omega}(s)\,x(s)}\,ds ,
	\]
and the last integrand is nonnegative and bounded by $\beta\pmax\,y$, so that
	\begin{equation}\label{eq:kminusxbar}
		k-\bar x_{t}\ \le\ \frac{k\beta\pmax}{\alpha}\,\bar y_{t}+\frac{k}{\alpha}\,r_{t}.
	\end{equation}
Next, from $|u|=-u+2u^{+}$ with $u=x-k$,
	\[
	\frac1t\int_{0}^{t}|x-k|
	=\bigl(k-\bar x_{t}\bigr)+e_{t},
	\qquad e_{t}:=\frac2t\int_{0}^{t}(x-k)^{+}.
	\]
By \cref{thm:dissipative}(i), for every $\vartheta>0$ there is $T_{\vartheta}$ with $x(t)\le k+\vartheta$ for $t\ge T_{\vartheta}$, so $e_{t}\le2\vartheta+\frac2t\int_{0}^{T_{\vartheta}}(x-k)^{+}$ and $\limsup_{t}e_{t}\le2\vartheta$. Hence $e_{t}\to0$, and \eqref{eq:kminusxbar} gives
	\begin{equation}\label{eq:xkbar}
		\frac1t\int_{0}^{t}\bigl|x(s)-k\bigr|\,ds
		\ \le\ \frac{k\beta\pmax}{\alpha}\,\bar y_{t}+d_{t},
		\qquad d_{t}:=\frac{k}{\alpha}\,r_{t}+e_{t}\ \xrightarrow[t\to+\infty]{}\ 0 .
	\end{equation}

\emph{Step 3. A lower bound for the growth rate of $y$.} As in \eqref{eq:unifclose}, $F(x,\varphi)\ge\Phi(\varphi)-c\beta\pmax|x-k|$. Since $y(t)>0$ for all $t$, \eqref{eq:representation} and \eqref{eq:xkbar} give
	\begin{equation}\label{eq:ygrowth}
		\frac1t\log\frac{y(t)}{y_{0}}
		=-\gamma+\frac1t\int_{0}^{t}F\bigl(x,\varphi_{\omega}\bigr)
		\ \ge\ \ell_{t}-C_{0}\,\bar y_{t}-c\beta\pmax\,d_{t},
		\qquad \ell_{t}:=-\gamma+\frac1t\int_{0}^{t}\Phi\bigl(\varphi_{\omega}(s)\bigr)ds ,
	\end{equation}
and $\ell_{t}\to\lambda$ by \cref{lem:ergodic}.

\emph{Step 4. Conclusion.} Let $\ell_{*}:=\liminf_{t\to+\infty}\bar y_{t}\ge0$ and suppose, for contradiction, that $\ell_{*}<\lambda/C_{0}$. Choose $t_{n}\to+\infty$ with $\bar y_{t_{n}}\to\ell_{*}$. By \eqref{eq:ygrowth},
	\[
	\liminf_{n\to\infty}\frac1{t_{n}}\log\frac{y(t_{n})}{y_{0}}
	\ \ge\ \lambda-C_{0}\ell_{*}\ >\ 0 ,
	\]
so $y(t_{n})\ge y_{0}\exp\bigl(\tfrac12(\lambda-C_{0}\ell_{*})t_{n}\bigr)\to+\infty$, which contradicts $\limsup_{t}y(t)\le c\Lambda$. Hence $\ell_{*}\ge\lambda/C_{0}$.
\end{proof}

\Cref{thm:extinction,thm:meanpersist} are exact complements of one another, since the sign of the same number $\lambda$ decides between them. The constant of \cref{thm:meanpersist} degenerates to $0$ as $\lambda\downarrow0$, which is the correct behaviour at the critical value, and it is conservative for the same reason as all our explicit constants, namely that the bound $F(x,\varphi)\ge\Phi(\varphi)-c\beta\pmax|x-k|$ is a worst-case one. The proof uses \cref{thm:preypersist} in an essential way, through Step~1. It is the unconditional persistence of the prey that makes $\log x$ bounded and underpins the whole argument.

Stronger notions of persistence, formulated for the Markov process $(v,z^{*})$ in terms of its occupation measures, are the subject of the theory of stochastic persistence \cite{SchreiberBenaimAtchade2011,Benaim2018,BenaimSchreiber2019,HeningNguyenChesson2021}, which covers degenerate noise and auxiliary environmental variables such as $z^{*}$; the expression \emph{persistence in the mean} is used here in the sense of \cite{LiuWang2011}. The invasion rate of that theory is defined through the ergodic measures of $(v,z^{*})$ carried by the extinction set, which here are $\delta_{(0,0)}\otimes\mathcal N(0,\sigma_{*}^{2})$ and $\delta_{(k,0)}\otimes\mathcal N(0,\sigma_{*}^{2})$ because the dynamics on the boundary $\{x=0\}\cup\{y=0\}$ are noise-free and $z^{*}$ is stationary and ergodic, and at the latter it is the number $\lambda$ of \eqref{eq:lambda}. What the present approach adds is of a different kind, namely explicit constants, exact rates and the worst-case statements of \cref{thm:Uext,thm:uniformpersist,cor:coexattractor}, which are uniform in the realisation. 

We now draw the consequence of \cref{thm:meanpersist} for the pullback attractor. Unlike the results of \cref{sec:extinction}, the argument is not pathwise, since it uses the invariance of $\PP$ under $\theta$ and the ergodicity of the base flow, and accordingly the conclusion holds on a $\theta$-invariant set of full measure that we do not identify.

\begin{Theorem}[The attractor meets the coexistence region]\label{thm:attractorpersist}
Assume \textbf{(H1)}--\textbf{(H3)} and \eqref{eq:P}, put $Y(\omega):=\max\{y\,:\,(x,y)\in\mathcal A(\omega)\}$ for $\omega\in\Omega_{0}$ and $Y:=0$ off $\Omega_{0}$. Then $Y$ is measurable and
	\begin{equation}\label{eq:meanheight}
		\E[Y]\ \ge\ \frac{\alpha\,\lambda}{c\,k\,\beta^{2}\pmax^{2}}\ >\ 0 .
	\end{equation}
Consequently, there is a $\theta$-invariant set $\Omega_{3}\subset\Omega_{0}$ of full measure such that $\mathcal A(\omega)\cap\Zo\neq\emptyset$ for every $\omega\in\Omega_{3}$.
\end{Theorem}

\begin{proof}
\emph{Step 1. Measurability.} If $D$ is a countable dense subset of $\Z$, then for $\omega\in\Omega_{0}$
	\[
	Y(\omega)=\sup_{(x_{d},y_{d})\in D}\bigl[y_{d}-\dist\bigl((x_{d},y_{d}),\mathcal A(\omega)\bigr)\bigr],
	\]
because each term is at most $Y(\omega)$, the second coordinate being $1$-Lipschitz, and the terms indexed by points of $D$ close to a maximiser approach $Y(\omega)$. Each $\omega\mapsto\dist(d,\mathcal A(\omega))$ is measurable by \cref{thm:attractor} and $\Omega_{0}$ is measurable, hence so is $Y$, which is the product of $\mathbf 1_{\Omega_{0}}$ with the supremum above.

\emph{Step 2. The lower bound \eqref{eq:meanheight}.} Fix $v_{0}=(k,1)\in\Zo$ and write $y(t,\omega):=y(t;0,\omega,v_{0})$, a jointly measurable function by \cref{prop:rds}. For every $\omega\in\Omega_{0}$ the pullback attraction of the bounded set $\{v_{0}\}$ gives $\dist(\Psi(t,\theta_{-t}\omega)v_{0},\mathcal A(\omega))\to0$, and since the second coordinate is $1$-Lipschitz,
	\[
	\limsup_{t\to+\infty}y(t,\theta_{-t}\omega)\ \le\ Y(\omega)\qquad\text{for every }
	\omega\in\Omega_{0}.
	\]
By \cref{cor:absorbing} with $B=\{v_{0}\}$ and $\eta=1$ there is a deterministic $T$ with $0\le y(t,\omega)\le c(\Lambda+1)$ for $t\ge T$ and all $\omega\in\Omega_{0}$, hence the reverse Fatou lemma applies along a sequence $t_{n}\to\infty$ with $t_{n}\ge T$ realising the upper limit below. Since $\theta_{-t}$ preserves $\PP$, the random variables $y(t,\theta_{-t}\,\cdot)$ and $y(t,\cdot)$ have the same law, and therefore
	\[
	\E[Y]\ \ge\ \limsup_{t\to+\infty}\E\bigl[y(t,\theta_{-t}\,\cdot)\bigr]
	=\limsup_{t\to+\infty}\E\bigl[y(t,\cdot)\bigr]
	\ \ge\ \liminf_{t\to+\infty}\frac1t\int_{0}^{t}\E\bigl[y(s,\cdot)\bigr]\,ds ,
	\]
the last inequality because the lower limit of the Ces\`aro means of a bounded function does not exceed the upper limit of the function itself. Let $t_{n}\to+\infty$ realise the lower limit on the right. By Fubini's theorem, Fatou's lemma along $(t_{n})$ and \cref{thm:meanpersist}, which holds for every $\omega\in\Omega_{0}$ and gives $\liminf_{t}\frac1t\int_{0}^{t}y(s,\omega)\,ds\ge\lambda/C_{0}$ with $C_{0}=ck\beta^{2}\pmax^{2}/\alpha$,
	\[
	\liminf_{t\to+\infty}\frac1t\int_{0}^{t}\E\bigl[y(s,\cdot)\bigr]\,ds
	=\lim_{n\to\infty}\E\Bigl[\frac1{t_{n}}\int_{0}^{t_{n}}y(s,\cdot)\,ds\Bigr]
	\ \ge\ \E\Bigl[\liminf_{n\to\infty}\frac1{t_{n}}\int_{0}^{t_{n}}y(s,\cdot)\,ds\Bigr]
	\ \ge\ \frac{\lambda}{C_{0}} ,
	\]
which is \eqref{eq:meanheight}.

\emph{Step 3. The invariant set.} Put $N_{Y}:=\{\omega\in\Omega_{0}\,:\,Y(\omega)=0\}=\{\omega\in\Omega_{0}\,:\,\mathcal A(\omega)\subset\{y=0\}\}$, a measurable set by Step~1. It is $\theta$-invariant. Indeed, if $Y(\omega)=0$ then $\mathcal A(\theta_{t}\omega)=\Psi(t,\omega)\mathcal A(\omega)\subset\{y=0\}$ for $t\ge0$ by the invariance of the face, and if $Y(\omega)>0$ then $\mathcal A(\theta_{t}\omega)$ contains the image of a point with positive second coordinate, which has positive second coordinate by \eqref{eq:representation}; thus $Y(\theta_{t}\omega)=0$ if and only if $Y(\omega)=0$ for $t\ge0$, and the same holds for $t<0$ by applying this to $\theta_{t}\omega$. By ergodicity $\PP(N_{Y})\in\{0,1\}$, and $\PP(N_{Y})=1$ is excluded by \eqref{eq:meanheight}. Hence $\Omega_{3}:=\Omega_{0}\setminus N_{Y}$ is $\theta$-invariant and of full measure. For $\omega\in\Omega_{3}$ the fibre $\mathcal A(\omega)$ contains a point with $y>0$, and such a point has $x>0$, because it lies on a complete trajectory in $K_{0}$ by \cref{thm:attractor}, the face $\{x=0\}$ is invariant in both time directions, and on that face $\dot y=-\gamma y$ admits no bounded complete solution with $y\ne0$. Hence $\mathcal A(\omega)$ meets $\Zo$ for every $\omega\in\Omega_{3}$.
\end{proof}

Under the worst-case criterion we show below that $\mathcal A(\omega)$ contains a second pullback attractor, the one of the compact subsets of $\Zo$, whose fibres are bounded away from both axes uniformly in $\omega$ (\cref{cor:coexattractor}).

\subsection{Uniform persistence and the coexistence attractor}

We use the term \emph{uniform persistence} in its usual sense, namely that the lower limit of each component is bounded below by a positive constant independent of the initial datum in $\Zo$ and, in the present setting, of the realisation $\omega\in\Omega_{0}$. For the general theory, and in particular for the acyclicity arguments that we avoid here in favour of explicit constants, we refer to \cite{FreedmanWaltman,HaleWaltman,SmithThieme,Thieme1993}. \Cref{thm:meanpersist} does not provide a positive lower bound for $\liminf_{t}y(t)$. For that we use the worst-case criterion \eqref{eq:Pprime}, which yields explicit constants, and we then show that the worst-case criterion is sharp for this purpose, up to the critical case $\lambda_{-}=0$. The skew-product semiflow over $\Omega_{0}$ is autonomous, and the persistence theory of skew-product and nonautonomous semiflows, for which see \cite{SmithThieme} and the references therein, applies to it in principle; what the explicit constants below add are the entrance times of \cref{lem:entrance}, which the abstract theory does not provide.

Assume \eqref{eq:Pprime}, set $\Delta:=\lambda_{-}>0$ and $\gamma_{1}:=\gamma+\frac{\Delta}{2}$. Since $F(\cdot,\pmin)$ is a strictly increasing bijection from $[0,+\infty)$ onto $[0,c\beta/a)$ by \cref{lem:Fprops} and $0<\gamma_{1}<\gamma+\Delta=F(k,\pmin)<c\beta/a$, there is a unique $x_{\dagger}\in(0,k)$ with $F(x_{\dagger},\pmin)=\gamma_{1}$, namely
\begin{equation}\label{eq:xdagger}
	x_{\dagger}=\frac{\gamma_{1}}{\pmin\,(c\beta-a\gamma_{1})}\ \in(0,k),
	\qquad c\beta-a\gamma_{1}>0 .
\end{equation}
Define further
\begin{equation}\label{eq:predconstants}
	\tilde y:=\frac{\alpha(k-x_{\dagger})}{2\beta\pmax k},
	\qquad
	r_{1}:=\frac{\alpha(k+x_{\dagger})}{2k},
	\qquad
	k_{1}:=\frac{k+x_{\dagger}}{2},
\end{equation}
and, recalling $\underline x$ from \eqref{eq:xlower} and observing that $\underline x\le x_{\natural}\le k/4<k_{1}$,
\begin{equation}\label{eq:tau}
	\tau_{*}(\chi):=\begin{cases}
		\dfrac1{r_{1}}\,\log\dfrac{x_{\dagger}\,(k_{1}-\chi)}{\chi\,(k_{1}-x_{\dagger})},
		& \text{if }\chi<x_{\dagger},\\[8pt]
		0,&\text{if }\chi\ge x_{\dagger},
	\end{cases}
	\qquad
	\tau_{*}:=\tau_{*}(\underline x),\qquad
	\underline y:=\tilde y\,e^{-\gamma\tau_{*}}>0 .
\end{equation}

\begin{Theorem}[Uniform persistence]\label{thm:uniformpersist}
Assume \textbf{(H1)}--\textbf{(H3)} and \eqref{eq:Pprime}, and let $\omega\in\Omega_{0}$. Then, for every $v_{0}\in\Zo$,
	\[
	\liminf_{t\to+\infty}x(t)\ \ge\ \underline x
	\qquad\text{and}\qquad
	\liminf_{t\to+\infty}y(t)\ \ge\ \underline y ,
	\]
with $\underline x$ and $\underline y$ given by \eqref{eq:xlower} and \eqref{eq:tau}. Both constants depend only on the parameters. In particular, the two species are uniformly persistent and therefore coexist. The same conclusions hold, with the same constants, for the solutions of \eqref{eq:control} with any measurable control $\varphi:[0,\infty)\to[\pmin,\pmax]$ such that $-\gamma+\Phi(\operatorname*{ess\,inf}\varphi)>0$, with $\pmin$ and $\pmax$ replaced by the essential bounds of $\varphi$ in all the constants.
\end{Theorem}

\begin{proof}
The first inequality is \cref{thm:preypersist}. For the second, fix $\chi\in(0,\underline x)$ and note that $\chi<\underline x\le k/4<k_{1}$. By \cref{thm:preypersist} there is $T_{1}\ge0$ with $x(t)\ge\chi$ for $t\ge T_{1}$, and $y(t)>0$ for all $t$ since $y_{0}>0$. We record two implications, valid for $t\ge T_{1}$.

\smallskip
\noindent\emph{(c) If $x(t)\ge x_{\dagger}$ then $\dot y(t)\ge\frac{\Delta}{2}y(t)$.} Indeed, by \cref{lem:Fprops}(i) and \eqref{eq:phibounds}, $F(x,\varphi_{\omega})\ge F(x_{\dagger},\pmin)=\gamma_{1}$, so $\dot y\ge y(\gamma_{1}-\gamma)=\frac{\Delta}{2}y$.

\smallskip
\noindent\emph{(d) If $y(t)\le\tilde y$ then $\dot x(t)\ge r_{1}\,x(t)\bigl(1-\frac{x(t)}{k_{1}}\bigr)$.} Indeed,
	\[
	\frac{\dot x}{x}=\alpha\Bigl(1-\frac{x}{k}\Bigr)
	-\frac{\beta\varphi_{\omega}y}{1+a\varphi_{\omega}x}
	\ \ge\ \alpha-\frac{\alpha}{k}x-\beta\pmax\tilde y
	= r_{1}-\frac{\alpha}{k}x
	=r_{1}\Bigl(1-\frac{x}{k_{1}}\Bigr),
	\]
because $k_{1}=r_{1}k/\alpha$ by \eqref{eq:predconstants}.

\smallskip
Let $[s_{1},s_{2}]$ be a connected component of the closed set $\{t\ge T_{1}\,:\,y(t)\le\tilde y\}$. On it (d) applies and $x(s_{1})\ge\chi$, so comparison with the logistic equation $\dot w=r_{1}w(1-w/k_{1})$, $w(s_{1})=\chi$, gives $x(t)\ge w(t)$ with
	\[
	w(t)=\frac{k_{1}\,\chi\,e^{r_{1}(t-s_{1})}}{k_{1}+\chi\bigl(e^{r_{1}(t-s_{1})}-1\bigr)},
	\]
an increasing function with $w(t)\ge x_{\dagger}$ for $t\ge s_{1}+\tau_{*}(\chi)$, with $\tau_{*}(\chi)$ as in \eqref{eq:tau}. Hence $x(t)\ge x_{\dagger}$ for $t\in[s_{1}+\tau_{*}(\chi),s_{2}]$, and by (c) $y$ is increasing there, with $\dot y\ge\frac{\Delta}{2}y$. In particular, $s_{2}<+\infty$, because $y$ grows exponentially and must leave the region $\{y\le\tilde y\}$. On the initial part $[s_{1},s_{1}+\tau_{*}(\chi)]$ we use only $\dot y\ge-\gamma y$. Combining the two parts,
	\[
	y(t)\ \ge\ y(s_{1})\,e^{-\gamma\tau_{*}(\chi)},\qquad t\in[s_{1},s_{2}] .
	\]
Exactly as in the proof of \cref{thm:preypersist}, every component with $s_{1}>T_{1}$ satisfies $y(s_{1})=\tilde y$, the finite component containing $T_{1}$, if any, is eventually left, and $y>\tilde y$ outside those components. Hence $\liminf_{t\to+\infty}y(t)\ge\tilde y\,e^{-\gamma\tau_{*}(\chi)}$. Finally, $\chi\mapsto\tau_{*}(\chi)$ is continuous and nonincreasing, so letting $\chi\uparrow\underline x$ gives $\liminf_{t}y(t)\ge\tilde y\,e^{-\gamma\tau_{*}}=\underline y$. The only properties of $\varphi_{\omega}$ used in this proof and in that of \cref{thm:preypersist} are the bounds \eqref{eq:phibounds} and the absorption of \cref{cor:absorbing}, which hold for every measurable control with values in a band, whence the last assertions of both theorems.
\end{proof}

\begin{Theorem}[No persistence of trajectories below the worst-case threshold]\label{thm:nouniform}
Assume \textbf{(H1)}--\textbf{(H3)} and $\lambda_{-}<0$. Then there is a $\theta$-invariant set $\Omega_{4}\subset\Omega_{0}$ of full measure such that, for every $\omega\in\Omega_{4}$ and every $v_{0}\in\Zo$,
	\[
	\liminf_{t\to+\infty}y(t)=0 .
	\]
\end{Theorem}

\begin{proof}
Since $-\gamma+\Phi(\pmin)=\lambda_{-}<0$ and $(\eta,\varphi)\mapsto\Phi_{\eta}(\varphi)$ is continuous, we may fix $\eta>0$ and $\varrho>0$ such that $-\gamma+\Phi_{\eta}(\pmin+\varrho)=:-\delta<0$. For $T\in\N$ consider the event
	\[
	Q_{T}:=\Bigl\{\omega\in\Omega_{1}\,:\ \sup_{t\in[0,T]}\bigl|\varphi_{\omega}(t)-\pmin\bigr|
	<\varrho\Bigr\},
	\]
which is measurable, because the supremum may be taken over rational times, and has positive probability by \cref{lem:longexc}. By Birkhoff's ergodic theorem applied to the indicator function of $Q_{T}$, there is a set $\Omega^{Q}_{T}$ of full measure, which may be taken $\theta$-invariant by the argument of \cref{lem:ergodic}, such that for every $\omega\in\Omega^{Q}_{T}$ the set $\{s\ge0\,:\,\theta_{s}\omega\in Q_{T}\}$ has density $\PP(Q_{T})>0$, and in particular is unbounded. Put $\Omega_{4}:=\Omega_{0}\cap\bigcap_{T\in\N}\Omega^{Q}_{T}$.

Let $\omega\in\Omega_{4}$ and $v_{0}\in\Zo$. By \cref{thm:dissipative} there is $t_{0}$ such that $x(t)\le k+\eta$ and $y(t)\le c\Lambda+1$ for $t\ge t_{0}$. Fix $T\in\N$ and pick $s\ge t_{0}$ with $\theta_{s}\omega\in Q_{T}$. Since $\varphi_{\omega}(s+r)=\varphi_{\theta_{s}\omega}(r)<\pmin+\varrho$ for $r\in[0,T]$, the monotonicity of $F$ in both variables gives, for $t\in[s,s+T]$,
	\[
	\frac{\dot y}{y}=-\gamma+F\bigl(x(t),\varphi_{\omega}(t)\bigr)
	\ \le\ -\gamma+\Phi_{\eta}(\pmin+\varrho)=-\delta ,
	\]
whence $y(s+T)\le y(s)e^{-\delta T}\le(c\Lambda+1)e^{-\delta T}$. Since such $s$ can be taken arbitrarily large, $\liminf_{t\to+\infty}y(t)\le(c\Lambda+1)e^{-\delta T}$ for every $T\in\N$, and letting $T\to\infty$ proves the claim.
\end{proof}

Thus every trajectory starting in the coexistence region comes arbitrarily close to the face $\{y=0\}$ at arbitrarily large times, so the predator is not even strongly persistent in the sense of \cite{SmithThieme}, let alone uniformly, although it remains weakly persistent when $\lambda>0$, since \cref{thm:meanpersist} gives $\limsup_{t}y(t)>0$. In particular, if $\lambda_{-}<0<\lambda$, the predator persists in the mean while every trajectory satisfies $\liminf_{t}y(t)=0$, and the criterion \eqref{eq:Pprime} in \cref{thm:uniformpersist} cannot be weakened to \eqref{eq:P}.

The explicit constants of \cref{thm:preypersist,thm:uniformpersist} are reached after times that depend on the initial datum only through upper and lower bounds on its coordinates, which allows us to obtain a pullback attractor bounded away from the axes. Let $\mathcal D_{+}$ denote the family of compact subsets of $\Zo$ and, for $\eta\in(0,\min\{\underline x,\underline y\})$, put
\begin{equation}\label{eq:Kplus}
	K^{+}_{\eta}:=K_{\eta}\cap\Bigl\{(x,y)\,:\ x\ge\underline x-\eta,\ \
	y\ge\underline y-\eta\Bigr\},
\end{equation}
a compact subset of $\Zo$.

\begin{lemma}[Uniform entrance times]\label{lem:entrance}
Assume \textbf{(H1)}--\textbf{(H3)} and \eqref{eq:Pprime}. For every $B\in\mathcal D_{+}$ and every $\eta\in(0,\min\{\underline x,\underline y\})$ there is $T^{+}_{B}(\eta)$, independent of $\omega$, such that $\Psi(t,\omega)B\subset K^{+}_{\eta}$ for all $t\ge T^{+}_{B}(\eta)$ and $\omega\in\Omega_{0}$.
\end{lemma}

\begin{proof}
Let $R_{1}:=\max_{B}x_{0}$, $R_{2}:=\max_{B}y_{0}$, $\varrho_{1}:=\min_{B}x_{0}>0$ and $\varrho_{2}:=\min_{B}y_{0}>0$. Let $\eta_{1}\in(0,\eta]$, to be fixed below, and let $T_{0}:=T_{B}(\eta_{1})$ be the absorption time of \eqref{eq:absorb}, so that $\Psi(t,\omega)B\subset K_{\eta_{1}}\subset K_{\eta}$ for $t\ge T_{0}$. We follow the proofs of \cref{thm:preypersist,thm:uniformpersist} and check that every time involved admits a bound depending only on $R_{1},R_{2},\varrho_{1},\varrho_{2},\eta_{1}$ and the parameters.

On $[0,T_{0}]$ we have $x\le\max\{R_{1},k\}$ and $y\le R_{2}e^{(c\beta/a)T_{0}}$ by \eqref{eq:apriori}, so the first equation of \eqref{eq:sys} gives $\dot x/x\ge-\alpha R_{1}/k-\beta\pmax R_{2}e^{(c\beta/a)T_{0}}$, since $\alpha(1-x/k)\ge-\alpha(R_{1}-k)^{+}/k$ whenever $x\le\max\{R_{1},k\}$, and therefore
	\[
	x(T_{0})\ \ge\ \varrho_{1}':=\varrho_{1}\exp\Bigl(-T_{0}\Bigl[\frac{\alpha R_{1}}{k}
	+\beta\pmax R_{2}e^{(c\beta/a)T_{0}}\Bigr]\Bigr)>0 .
	\]
In the proof of \cref{thm:preypersist}, with $\eta_{1}$ in place of $\eta$, the component of $J$ containing $T_{0}$, if any, satisfies $x\ge \varrho_{1}'e^{-2\beta\pmax y^{+}_{\eta_{1}}/\gamma} =:\varrho_{1}''$ by \eqref{eq:noworse2}, and on it $x$ grows at rate at least $\alpha/4$ after time $T_{0}+\tau_{\eta_{1}}$ by (b), so that it ends before
	\[
	T_{1}:=T_{0}+\tau_{\eta_{1}}+\frac4\alpha\log^{+}\frac{x_{\natural}}{\varrho_{1}''} .
	\]
From then on, every component of $J$ starts at the level $x_{\natural}$ and $x\ge\underline x_{\eta_{1}}$ on it by \eqref{eq:noworse2}, while $x>x_{\natural}$ outside $J$. Hence $x(t)\ge\underline x_{\eta_{1}}$ for all $t\ge T_{1}$, with $\underline x_{\eta_{1}}=x_{\natural}\exp(-2\beta\pmax c(\Lambda+\eta_{1})/\gamma)\uparrow\underline x$ as $\eta_{1}\downarrow0$.

Next we run the proof of \cref{thm:uniformpersist} with $\chi:=\underline x_{\eta_{1}}$ and $T_{1}$ as above, with $\tau_{*}(\chi)$ as in that proof. Since $\dot y\ge-\gamma y$, we have $y(T_{1})\ge \varrho_{2}e^{-\gamma T_{1}}$. The component of $\{t\ge T_{1}\,:\,y(t)\le\tilde y\}$ containing $T_{1}$, if any, satisfies $y\ge \varrho_{2}e^{-\gamma T_{1}}e^{-\gamma\tau_{*}(\chi)}$, and on it $y$ grows at rate at least $\Delta/2$ after time $T_{1}+\tau_{*}(\chi)$ by (c), so that it ends before
	\[
	T_{2}:=T_{1}+\tau_{*}(\chi)+\frac2\Delta\log^{+}
	\frac{\tilde y\,e^{\gamma(T_{1}+\tau_{*}(\chi))}}{\varrho_{2}} .
	\]
Hence $y(t)\ge\tilde y\,e^{-\gamma\tau_{*}(\chi)}=:\underline y_{\chi}$ for all $t\ge T_{2}$, exactly as in that proof, and $\underline y_{\chi}\uparrow\underline y$ as $\chi\uparrow\underline x$ because $\tau_{*}(\cdot)$ is continuous and nonincreasing. We now fix $\eta_{1}\in(0,\eta]$ so small that $\underline x_{\eta_{1}}\ge\underline x-\eta$ and $\underline y_{\chi}\ge\underline y-\eta$, and we set $T^{+}_{B}(\eta):=T_{2}$. All the quantities involved depend only on $R_{1},R_{2},\varrho_{1},\varrho_{2},\eta$ and the parameters, and the only property of the noise used is \eqref{eq:phibounds}, so $T^{+}_{B}(\eta)$ is independent of $\omega$.
\end{proof}

\begin{corollary}[Coexistence attractor]\label{cor:coexattractor}
Assume \textbf{(H1)}--\textbf{(H3)} and \eqref{eq:Pprime}. Then $\Psi$ possesses a $\mathcal D_{+}$-pullback attractor $\{\mathcal A^{\circ}(\omega)\}$, unique up to a $\PP$-null set among the invariant random compact subsets of $\Zo$ that pullback attract every $B\in\mathcal D_{+}$, and
	\begin{equation}\label{eq:Aplus}
		\mathcal A^{\circ}(\omega)\ \subset\ \mathcal A(\omega)\cap
		\bigl\{x\ge\underline x,\ y\ge\underline y\bigr\}
	\end{equation}
for every $\omega\in\Omega_{0}$. In particular, $\mathcal A^{\circ}(\omega)$ is bounded away from the two axes, uniformly in $\omega$.
\end{corollary}

\begin{proof}
The interior $\Zo$ is a Polish space, being an open subset of $\R^{2}$, and it is positively invariant under $\Psi$ by \cref{thm:wellposed}, so the restriction of $\Psi$ to $\Zo$ is a random dynamical system on $\Zo$, and $\mathcal D_{+}$ is the family of compact subsets of this phase space. Fix $\eta\in(0,\min\{\underline x,\underline y\})$. Since $K^{+}_{\eta}\in\mathcal D_{+}$, \cref{lem:entrance} gives a deterministic $T_{\eta}$ with $\Psi(s,\theta_{-s}\omega)K^{+}_{\eta}\subset K^{+}_{\eta}$ for $s\ge T_{\eta}$ and $\omega\in\Omega_{0}$, so that, for $t\ge T_{\eta}$, the closed sets $\overline{\bigcup_{s\ge t}\Psi(s,\theta_{-s}\omega)K^{+}_{\eta}}$ are nonempty compact subsets of $K^{+}_{\eta}$ decreasing in $t$, and we define $\mathcal A^{\circ}(\omega):=\mathcal O_{K^{+}_{\eta}}(\omega)$ by \eqref{eq:pullbackomega}, a nonempty compact subset of $K^{+}_{\eta}\subset\Zo$, and $\mathcal A^{\circ}(\omega):=K^{+}_{\eta}$ off $\Omega_{0}$. The set $K^{+}_{\eta}$ need not be positively invariant, so \eqref{eq:attractorformula} is not available for it, but the argument of \cite[Theorem 3.11]{CrauelFlandoli} applies verbatim to the family $\mathcal D_{+}$ in place of the bounded sets, with \cref{lem:entrance} in the role of the absorption, and gives the three properties of a $\mathcal D_{+}$-pullback attractor. Measurability holds because $\dist(v,\mathcal A^{\circ}(\omega))$ is the limit, as $n\to\infty$, of the infima over a countable dense subset of $K^{+}_{\eta}$ and over rational $s\ge n$ of the measurable functions of \cref{prop:rds}, the restriction to rational $s$ being legitimate because $s\mapsto\Psi(s,\theta_{-s}\omega)v'$ is continuous, the field being continuous in time and locally Lipschitz in the state. Invariance holds because $\Psi(\tau,\omega)$ is continuous and the cocycle property maps the family defining $\mathcal O_{K^{+}_{\eta}}(\omega)$ onto the family defining $\mathcal O_{K^{+}_{\eta}}(\theta_{\tau}\omega)$ up to a finite initial interval, the reverse inclusion following from the sequential compactness of $K^{+}_{\eta}$ by extracting convergent preimages. Attraction of every $B\in\mathcal D_{+}$ holds because $\Psi(s,\theta_{-s}\omega)B\subset\Psi(s-T^{+}_{B}(\eta),\theta_{-(s-T^{+}_{B}(\eta))}\omega)K^{+}_{\eta}$ for $s\ge T^{+}_{B}(\eta)$ by \cref{lem:entrance} and the cocycle property, and the Hausdorff semidistance from a decreasing family of compact sets to its intersection tends to zero. The uniqueness in the stated sense is Crauel's theorem recalled before \cref{thm:attractor}, applied on the Polish space $\Zo$, whose deterministic compact sets are exactly the members of $\mathcal D_{+}$; the measurability notion used here is equivalent to Effros measurability and does not depend on the choice of a complete metric on $\Zo$ \cite[Remark 5.9]{Crauel1999}. The inclusion \eqref{eq:Aplus} is obtained pathwise. For $\eta'\in(0,\min\{\underline x,\underline y\})$ and $t\ge T^{+}_{K^{+}_{\eta}}(\eta')$, the invariance of $\mathcal A^{\circ}$ and \cref{lem:entrance} give $\mathcal A^{\circ}(\omega)=\Psi(t,\theta_{-t}\omega)\mathcal A^{\circ}(\theta_{-t}\omega)\subset\Psi(t,\theta_{-t}\omega)K^{+}_{\eta}\subset K^{+}_{\eta'}$ for every $\omega\in\Omega_{0}$, since $\theta_{-t}\omega\in\Omega_{0}$ and the entrance time is deterministic; letting $\eta'\downarrow0$ gives $\mathcal A^{\circ}(\omega)\subset\{x\ge\underline x,\ y\ge\underline y\}$. Finally, since $\mathcal A^{\circ}(\theta_{-t}\omega)\subset K^{+}_{\eta}\subset K_{\eta}$ for every $t\ge0$, the same invariance gives $\mathcal A^{\circ}(\omega)\subset\Psi(t,\theta_{-t}\omega)K_{\eta}$ for every $t\ge0$, and \eqref{eq:attractorformula} with this $\eta$ gives $\mathcal A^{\circ}(\omega)\subset\mathcal A(\omega)$.
\end{proof}

The gap between \eqref{eq:P} and \eqref{eq:Pprime} is a genuine effect of the noise and not an artefact of the method, as \cref{thm:nouniform} shows, and it has a simple explanation. On the extinction side the relevant bound is $\limsup_{t}x(t)\le k$, which is attained along the prey-only dynamics. Replacing $F(x,\varphi_{\omega})$ by $\Phi(\varphi_{\omega})$ and averaging is therefore lossless, and the criterion \eqref{eq:E} is exact. On the persistence side the predator can be kept away from zero only if the prey stays above a level at which the predator's growth rate is positive, and whether this happens depends on the joint history of $(x(t),\varphi_{\omega}(t))$ and not on the law of $\varphi_{\omega}$ alone. As $\varepsilon\downarrow0$ the four conditions of \cref{ssec:criteria} collapse to the classical dichotomy of the deterministic model \eqref{eq:det}, decided by the sign of $\lambda_{\rm det}=\Phi(p)-\gamma$.

\section{Forward versus pullback when the predator persists}\label{sec:fibres}

In the extinction regime the gap between the two descriptions is produced by the ergodic average, which forces the predator component of every complete trajectory to vanish while the control system can sustain the frozen equilibrium $v^{*}_{\pmax}$. When $\lambda>0$ no such constraint is available, and the question becomes how large the union of the fibres of the pullback attractor is. In this section we bound it from below (\cref{thm:lowerbound}).

\subsection{The essential closure of the fibres}\label{ssec:essclosure}

The closure of the union of the fibres $\mathcal A(\omega)$ over $\omega\in\Omega_{0}$ depends on the choice of the full-measure set $\Omega_{0}$, since the fibres over a null set of realisations may lie outside the closure of the others. The intrinsic object is the following. By \cref{thm:attractor}, the map $\omega\mapsto\dist(v,\mathcal A(\omega))$ is measurable for every $v\in\Z$, and we set
\begin{equation}\label{eq:Abar}
	\overline{\mathcal A}:=\Bigl\{v\in\Z\,:\ \PP\bigl(\dist(v,\mathcal A(\omega))<r\bigr)>0
	\ \text{ for every }r>0\Bigr\},
\end{equation}
the \emph{essential closure} of the union of the fibres. Let $M_{1}$ denote the maximum of $|\mathcal F|$ on $K_{1}\times[\pmin,\pmax]$.

\begin{lemma}\label{lem:essclosure}
Assume \textbf{(H1)}--\textbf{(H3)}. Then $\overline{\mathcal A}$ is a nonempty compact subset of $K_{0}$, $\mathcal A(\omega)\subset\overline{\mathcal A}$ for $\PP$-almost every $\omega$, and
	\begin{equation}\label{eq:essclosure}
		\overline{\mathcal A}=\bigcap_{\Omega'}\ \overline{\bigcup_{\omega\in\Omega'}\mathcal A(\omega)} ,
	\end{equation}
the intersection being taken over all measurable $\theta$-invariant sets $\Omega'\subset\Omega_{0}$ of full measure. In particular $\overline{\mathcal A}\subset\mathcal U$.
\end{lemma}

\begin{proof}
If $v\notin\overline{\mathcal A}$ there is $r>0$ with $\PP(\dist(v,\mathcal A(\omega))<r)=0$, and for $|v'-v|<r/2$ the event $\{\dist(v',\mathcal A(\omega))<r/2\}$ is contained in the former, so $v'\notin\overline{\mathcal A}$. Hence $\overline{\mathcal A}$ is closed, it is nonempty because it contains almost every fibre by the next paragraph and the fibres are nonempty, and it is contained in $K_{0}$ because $\dist(v,\mathcal A(\omega))\ge\dist(v,K_{0})$ for every $\omega\in\Omega_{0}$ by \cref{thm:attractor}, so that for $v\notin K_{0}$ the event $\{\dist(v,\mathcal A(\omega))<\dist(v,K_{0})\}$ lies in the null set $\Omega\setminus\Omega_{0}$. The open set $\Z\setminus\overline{\mathcal A}$ is covered by the balls $B(v,r_{v})$, $v\notin\overline{\mathcal A}$, with $\PP(\dist(v,\mathcal A(\omega))<r_{v})=0$, and by the Lindel\"of property by countably many of them; the union $N_{*}$ of the corresponding null events is null, and for $\omega\notin N_{*}$ the fibre $\mathcal A(\omega)$ meets none of these balls, that is, $\mathcal A(\omega)\subset\overline{\mathcal A}$.

Let $\Omega'\subset\Omega_{0}$ be $\theta$-invariant and of full measure. If $v\in\overline{\mathcal A}$ and $r>0$, the event $\{\dist(v,\mathcal A(\omega))<r\}$ has positive probability, hence meets $\Omega'$, so $v$ lies in the closure of $\bigcup_{\omega\in\Omega'}\mathcal A(\omega)$. Conversely, let $v\notin\overline{\mathcal A}$ and let $r>0$ be such that $N_{0}:=\{\omega\in\Omega_{0}\,:\,\dist(v,\mathcal A(\omega))<r\}$ is null. The set
	\[
	N':=\bigl\{\omega\in\Omega_{0}\,:\ \bigl|\{s\in\R\,:\,\theta_{s}\omega\in N_{0}\}\bigr|>0\bigr\}
	\]
is measurable by Fubini's theorem, $\theta$-invariant because replacing $\omega$ by $\theta_{\tau}\omega$ shifts the variable $s$, and null, since $\int_{\Omega}|\{s\in[-n,n]\,:\,\theta_{s}\omega\in N_{0}\}|\,d\PP(\omega)=\int_{-n}^{n}\PP(\theta_{s}^{-1}N_{0})\,ds=0$ for every $n\in\N$ by the invariance of $\PP$. Put $\Omega':=\Omega_{0}\setminus N'$, a $\theta$-invariant set of full measure. For $\omega\in\Omega_{0}$ and $s\in\R$ the invariance of the pullback attractor and the bound $|\mathcal F|\le M_{1}$ on $K_{1}\times[\pmin,\pmax]$ give $\dist(v,\mathcal A(\theta_{s}\omega))\le\dist(v,\mathcal A(\omega))+M_{1}|s|$, because every point of $\mathcal A(\omega)$ has its image in $\mathcal A(\theta_{s}\omega)$ when $s\ge0$ and is the image of a point of $\mathcal A(\theta_{s}\omega)$ when $s<0$, under a solution map over a time interval of length $|s|$ along which the solution stays in $K_{0}$. Hence, if $\dist(v,\mathcal A(\omega))<r$, then $\theta_{s}\omega\in N_{0}$ for all $s$ in a nondegenerate interval and $\omega\in N'$. Therefore $\dist(v,\mathcal A(\omega))\ge r$ for every $\omega\in\Omega'$, and $v$ does not belong to the closure of $\bigcup_{\omega\in\Omega'}\mathcal A(\omega)$. This proves \eqref{eq:essclosure}, and the last assertion follows from \eqref{eq:sandwich} with $\Omega'=\Omega_{0}$, $\mathcal U$ being closed.
\end{proof}

In the extinction regime every fibre over $\Omega_{0}$ is the segment $[0,k]\times\{0\}$ by \cref{thm:attractorext}, so that $\overline{\mathcal A}$ coincides with the closure of the union of the fibres and \cref{thm:charext} holds verbatim with $\overline{\mathcal A}$ in place of that closure. The two comparisons thus bear on the same object.

\subsection{The past and the future of a realisation}\label{ssec:pastfuture}

The pullback attractor at $\omega$ depends only on the past of $\omega$, whereas the shadowing of \cref{lem:longexc} is, once $z^{*}(\omega)$ is given, a statement about its future, and the two can be prescribed independently. We make this precise. For $\omega\in\Omega$ let $\omega_{-}$ and $\omega_{+}$ denote the restrictions of $\omega$ to $(-\infty,0]$ and to $[0,\infty)$, and let $\mathscr F_{-}$ and $\mathscr F_{+}$ be the $\sigma$-algebras they generate. Since $\mathscr F$ is generated by the evaluations and the increments of $\omega$ on $(-\infty,0]$ and on $[0,\infty)$ are independent under the two-sided Wiener measure, the map $\omega\mapsto(\omega_{-},\omega_{+})$ identifies $(\Omega,\mathscr F,\PP)$ with the product of the two marginal spaces $(\Omega_{-},\mathscr F_{-},\PP_{-})$ and $(\Omega_{+},\mathscr F_{+},\PP_{+})$; in particular an event $B\in\mathscr F_{+}$ is of the form $\{\omega\,:\,\omega_{+}\in B_{+}\}$ with $\PP_{+}(B_{+})=\PP(B)$. For a measurable $S\subset\Omega$ and $\omega\in\Omega$ we consider the future section of $S$ at $\omega$ and its measure,
\begin{equation}\label{eq:section}
	S^{+}(\omega):=\bigl\{\omega'\in S\,:\,\omega'_{-}=\omega_{-}\bigr\},
	\qquad
	\mu_{S}(\omega):=\PP_{+}\bigl(\{\omega'_{+}\,:\,\omega'\in S^{+}(\omega)\}\bigr).
\end{equation}
By Fubini's theorem $\mu_{S}$ is $\mathscr F_{-}$-measurable and $\int_{\Omega}\mu_{S}\,d\PP=\PP(S)$; in particular $\mu_{S}=1$ almost surely whenever $\PP(S)=1$.

\begin{lemma}\label{lem:pastmeas}
Let $\omega,\omega'\in\Omega_{1}$ with $\omega'_{-}=\omega_{-}$. Then $\varphi_{\omega'}(t)=\varphi_{\omega}(t)$ for every $t\le0$, and $\mathcal A(\omega')=\mathcal A(\omega)$ whenever $\omega,\omega'\in\Omega_{0}$. Moreover, for $t\ge0$,
	\begin{equation}\label{eq:futureOU}
		z^{*}(\theta_{t}\omega')=e^{-\kappa t}z^{*}(\omega)
		+\sigma\Bigl(\omega'(t)-\kappa\!\int_{0}^{t}\!e^{-\kappa(t-u)}\omega'(u)\,du\Bigr).
	\end{equation}
\end{lemma}

\begin{proof}
For $s\le0$ the integral \eqref{eq:OU} defining $z^{*}(\theta_{s}\omega)$ involves $(\theta_{s}\omega)(r)=\omega(r+s)-\omega(s)$ for $r\le0$ only, that is, the values of $\omega$ on $(-\infty,s]\subset(-\infty,0]$, so $z^{*}(\theta_{s}\omega')=z^{*}(\theta_{s}\omega)$ and $\varphi_{\omega'}(s)=\varphi_{\omega}(s)$. In \eqref{eq:attractorformula} the map $\Psi(t,\theta_{-t}\omega)$ is the solution operator of \eqref{eq:sys} on $[-t,0]$, which depends on $\varphi_{\omega}$ on $[-t,0]$ only, whence $\mathcal A(\omega')=\mathcal A(\omega)$. Finally, \eqref{eq:futureOU} is \eqref{eq:varconst} written for $\omega'$, with $z^{*}(\omega')=z^{*}(\omega)$.
\end{proof}

\begin{lemma}\label{lem:goodset}
Let $\Omega'\subset\Omega_{0}$ be $\theta$-invariant and of full measure. There is a $\theta$-invariant set $\Omega'^{\sharp}\subset\Omega'$ with $\PP(\Omega'^{\sharp})=1$ such that, for every $\omega\in\Omega'^{\sharp}$,
	\begin{equation}\label{eq:goodset}
		\mu_{\Omega'^{\sharp}}(\theta_{s}\omega)=1\qquad\text{for Lebesgue-a.e. }s\in\R .
	\end{equation}
\end{lemma}

\begin{proof}
Put $\Omega^{(0)}:=\Omega'$ and, for $n\ge0$,
	\[
	\Omega^{(n+1)}:=\Bigl\{\omega\in\Omega^{(n)}\,:\ \mu_{\Omega^{(n)}}(\theta_{s}\omega)=1
	\ \text{for Lebesgue-a.e. }s\in\R\Bigr\} .
	\]
We show by induction that $\Omega^{(n)}$ is measurable, $\theta$-invariant and of full measure. This holds for $n=0$. Assume it for $n$. The map $(s,\omega)\mapsto\mu_{\Omega^{(n)}}(\theta_{s}\omega)$ is jointly measurable, being the composition of the measurable flow with a measurable function, so $\Omega^{(n+1)}$ is measurable by Fubini's theorem, and it is $\theta$-invariant because replacing $\omega$ by $\theta_{r}\omega$ shifts the variable $s$ by $r$ and $\Omega^{(n)}$ is invariant. Since $\int\mu_{\Omega^{(n)}}\,d\PP=\PP(\Omega^{(n)})=1$ and $\mu_{\Omega^{(n)}}\le1$, we have $\mu_{\Omega^{(n)}}=1$ almost surely, and the invariance of $\PP$ under each $\theta_{s}$ together with Fubini's theorem gives, for every $N\in\N$,
	\[
	\int_{\Omega}\Bigl|\bigl\{s\in[-N,N]\,:\,\mu_{\Omega^{(n)}}(\theta_{s}\omega)<1\bigr\}\Bigr|\,d\PP(\omega)
	=\int_{-N}^{N}\PP\bigl(\mu_{\Omega^{(n)}}\circ\theta_{s}<1\bigr)\,ds=0 ,
	\]
so that $\PP(\Omega^{(n+1)})=1$. Let $\Omega'^{\sharp}:=\bigcap_{n\ge0}\Omega^{(n)}$, a $\theta$-invariant set of full measure contained in $\Omega'$. If $\omega\in\Omega'^{\sharp}$, then for every $n$ the set of $s$ with $\mu_{\Omega^{(n)}}(\theta_{s}\omega)<1$ is Lebesgue-null, hence so is the union of these sets over $n$, and for every $s$ outside this union the sets $\bigl(\Omega^{(n)}\bigr)^{+}(\theta_{s}\omega)$ have full $\PP_{+}$-measure for all $n$. Since they decrease to $\bigl(\Omega'^{\sharp}\bigr)^{+}(\theta_{s}\omega)$, continuity from above of $\PP_{+}$ gives \eqref{eq:goodset}. The iteration is needed because a single step would only relate the future sections of $\Omega'$ to the shifts of a smaller set, whereas \eqref{eq:goodset} refers to $\Omega'^{\sharp}$ itself.
\end{proof}

The next lemma is the version of \cref{lem:longexc} in which the past of the realisation is prescribed.

\begin{lemma}[Shadowing with prescribed past]\label{lem:shadowpast}
Let $\Omega''\subset\Omega_{0}$ be measurable, let $\omega\in\Omega''$ satisfy $\mu_{\Omega''}(\omega)=1$, let $T>0$, $\varrho>0$, and let $\varphi_{0}:[0,T]\to[\pmin,\pmax]$ be continuous with $\varphi_{0}(0)=\varphi_{\omega}(0)$. Then there is $\omega'\in\Omega''$ with $\omega'_{-}=\omega_{-}$ and
	\begin{equation}\label{eq:shadowpast}
		\sup_{t\in[0,T]}\bigl|\varphi_{\omega'}(t)-\varphi_{0}(t)\bigr|<\varrho .
	\end{equation}
\end{lemma}

\begin{proof}
Since $\varphi_{\omega}(0)=\psi(z^{*}(\omega))$ lies in the open interval $(\pmin,\pmax)$, for $\delta\in(0,1)$ with $2\delta\varepsilon<\varrho/2$ the function
	\[
	\tilde\varphi_{0}(t):=\varphi_{0}(0)+(1-\delta)\bigl(\varphi_{0}(t)-\varphi_{0}(0)\bigr),
	\qquad t\in[0,T],
	\]
is continuous, takes values in $(\pmin,\pmax)$, satisfies $\tilde\varphi_{0}(0)=\varphi_{0}(0)$ and $\sup_{[0,T]}|\tilde\varphi_{0}-\varphi_{0}|\le2\delta\varepsilon<\varrho/2$. Put $M_{0}:=\psi^{-1}\circ\tilde\varphi_{0}$, a continuous function with $M_{0}(0)=z^{*}(\omega)$, and let $\varrho'>0$ be such that $|z-M_{0}(t)|<\varrho'$ implies $|\psi(z)-\tilde\varphi_{0}(t)|<\varrho/2$ for all $t\in[0,T]$. Define $r$, $h$ and $\delta_{0}$ from $M_{0}$ exactly as in the proof of \cref{lem:longexc}, and let $B$ be the event defined there, which belongs to $\mathscr F_{+}$ and has positive probability. Since $B\in\mathscr F_{+}$, its future section at $\omega$ has $\PP_{+}$-measure $\PP(B)>0$, and since $\mu_{\Omega''}(\omega)=1$ the two sections intersect, so there is $\omega'\in\Omega''$ with $\omega'_{-}=\omega_{-}$ and $\omega'\in B$. By \eqref{eq:futureOU} and \eqref{eq:volterra}, writing $e:=\omega'-h$ on $[0,T]$,
	\[
	z^{*}(\theta_{t}\omega')-M_{0}(t)
	=e^{-\kappa t}\bigl(z^{*}(\omega)-M_{0}(0)\bigr)
	+\sigma\Bigl(e(t)-\kappa\!\int_{0}^{t}\!e^{-\kappa(t-u)}e(u)\,du\Bigr)
	=\sigma\Bigl(e(t)-\kappa\!\int_{0}^{t}\!e^{-\kappa(t-u)}e(u)\,du\Bigr),
	\]
whose absolute value is at most $2\sigma\delta_{0}=\varrho'/2$. Hence $|\varphi_{\omega'}(t)-\tilde\varphi_{0}(t)|<\varrho/2$ on $[0,T]$, and \eqref{eq:shadowpast} follows.
\end{proof}

We shall also need the continuous dependence of the solutions of \eqref{eq:control} on the control in the $L^{1}$ norm, which is proved exactly as \eqref{eq:contdep}. Let $L$ and $L'$ denote constants such that $|\mathcal F(v,\varphi)-\mathcal F(v',\varphi')|\le L|v-v'|+L'|\varphi-\varphi'|$ on $K_{1}\times[\pmin,\pmax]$, as in the proof of \cref{prop:rds}.

\begin{lemma}\label{lem:L1dep}
Let $t>0$, let $\varphi,\varphi':[0,t]\to[\pmin,\pmax]$ be measurable, and let $v$, $v'$ be the solutions of \eqref{eq:control} with controls $\varphi$, $\varphi'$ starting from $v_{0},v_{0}'\in K_{1}$ at time $0$. Then
	\begin{equation}\label{eq:L1dep}
		\sup_{s\in[0,t]}\bigl|v(s)-v'(s)\bigr|
		\ \le\ e^{Lt}\Bigl(|v_{0}-v_{0}'|+L'\!\int_{0}^{t}\bigl|\varphi(s)-\varphi'(s)\bigr|\,ds\Bigr),
	\end{equation}
with $L$ and $L'$ as above.
\end{lemma}

\begin{proof}
Both solutions remain in $K_{1}$ by \cref{cor:absorbing}. Subtracting the integral forms of the two equations and using $|\mathcal F(v,\varphi)-\mathcal F(v',\varphi')|\le L|v-v'|+L'|\varphi-\varphi'|$ on $K_{1}\times[\pmin,\pmax]$, we obtain $|v(s)-v'(s)|\le|v_{0}-v_{0}'|+L'\int_{0}^{t}|\varphi-\varphi'|+L\int_{0}^{s}|v-v'|$, and Gronwall's inequality gives \eqref{eq:L1dep}.
\end{proof}

\subsection{The lower bound}\label{ssec:lowerbound}

For $\bar\varphi\in[\pmin,\pmax]$ with $\Phi(\bar\varphi)>\gamma$ the frozen system has the interior equilibrium
\begin{equation}\label{eq:vfrozen}
	v^{*}_{\bar\varphi}:=\bigl(x^{*}_{\bar\varphi},y^{*}_{\bar\varphi}\bigr),
	\qquad
	x^{*}_{\bar\varphi}=\frac{\gamma}{\bar\varphi(c\beta-a\gamma)},
	\qquad
	y^{*}_{\bar\varphi}=\frac{\alpha}{\beta\bar\varphi}\Bigl(1-\frac{x^{*}_{\bar\varphi}}{k}\Bigr)\bigl(1+a\bar\varphi x^{*}_{\bar\varphi}\bigr),
\end{equation}
which is \eqref{eq:detequil} with $p$ replaced by $\bar\varphi$ and reduces to \eqref{eq:vplus} for $\bar\varphi=\pmax$, and it is the only equilibrium of $S_{\bar\varphi}$ in $\Zo$. We say that the frozen system with accessible fraction $\bar\varphi$ is \emph{globally stable}, and write $\bar\varphi\in\Sigma_{\rm gs}$, if $\Phi(\bar\varphi)>\gamma$ and $v^{*}_{\bar\varphi}$ is asymptotically stable for the semiflow $S_{\bar\varphi}$ and attracts every point of $\Zo$. The prey isocline of the frozen system is the graph of $y=\alpha(1-x/k)(1+a\bar\varphi x)/(\beta\bar\varphi)$, a parabola whose maximum is attained at $x_{\rm v}(\bar\varphi):=\tfrac12\bigl(k-1/(a\bar\varphi)\bigr)$, and by the geometric criterion of Cheng, Hsu and Lin \cite[Theorem 1]{ChengHsuLin1981}, which applies to the Holling type~II response, $\bar\varphi\in\Sigma_{\rm gs}$ whenever $\Phi(\bar\varphi)>\gamma$ and $x^{*}_{\bar\varphi}>x_{\rm v}(\bar\varphi)$. After multiplication by $\bar\varphi>0$ the two conditions $x^{*}_{\bar\varphi}<k$ and $x^{*}_{\bar\varphi}>x_{\rm v}(\bar\varphi)$ become linear in $\bar\varphi$, which gives
\begin{equation}\label{eq:Sigmags}
	\Sigma_{\rm gs}\ \supset\ \Bigl(\frac{\gamma}{k(c\beta-a\gamma)},\ \frac1k\Bigl(\frac{2\gamma}{c\beta-a\gamma}+\frac1a\Bigr)\Bigr)\cap[\pmin,\pmax]
	\qquad\text{whenever }c\beta>a\gamma ,
\end{equation}
so that $\Sigma_{\rm gs}$ is nonempty whenever $\lambda_{+}>0$ and $\pmin<\frac1k\bigl(\frac{2\gamma}{c\beta-a\gamma}+\frac1a\bigr)$, and it is empty when $\lambda_{+}\le0$; in particular it contains every $\bar\varphi$ with $\Phi(\bar\varphi)>\gamma$ when $ak\bar\varphi\le1$, since then $x_{\rm v}(\bar\varphi)\le0$.

\begin{Theorem}[Lower bound for the union of the fibres]\label{thm:lowerbound}
Assume \textbf{(H1)}--\textbf{(H3)} and let $\overline{\mathcal A}$ be the essential closure \eqref{eq:Abar}. Then the following properties hold.
\begin{enumerate}
\item $S_{\varphi}(t,0)\,\overline{\mathcal A}\subset\overline{\mathcal A}$ for every $\varphi\in\mathcal C$ and $t\ge0$; in particular $\Psi(t,\omega)\overline{\mathcal A}\subset\overline{\mathcal A}$ for every $\omega\in\Omega_{0}$.
\item If $\lambda>0$, then $\overline{\mathcal A}\cap\Zo\ne\emptyset$, and for every $u_{0}\in\overline{\mathcal A}\cap\Zo$ and every $\bar\varphi\in[\pmin,\pmax]$ with $\Phi(\bar\varphi)>\gamma$ the $\omega$-limit set $\mathcal L_{\bar\varphi}(u_{0})$ of $u_{0}$ under the frozen semiflow $S_{\bar\varphi}$ is contained in $\overline{\mathcal A}\cap\Zo$, and it is either the equilibrium $v^{*}_{\bar\varphi}$ or a periodic orbit of $S_{\bar\varphi}$. In particular $v^{*}_{\bar\varphi}\in\overline{\mathcal A}$ for every $\bar\varphi\in\Sigma_{\rm gs}$.
\item If $\lambda>0$ and $u_{0}\in\overline{\mathcal A}\cap\Zo$, let
		\begin{equation}\label{eq:Edef}
			\mathcal E_{u_{0}}:=\bigl([0,k]\times\{0\}\bigr)\cup\bigcup_{\bar\varphi\in[\pmin,\pmax],\ \Phi(\bar\varphi)>\gamma}\mathcal L_{\bar\varphi}(u_{0}),
			\qquad
			\mathcal R_{u_{0}}:=\bigl\{S_{\varphi}(t,0)e\,:\,e\in\mathcal E_{u_{0}},\ \varphi\in\mathcal C,\ t\ge0\bigr\}
		\end{equation}
		be the union of the segment with these limit sets, the set of $\bar\varphi$ in the union being the nonempty interval $(\gamma/(k(c\beta-a\gamma)),\pmax]\cap[\pmin,\pmax]$ because $\lambda>0$ forces $\Phi(\pmax)>\gamma$, and the set of states reachable from $\mathcal E_{u_{0}}$ under the control system. Then
		\begin{equation}\label{eq:lowerbound}
			\overline{\mathcal R_{u_{0}}}\ \subset\ \overline{\mathcal A}\ \subset\ \mathcal U .
		\end{equation}
\item $\overline{\mathcal A}=\overline{\bigcup_{\omega\in\Omega_{0}}\Psi(t,\omega)\overline{\mathcal A}}$ for every $t\ge0$.
\end{enumerate}
\end{Theorem}

\begin{proof}
(i) We first prove that, for every $\theta$-invariant $\Omega'\subset\Omega_{0}$ of full measure, the closed set
	\[
	\overline{\mathcal A}_{\Omega'^{\sharp}}:=\overline{\bigcup_{\omega\in\Omega'^{\sharp}}\mathcal A(\omega)}
	\]
satisfies $S_{\varphi}(t,0)\overline{\mathcal A}_{\Omega'^{\sharp}}\subset\overline{\mathcal A}_{\Omega'^{\sharp}}$, where $\Omega'^{\sharp}$ is the set of \cref{lem:goodset}. Let $\omega\in\Omega'^{\sharp}$, $u\in\mathcal A(\omega)$, $\varphi\in\mathcal C$ and $t>0$. By \eqref{eq:goodset} there are $s_{n}\downarrow0$ with $\mu_{\Omega'^{\sharp}}(\theta_{-s_{n}}\omega)=1$. Put $\omega_{n}:=\theta_{-s_{n}}\omega\in\Omega'^{\sharp}$. By the invariance of the pullback attractor, $\mathcal A(\omega)=\Psi(s_{n},\omega_{n})\mathcal A(\omega_{n})$, so there is $u_{n}\in\mathcal A(\omega_{n})$ with $\Psi(s_{n},\omega_{n})u_{n}=u$, and $|u_{n}-u|\le M_{1}s_{n}\to0$. Next, since continuous functions are dense in $L^{1}(0,t)$ and the projection of $\R$ onto $[\pmin,\pmax]$ is $1$-Lipschitz, there are continuous $\varphi_{n}:[0,t]\to[\pmin,\pmax]$ with $\int_{0}^{t}|\varphi_{n}-\varphi|<1/n$, and after replacing $\varphi_{n}$ on a short initial interval by the linear interpolation between $\varphi_{\omega_{n}}(0)$ and its value at the end of that interval, which changes it by less than $1/n$ in $L^{1}(0,t)$, we may assume $\varphi_{n}(0)=\varphi_{\omega_{n}}(0)$ and $\int_{0}^{t}|\varphi_{n}-\varphi|<2/n$. \Cref{lem:shadowpast}, applied with $\Omega''=\Omega'^{\sharp}$, which is measurable as a countable intersection of measurable sets, provides $\omega_{n}'\in\Omega'^{\sharp}$ with $(\omega_{n}')_{-}=(\omega_{n})_{-}$ and $\sup_{[0,t]}|\varphi_{\omega_{n}'}-\varphi_{n}|<1/n$. By \cref{lem:pastmeas}, $\mathcal A(\omega_{n}')=\mathcal A(\omega_{n})\ni u_{n}$, and by invariance $\Psi(t,\omega_{n}')u_{n}\in\mathcal A(\theta_{t}\omega_{n}')\subset\overline{\mathcal A}_{\Omega'^{\sharp}}$, since $\theta_{t}\omega_{n}'\in\Omega'^{\sharp}$. On the other hand, $\Psi(t,\omega_{n}')u_{n}$ is the solution of \eqref{eq:control} with control $\varphi_{\omega_{n}'}$ starting from $u_{n}$, and $\int_{0}^{t}|\varphi_{\omega_{n}'}-\varphi|<(t+2)/n$, so \cref{lem:L1dep} gives
	\[
	\bigl|\Psi(t,\omega_{n}')u_{n}-S_{\varphi}(t,0)u\bigr|
	\ \le\ e^{Lt}\Bigl(M_{1}s_{n}+\frac{L'(t+2)}{n}\Bigr)\ \xrightarrow[n\to\infty]{}\ 0 .
	\]
Since $\overline{\mathcal A}_{\Omega'^{\sharp}}$ is closed, $S_{\varphi}(t,0)u\in\overline{\mathcal A}_{\Omega'^{\sharp}}$, and the same holds for every $u\in\overline{\mathcal A}_{\Omega'^{\sharp}}$ by \eqref{eq:L1dep} with $\varphi'=\varphi$. Now let $u\in\overline{\mathcal A}$. By \eqref{eq:essclosure}, $u\in\overline{\mathcal A}_{\Omega'^{\sharp}}$ for every $\Omega'$, hence $S_{\varphi}(t,0)u\in\overline{\mathcal A}_{\Omega'^{\sharp}}\subset\overline{\bigcup_{\omega\in\Omega'}\mathcal A(\omega)}$ for every $\theta$-invariant $\Omega'\subset\Omega_{0}$ of full measure, and \eqref{eq:essclosure} gives $S_{\varphi}(t,0)u\in\overline{\mathcal A}$. The last assertion is the case $\varphi=\varphi_{\omega}$.

(ii) Let $\lambda>0$. By \cref{thm:attractorpersist} there is a set of full measure on which $\mathcal A(\omega)\cap\Zo\neq\emptyset$, and by \cref{lem:essclosure} there is a set of full measure on which $\mathcal A(\omega)\subset\overline{\mathcal A}$; for $\omega$ in both, any point of $\mathcal A(\omega)\cap\Zo$ lies in $\overline{\mathcal A}\cap\Zo$. Now let $u_{0}\in\overline{\mathcal A}\cap\Zo$ and let $\bar\varphi\in[\pmin,\pmax]$ satisfy $\Phi(\bar\varphi)>\gamma$. By (i) applied to the constant control $\varphi\equiv\bar\varphi$, the forward orbit $\{S_{\bar\varphi}(t)u_{0}\,:\,t\ge0\}$ lies in $\overline{\mathcal A}$, and since this set is closed so does the $\omega$-limit set $\mathcal L_{\bar\varphi}(u_{0})$, which is nonempty, compact, connected and invariant under $S_{\bar\varphi}$ because the orbit is bounded by \cref{cor:absorbing}. The control versions of \cref{thm:preypersist,thm:uniformpersist}, applied to the constant control $\varphi\equiv\bar\varphi$, for which the essential bounds are both $\bar\varphi$ and $-\gamma+\Phi(\bar\varphi)>0$, give positive lower bounds for $\liminf_{t}x(t)$ and $\liminf_{t}y(t)$ along the orbit of $u_{0}$. Hence $\mathcal L_{\bar\varphi}(u_{0})\subset\Zo$. The Jacobian matrix of $\mathcal F(\cdot,\bar\varphi)$ at $v^{*}_{\bar\varphi}$ has the form
	\[
	\begin{pmatrix}
	\ast & -G(x^{*}_{\bar\varphi},\bar\varphi)\\[2pt]
	c\,y^{*}_{\bar\varphi}\,\partial_{x}G(x^{*}_{\bar\varphi},\bar\varphi) & 0
	\end{pmatrix},
	\]
whose determinant $c\,y^{*}_{\bar\varphi}G(x^{*}_{\bar\varphi},\bar\varphi)\,\partial_{x}G(x^{*}_{\bar\varphi},\bar\varphi)$ is positive, so $v^{*}_{\bar\varphi}$ is a node, a focus or a centre of the analytic planar vector field $\mathcal F(\cdot,\bar\varphi)$ \cite[Section 2.10, Theorems 4 and 5 and the Corollary to Theorem 5]{Perko2001}, and no orbit is homoclinic to it. Indeed, a neighbourhood of a centre is filled with periodic orbits, whereas a node or a focus, hyperbolic or not, is asymptotically stable either for $\mathcal F(\cdot,\bar\varphi)$ or for $-\mathcal F(\cdot,\bar\varphi)$. Reversing time if necessary, we may assume the former, and then $v^{*}_{\bar\varphi}$ is the $\alpha$-limit of no orbit $\Gamma\neq\{v^{*}_{\bar\varphi}\}$, because if $N$ is a neighbourhood of $v^{*}_{\bar\varphi}$ missing a point $\Gamma(t_{0})$ and $N'\subset N$ is a neighbourhood whose forward orbit stays in $N$, then $\Gamma(s)\in N'$ for some $s<t_{0}$, whence $\Gamma(t_{0})\in N$, a contradiction. Since $v^{*}_{\bar\varphi}$ is the only equilibrium in $\Zo$, the generalised Poincar\'e--Bendixson theorem \cite[Section 3.7, Theorem 2]{Perko2001} shows that $\mathcal L_{\bar\varphi}(u_{0})$ is either $\{v^{*}_{\bar\varphi}\}$ or a periodic orbit. If $\bar\varphi\in\Sigma_{\rm gs}$, then $S_{\bar\varphi}(t)u_{0}\to v^{*}_{\bar\varphi}$ and $\mathcal L_{\bar\varphi}(u_{0})=\{v^{*}_{\bar\varphi}\}$.

(iii) The segment $[0,k]\times\{0\}$ is contained in $\mathcal A(\omega)$ for every $\omega\in\Omega_{0}$, by the first paragraph of the proof of \cref{thm:attractorext}, which does not use \eqref{eq:E}, hence in $\overline{\mathcal A}$ by \cref{lem:essclosure}. Together with (ii) this gives $\mathcal E_{u_{0}}\subset\overline{\mathcal A}$, then $\mathcal R_{u_{0}}\subset\overline{\mathcal A}$ by (i), and $\overline{\mathcal R_{u_{0}}}\subset\overline{\mathcal A}$ because $\overline{\mathcal A}$ is closed. The second inclusion in \eqref{eq:lowerbound} is \cref{lem:essclosure}.

(iv) The inclusion $\supset$ follows from (i). Conversely, let $v\in\overline{\mathcal A}$ and $r>0$. The event $\{\dist(v,\mathcal A(\omega))<r\}$ has positive probability, whereas by \cref{lem:essclosure} and the invariance of $\PP$ the set of $\omega\in\Omega_{0}$ with $\mathcal A(\theta_{-t}\omega)\not\subset\overline{\mathcal A}$ is null. For $\omega$ in the former event and not in the latter set, $\mathcal A(\omega)=\Psi(t,\theta_{-t}\omega)\mathcal A(\theta_{-t}\omega)\subset\Psi(t,\theta_{-t}\omega)\overline{\mathcal A}$ contains a point at distance less than $r$ from $v$. The proof is complete.
\end{proof}

\Cref{thm:lowerbound} is the counterpart, when the predator persists, of \cref{thm:charext}, and it confines the possible difference between the two descriptions to the states of $\mathcal U$ that lie outside $\overline{\mathcal R_{u_{0}}}$. Whether such states exist is a question about the bounded complete solutions of the planar control system \eqref{eq:control} alone, since $\mathcal R_{u_{0}}$ and $\mathcal U$ are both defined in its terms once the sets $\mathcal L_{\bar\varphi}(u_{0})$ are known, whereas the question whether $\overline{\mathcal A}=\mathcal U$ still involves the noise through $\overline{\mathcal A}$. Properties (i) and (iv), the latter being the analogue of \eqref{eq:familyinv}, hold for all values of the parameters. What distinguishes $\overline{\mathcal A}$ from $\mathcal U$ is not their behaviour under the family of realisations, but the states they contain.

\section{The cost of fluctuation}\label{sec:cost}

We now compare the random model with the deterministic one having the same mean refuge, that is, we compare $\lambda$ with $\lambda_{\rm det}=-\gamma+\Phi(p)$. By \eqref{eq:meanpreserved} the two models have exactly the same mean accessible fraction, so any discrepancy between $\lambda$ and $\lambda_{\rm det}$ is attributable to the fluctuation alone. Recall from \eqref{eq:Vdef} the variance $V=\Var(\psi(z^{*}))$ of the accessible fraction, which satisfies $0<V<\varepsilon^{2}$.

\begin{Theorem}[The fluctuation is always detrimental to the predator]
\label{thm:cost} Assume \textbf{(H1)}--\textbf{(H3)}. Then
	\begin{equation}\label{eq:jensen}
		\lambda\ <\ \lambda_{\rm det},
	\end{equation}
and the deficit satisfies the two-sided bounds
	\begin{equation}\label{eq:cost}
		\frac{ac\beta k^{2}}{\bigl(1+ak\pmax\bigr)^{3}}\ V
		\ \ \le\ \ \lambda_{\rm det}-\lambda\ \ \le\ \
		\frac{ac\beta k^{2}}{\bigl(1+ak\pmin\bigr)^{3}}\ V .
	\end{equation}
\end{Theorem}

\begin{proof}
Write $X:=\psi(z^{*})$, a random variable with values in $(\pmin,\pmax)$, $\E[X]=p$ by \cref{lem:ergodic} and $\Var(X)=V>0$, the variance being positive because $\psi$ is strictly increasing and $z^{*}$ is a nondegenerate Gaussian variable. By \cref{lem:Fprops}(iii), $\Phi$ is of class $C^{2}$ and strictly concave on $[\pmin,\pmax]$. Define
	\[
	R(u):=\Phi(u)-\Phi(p)-\Phi'(p)(u-p),\qquad u\in[\pmin,\pmax],
	\]
a continuous function of $u$ alone, so that no measurable selection of an intermediate point is needed. Taylor's formula with integral remainder gives, for every $u\in[\pmin,\pmax]$,
	\[
	-R(u)=(u-p)^{2}\int_{0}^{1}(1-r)\,
	\bigl[-\Phi''\bigl(p+r(u-p)\bigr)\bigr]\,dr ,
	\]
the point $p+r(u-p)$ lying in $[\pmin,\pmax]$ for every $r\in[0,1]$. By \eqref{eq:Phider}, $-\tfrac12\Phi''(\vartheta)=ac\beta k^{2}/(1+ak\vartheta)^{3}$ is positive and decreasing in $\vartheta$, and $\int_{0}^{1}(1-r)\,dr=\tfrac12$, so that
	\[
	\frac{ac\beta k^{2}}{(1+ak\pmax)^{3}}\,(u-p)^{2}
	\ \le\ -R(u)\ \le\
	\frac{ac\beta k^{2}}{(1+ak\pmin)^{3}}\,(u-p)^{2},
	\qquad u\in[\pmin,\pmax].
	\]
Evaluating at $u=X$, taking expectations and using $\E[X-p]=0$, which gives
	\begin{equation}\label{eq:taylorexp}
		\lambda_{\rm det}-\lambda=\Phi(p)-\E[\Phi(X)]=\E\bigl[-R(X)\bigr],
	\end{equation}
we obtain \eqref{eq:cost}. The left-hand side of \eqref{eq:cost} is positive because $a,c,\beta,k>0$ and $V>0$, which proves \eqref{eq:jensen}.
\end{proof}

\Cref{thm:cost} is an instance of nonlinear averaging. The effect of the fluctuation on the invasion exponent is entirely mediated by the curvature of the functional response, because a Holling type~II response saturates, and an unusually large accessible fraction buys the predator less than an equally unusual small one costs it. Averaged over a zero-mean fluctuation, this asymmetry is a net loss whose size lies between $\tfrac12\min|\Phi''|\,V$ and $\tfrac12\max|\Phi''|\,V$. The content of \cref{thm:cost} is that for the invasion exponent of the predator this correction holds as the rigorous two-sided inequality \eqref{eq:cost} rather than as a second-order approximation, and that it governs the dynamics through \cref{thm:extinction,thm:meanpersist}. The next proposition shows that the bounds in \eqref{eq:cost} do not simply grow with the saturation of the response. They are largest at an intermediate degree of saturation and vanish at both ends.

\begin{proposition}[Where the bounds are largest]\label{prop:amax}
For fixed $\bar\varphi>0$ and fixed $c,\beta,k>0$, the map
	\begin{equation}\label{eq:kappaofa}
		a\ \longmapsto\ \varpi(a,\bar\varphi):=\frac{ac\beta k^{2}}{(1+ak\bar\varphi)^{3}}
		=-\tfrac12\Phi''(\bar\varphi)
	\end{equation}
is strictly increasing on $(0,\frac{1}{2k\bar\varphi}]$, strictly decreasing on $[\frac{1}{2k\bar\varphi},\infty)$, and tends to $0$ both as $a\downarrow0$ and as $a\to+\infty$, its maximum $\frac{4c\beta k}{27\bar\varphi}$ being attained at $ak\bar\varphi=\tfrac12$.
\end{proposition}

\begin{proof}
A direct differentiation gives $\partial_{a}\varpi(a,\bar\varphi)=c\beta k^{2}(1-2ak\bar\varphi)/(1+ak\bar\varphi)^{4}$, whose sign is that of $1-2ak\bar\varphi$, and $a^{2}\varpi(a,\bar\varphi)\to c\beta/(k\bar\varphi^{3})$ as $a\to+\infty$.
\end{proof}

When $a\downarrow0$ the response becomes linear and the effects of an unusually open and an unusually closed refuge cancel, whereas when $a\to+\infty$ the response is fully saturated and the refuge, variable or not, ceases to be felt. The location of the peak depends on what is held fixed as $a$ varies. With the attack rate $\beta$ fixed, as in \cref{prop:amax}, the bounds peak when the half-saturation density $1/(a\bar\varphi)$ equals twice the carrying capacity; if instead the maximal ingestion rate $\beta/a$ is held fixed, then $\varpi(a,\bar\varphi)=(\beta/a)\,ca^{2}k^{2}/(1+ak\bar\varphi)^{3}$ peaks at $ak\bar\varphi=2$, that is, when the half-saturation density equals half the carrying capacity.

We next show that the effect increases with the amplitude of the fluctuation. We use the following classical correlation inequality, whose short proof we include for completeness.

\begin{lemma}[Chebyshev's association inequality]\label{lem:cheby}
Let $X$ be a real random variable and let $f$ be nonincreasing and $g$ nondecreasing, both of them bounded and measurable. Then $\Cov(f(X),g(X))\le0$. If moreover $f$ is strictly decreasing, $g$ is strictly increasing and $\Var(X)>0$, then the inequality is strict.
\end{lemma}

\begin{proof}
Let $X'$ be an independent copy of $X$. For all real $u,v$ the product $(f(u)-f(v))(g(u)-g(v))$ is nonpositive, because $f$ and $g$ are monotone in opposite senses. Hence
	\[
	\Cov\bigl(f(X),g(X)\bigr)
	=\tfrac12\,\E\Bigl[\bigl(f(X)-f(X')\bigr)\bigl(g(X)-g(X')\bigr)\Bigr]\ \le\ 0 .
	\]
If $f$ and $g$ are strictly monotone, the integrand is negative whenever $X\ne X'$, an event of positive probability when $\Var(X)>0$.
\end{proof}

\begin{proposition}[Monotonicity in the amplitude]\label{prop:monotone}
Fix $\alpha,k,\beta,a,c,\gamma,m,\kappa,\sigma$ and regard $\lambda$ as a function of the amplitude $\varepsilon\in(0,\min\{m,1-m\})$. Then $\varepsilon\mapsto\lambda(\varepsilon)$ is of class $C^{1}$, strictly decreasing, and
	\begin{equation}\label{eq:dlambda}
		\frac{d\lambda}{d\varepsilon}(\varepsilon)
		=\Cov\Bigl(\Phi'\bigl(p+\varepsilon\,\Xi\bigr),\ \Xi\Bigr)\ <\ 0,
		\qquad \Xi:=\frac{2}{\pi}\arctan z^{*}\in(-1,1),\ \ \E[\Xi]=0 .
	\end{equation}
Moreover, $\lambda(\varepsilon)\to\lambda_{\rm det}$ as $\varepsilon\downarrow0$.
\end{proposition}

\begin{proof}
With $\Xi$ as in \eqref{eq:dlambda} we have $\psi(z^{*})=p+\varepsilon\,\Xi$ and $\lambda(\varepsilon)=-\gamma+\E[\Phi(p+\varepsilon\Xi)]$. The integrand is bounded by $c\beta/a$ and its $\varepsilon$-derivative $\Phi'(p+\varepsilon\Xi)\,\Xi$ is bounded by $c\beta k$ in absolute value, so differentiation under the expectation is legitimate and
	\[
	\frac{d\lambda}{d\varepsilon}=\E\bigl[\Phi'(p+\varepsilon\Xi)\,\Xi\bigr]
	=\Cov\bigl(\Phi'(p+\varepsilon\Xi),\Xi\bigr),
	\]
the last equality because $\E[\Xi]=0$. By \cref{lem:Fprops}(iii) the map $\Phi'$ is strictly decreasing, so $u\mapsto\Phi'(p+\varepsilon u)$ is strictly decreasing while $u\mapsto u$ is strictly increasing, and $\Xi$ is nondegenerate because $z^{*}$ is. Hence \cref{lem:cheby} gives \eqref{eq:dlambda} with strict inequality. The derivative is continuous in $\varepsilon$ by dominated convergence, with the same majorant, which gives the $C^{1}$ regularity, and $\lambda(\varepsilon)\to-\gamma+\Phi(p)=\lambda_{\rm det}$ as $\varepsilon\downarrow0$ by the same argument.
\end{proof}

The following consequence shows that a bounded zero-mean fluctuation of the refuge can drive to extinction a predator that would persist in the deterministic model, in a regime in which no worst-case argument detects the phenomenon. It should be distinguished from the extinction studied in \cite{BenaimLobry2016}, induced by switching between environments each of which would sustain the species, where the resident state moves with the environment and the invasion rate averages the invader's growth rate over a nontrivial boundary distribution. Here the prey-only state is the same in every environment, the invasion exponent is the plain average of the exponents $-\gamma+\Phi(\bar\varphi)$ of the frozen systems, and by \eqref{eq:lambdaorder} no switching can drive to extinction a predator that every frozen environment sustains; the reversal below is produced by the curvature of the response alone.

\begin{corollary}[Noise-induced extinction]\label{cor:noiseinduced}
Assume \textbf{(H1)}--\textbf{(H3)} and
	\begin{equation}\label{eq:nie}
		0\ <\ \lambda_{\rm det}\ <\ \frac{ac\beta k^{2}}{(1+ak\pmax)^{3}}\ V .
	\end{equation}
Then $\lambda<0$, hence by \cref{thm:extinction} the predator becomes extinct exponentially fast and $x(t)\to k$, although in the deterministic model \eqref{eq:det} with the same mean refuge $m$ the predator invades and a unique interior equilibrium \eqref{eq:detequil} exists. Moreover, in this regime condition \eqref{eq:Eprime} necessarily fails.
\end{corollary}

\begin{proof}
By the left inequality in \eqref{eq:cost} and \eqref{eq:nie}, $\lambda\le\lambda_{\rm det}-ac\beta k^{2}V/(1+ak\pmax)^{3}<0$. That the deterministic model has an interior equilibrium is the statement $\lambda_{\rm det}>0$, discussed after \eqref{eq:detequil}. Finally, $\Phi(\pmax)>\Phi(p)>\gamma$ by strict monotonicity of $\Phi$, so \eqref{eq:Eprime} fails.
\end{proof}

\subsection{The reversal window}\label{ssec:reversal}

\Cref{cor:noiseinduced} gives a closed-form sufficient condition for the fluctuation to reverse the outcome of the interaction. In fact the set of parameters on which this happens can be described exactly, and its size is precisely the deficit of \cref{thm:cost}.

Fix $\alpha,k,\beta,a,c,m,\varepsilon,\kappa,\sigma$ and regard $\gamma$ as free. By \cref{def:lambda}, the predator invades in the deterministic model \eqref{eq:det} and becomes extinct in the random one, that is, $\lambda_{\rm det}>0>\lambda$, if and only if
\begin{equation}\label{eq:window}
	\gamma\ \in\ \mathcal W:=\Bigl(\,\E\bigl[\Phi(\psi(z^{*}))\bigr],\ \Phi(p)\,\Bigr),
\end{equation}
the dynamical translation being \cref{thm:extinction} on one side and the classical analysis of \eqref{eq:det} on the other. The length of this open interval is exactly $\lambda_{\rm det}-\lambda$, so that by \eqref{eq:cost}
\begin{equation}\label{eq:windowlength}
	\frac{ac\beta k^{2}}{(1+ak\pmax)^{3}}\,V
	\ \le\ |\mathcal W|\ \le\
	\frac{ac\beta k^{2}}{(1+ak\pmin)^{3}}\,V ,
\end{equation}
and $\mathcal W$ is nonempty for every choice of the remaining parameters.

A reversal is therefore always possible, but it requires a predator whose deterministic exponent is already small. The lower bound in \eqref{eq:windowlength} is the informative one when $ak\pmax$ is large, since the ratio of the two bounds is $((1+ak\pmax)/(1+ak\pmin))^{3}$. This ratio tends to $1$ as $\varepsilon\downarrow0$, because $\pmin$ and $\pmax$ both tend to $p$, so that for narrow bands the two bounds pin the deficit down to $-\tfrac12\Phi''(p)\,V$ up to a factor arbitrarily close to $1$. At the parameters of \cref{ex:nie} the window is $\mathcal W=(0.55963,\,0.56842)$ and the value $\gamma=0.5679522$ chosen there lies inside it.

\subsection{Dependence on the driving process}\label{ssec:kapsig}

The parameters of the Ornstein--Uhlenbeck process enter $\lambda$ only through the stationary standard deviation $\sigma_{*}=\sigma/\sqrt{2\kappa}$, by \eqref{eq:lambdaintegral}. Increasing $\sigma_{*}$ spreads the law of $z^{*}$, and the next proposition shows that $\lambda$ decreases accordingly. The argument of \cref{prop:monotone} does not transfer verbatim, because $z\mapsto\frac{\partial}{\partial\sigma_{*}}\arctan(\sigma_{*}z)$ is not monotone in $z$, and the association inequality is replaced by a symmetrisation. Write $V(\sigma_{*})$ and $\lambda(\sigma_{*})$ for the quantities \eqref{eq:Vdef} and \eqref{eq:lambda} regarded as functions of $\sigma_{*}\in(0,\infty)$, all the other parameters being fixed.

\begin{proposition}[Monotonicity in the stationary variance]\label{prop:kapsig}
The maps $\sigma_{*}\mapsto V(\sigma_{*})$ and $\sigma_{*}\mapsto-\lambda(\sigma_{*})$ are of class $C^{1}$ and strictly increasing on $(0,\infty)$. More precisely, if $\mathsf N$ denotes a standard Gaussian random variable,
	\begin{equation}\label{eq:dlambdasigma}
		\frac{d\lambda}{d\sigma_{*}}(\sigma_{*})
		=\frac{2\varepsilon}{\pi}\,
		\E\Bigl[\Phi'\Bigl(p+\tfrac{2\varepsilon}{\pi}\arctan(\sigma_{*}\mathsf N)\Bigr)
		\frac{\mathsf N}{1+\sigma_{*}^{2}\mathsf N^{2}}\Bigr]\ <\ 0 .
	\end{equation}
Moreover,
	\begin{equation}\label{eq:lambdasigmalimits}
		\lim_{\sigma_{*}\downarrow0}\lambda(\sigma_{*})=\lambda_{\rm det},
		\qquad
		\lim_{\sigma_{*}\to+\infty}\lambda(\sigma_{*})
		=-\gamma+\tfrac12\bigl[\Phi(\pmin)+\Phi(\pmax)\bigr].
	\end{equation}
\end{proposition}

\begin{proof}
Since $z^{*}$ has the law of $\sigma_{*}\mathsf N$, we have $V(\sigma_{*})=\frac{4\varepsilon^{2}}{\pi^{2}}\E[\arctan^{2}(\sigma_{*}\mathsf N)]$ and $\lambda(\sigma_{*})=-\gamma+\E[\Phi(p+\frac{2\varepsilon}{\pi}\arctan(\sigma_{*}\mathsf N))]$. For each fixed $n\neq0$ the map $\sigma_{*}\mapsto\arctan^{2}(\sigma_{*}n)$ is strictly increasing, so $V$ is strictly increasing, and its $\sigma_{*}$-derivative $\frac{4\varepsilon^{2}}{\pi^{2}}\E[2\arctan(\sigma_{*}\mathsf N)\,\mathsf N/(1+\sigma_{*}^{2}\mathsf N^{2})]$ has the integrable majorant $\frac{4\varepsilon^{2}}{\pi}|\mathsf N|$, independent of $\sigma_{*}$, so $V$ is of class $C^{1}$ by dominated convergence. The $\sigma_{*}$-derivative of the integrand defining $\lambda$ is
	\[
	\frac{2\varepsilon}{\pi}\,\Phi'\Bigl(p+\tfrac{2\varepsilon}{\pi}\arctan(\sigma_{*}\mathsf N)\Bigr)
	\frac{\mathsf N}{1+\sigma_{*}^{2}\mathsf N^{2}},
	\]
which is bounded in absolute value by $\frac{2\varepsilon}{\pi}c\beta k|\mathsf N|$, an integrable majorant independent of $\sigma_{*}$. Hence differentiation under the expectation is legitimate and gives the identity in \eqref{eq:dlambdasigma}. To see that the expectation is negative, put $h(n):=\Phi'(p+\frac{2\varepsilon}{\pi}\arctan(\sigma_{*}n))$ and $g(n):=n/(1+\sigma_{*}^{2}n^{2})$. The function $h$ is strictly decreasing, because $\Phi'$ is strictly decreasing by \eqref{eq:Phider} and $\arctan$ is strictly increasing, while $g$ is odd and positive on $(0,\infty)$. Splitting the integral according to the sign of $\mathsf N$ and using the symmetry of its law,
	\[
	\E\bigl[h(\mathsf N)g(\mathsf N)\bigr]
	=\frac{1}{\sqrt{2\pi}}\int_{0}^{\infty}\bigl[h(n)-h(-n)\bigr]\,g(n)\,e^{-n^{2}/2}\,dn\ <\ 0 ,
	\]
since the integrand is negative for every $n>0$. The derivative is continuous in $\sigma_{*}$ by dominated convergence, with the same majorant. Finally, as $\sigma_{*}\downarrow0$ the integrand defining $\lambda$ converges pointwise to $\Phi(p)$, and as $\sigma_{*}\to+\infty$ it converges pointwise to $\Phi(\pmax)$ on $\{\mathsf N>0\}$ and to $\Phi(\pmin)$ on $\{\mathsf N<0\}$. Both limits in \eqref{eq:lambdasigmalimits} follow from dominated convergence, the integrand being bounded by $c\beta/a$.
\end{proof}

The second limit in \eqref{eq:lambdasigmalimits} is the exponent of a refuge that switches between its two extreme values with equal frequencies; it is the infimum of $\lambda$ for a given amplitude, and it is not attained. Together with \cref{prop:monotone}, \cref{prop:kapsig} says that the deficit $\lambda_{\rm det}-\lambda$ is strictly increasing in each of the two parameters that measure the size of the fluctuation, and the range of $\lambda$ is the open interval
\[
\Bigl(-\gamma+\tfrac12\bigl[\Phi(\pmin)+\Phi(\pmax)\bigr],\ \lambda_{\rm det}\Bigr).
\]
Since $\varepsilon$ and $\sigma_{*}$ enter $\lambda$ only through the law of $\psi(z^{*})$, the threshold cannot distinguish a wide band with a concentrated driving process from a narrow band with a spread one, whereas only the band itself is visible to the worst-case exponents.

\subsection{Ecological reading}\label{ssec:ecology}

Fix the mean refuge level $m$ and let the actual refuge fluctuate around it within a band of half-width $\varepsilon$. Then the predator is strictly worse off, as measured by its invasion exponent, than under the constant refuge $m$ (\cref{thm:cost}), and increasingly so as the band widens or the environment becomes more variable (\cref{prop:monotone,prop:kapsig}), while the correlation time of the environment is irrelevant to the threshold, by \eqref{eq:lambdaintegral}, once the stationary variance of the driving process is fixed. This feature is specific to the structure of the problem. The invasion exponent at the prey-only state is a time average of a function of the environment alone, because the prey-only dynamics are unaffected by the noise, and time averages of a stationary process are blind to its autocorrelation, whereas the colour of environmental noise is known to matter when the environment acts on a population whose own dynamics are nonlinear \cite{RipaLundberg1996,PetcheyGonzalezWilson1997,RuokolainenLindenKaitalaFowler2009}. The correlation time does govern the speed at which the averages of \cref{lem:ergodic} converge, hence the horizon over which the sign of $\lambda$ becomes visible, and near the threshold, as in \cref{ex:nie}, this is the ecologically relevant time scale.

In the treatments of environmental stochasticity in which the noise enters through the growth rates \cite{May1973,Turelli1977}, the effect is driven by the variance of the per-capita increments, whereas here the per-capita rates are bounded, the mean environment is held exactly fixed, and the invasion exponent is unaffected in the linear limit. \Cref{prop:amax}, read with the maximal ingestion rate $\beta/a$ held fixed as discussed after it, yields a comparative prediction. Of three predator species with the same maximal ingestion rate but differing in the saturation level $akp$ of the accessible prey at carrying capacity, the one operating in the linear part of its response curve and the one deep in saturation should both be comparatively insensitive to the variability of the refuge, while the one operating in the transition region should be markedly sensitive to it. The effect is quantitatively small unless the band is wide, since the deficit is of second order in $\varepsilon$, and \cref{ex:nie} accordingly uses a band covering most of the unit interval and a predator whose deterministic exponent is below $10^{-3}$.

\section{Open questions}\label{sec:open}

By \cref{thm:lowerbound}, when the predator persists the forward and the pullback descriptions coincide, in the sense that $\overline{\mathcal A}=\mathcal U$, as soon as the value at time zero of every bounded complete solution of the planar control system \eqref{eq:control} lies in the closure of the states reachable from the segment and from the frozen limit sets of that theorem. Whether this holds, and whether the two descriptions can differ at all in this regime, we leave open, together with the related question of whether $\mathcal U$ coincides with the closure of the union of the frozen attractors. The critical case $\lambda_{+}=0$ is covered by \cref{thm:charext}, whereas the critical cases $\lambda=0$ and $\lambda_{-}=0$ are not treated, the former being a central-limit rather than a law-of-large-numbers regime in which the correlation time of the environment will matter. The structure of the pullback attractor in the regime $\lambda_{-}<0<\lambda$, where the predator persists in the mean but not uniformly, remains to be described; the occupation-measure theory of \cite{SchreiberBenaimAtchade2011,BenaimSchreiber2019,HeningNguyenChesson2021} is the natural tool there.

\section*{Acknowledgements}
\noindent This work is part of the doctoral research of E.~F. Cruzado-Quispe, carried out within the Doctoral Programme in Mathematics of the Universidad Nacional del Santa under the supervision of P.~N. Seminario-Huertas and J.~L\'opez-de-la-Cruz.

\section*{Statements and Declarations}
\begin{sloppypar}
\noindent\textbf{Funding.} J.~L\'opez-de-la-Cruz and P.~N. Seminario-Huertas have been supported by the Madrid Government (Comunidad de Madrid, Spain) under the Multiannual Agreement 2023--2026 with the Universidad Polit\'ecnica de Madrid, Line~A, Emerging PhD Researchers (Grant agreement DOCTORES-EMERGENTES-24-UMDLRU-108-XVLWPC). J.~L\'opez-de-la-Cruz has also been supported by the project PID2024-156228NB-I00, funded by \mbox{MICIU/AEI}/\allowbreak 10.13039/\allowbreak 501100011033 and by ERDF/EU (Ministerio de Ciencia, Innovaci\'on y Universidades, Agencia Estatal de Investigaci\'on and European Regional Development Fund, Spain), and P.~N. Seminario-Huertas by the project PID2025-174633NA-I00, funded by \mbox{MICIU/AEI}/\allowbreak 10.13039/\allowbreak 501100011033 (Ministerio de Ciencia, Innovaci\'on y Universidades and Agencia Estatal de Investigaci\'on, Spain). The funding sources had no involvement in the study design, in the analysis, in the writing of the article or in the decision to submit it for publication.
\end{sloppypar}

\smallskip
\noindent\textbf{Competing interests.} The authors have no competing interests to declare that are relevant to the content of this article.

\smallskip
\noindent\textbf{Data availability.} Data sharing is not applicable to this article, as no datasets were generated or analysed during the current study.

\end{document}